\documentclass[leqno,12pt]{amsart}

\usepackage{amscd}              
\usepackage[all]{xy}            
\xyoption{curve}

\CompileMatrices

\makeatletter
\def\overunderbraces #1#2#3{{%
 \baselineskip\z@skip \lineskip4\p@ \lineskiplimit4\p@
 \displaystyle  
 \setbox\z@\vbox{\ialign{&\hfil${}##{}$\hfil\cr
   \global\let\br\br@label #1\cr 
   \global\let\br\br@down #1\cr   
   #2\cr 
 }}
 \dimen@-\ht\z@ 
 \setbox\z@\vbox{\ialign{&\hfil${}##{}$\hfil\cr
   \global\let\br\br@label #1\cr 
   \global\let\br\br@down #1\cr   
   #2\cr 
   \global\let\br\br@up #3\cr 
   \global\let\br\br@label #3\cr   
 }}
 \advance\dimen@\ht\z@ 
 \lower\dimen@\hbox{\box\z@} 
}}
\def\br@up#1#2{\multispan{#1}\upbracefill}
\def\br@down#1#2{\multispan{#1}\downbracefill}
\def\br@label#1#2{\multispan{#1}\hidewidth $#2$\hidewidth}
\makeatother

\newcommand{\C}{{\mathbb{C}}}                    
\newcommand{\card}{\sharp}
\newcommand{\coker}{\text{coker}}
\newcommand{\comp}{\circ}
                  
\newcommand{\deltabf}{\text{\boldmath $\delta$}}
\newcommand{\dual}{\vee}
\newcommand{\eins}{{\mathbb{I}}}
\newcommand{\E}{{\mathcal{E}}}
\newcommand{\F}{{\mathcal{F}}}
\newcommand{\Fl}{\text{Fl}}
\newcommand{\Gr}{\text{Gr}}
\newcommand{\Gl}{\text{Gl}}
\newcommand{\Gln}{\text{Gl}_n}
\newcommand{\Hh}{{\mathcal{H}}}
\newcommand{\id}{\text{id}}
\newcommand{\im}{\text{im}}
\newcommand{\injto}{\hookrightarrow}
\newcommand{\Isom}{\text{Isom}}
\newcommand{\isomorph}{\cong}           
\newcommand{\isomto}{\overset{\sim}{\rightarrow}}  
\newcommand{\I}{{\mathcal{I}}}
\newcommand{\Isomto}{\overset{\sim}{\longrightarrow}}
\newcommand{\J}{{\mathcal{J}}}
\newcommand{\K}{{\mathcal{K}}}
\newcommand{\KGl}{\text{KGl}}
\newcommand{\KGlbf}{\text{\bf KGl}}
\newcommand{\KGln}{\text{KGl}_n}
\newcommand{\Ll}{{\mathcal{L}}}                  
\newcommand{\m}{{\mathfrak{m}}}                  
\newcommand{\M}{{\mathcal{M}}}
\newcommand{\Nn}{{\mathcal{N}}}
\newcommand{\N}{{\mathbb{N}}}                    
\newcommand{\Nc}{{\mathcal{N}}}
\newcommand{\ot}{\leftarrow}
\newcommand{\Oo}{{\mathcal{O}}}                  
\newcommand{\Pc}{{\mathcal{P}}}
\newcommand{\PGl}{\text{PGl}}
\newcommand{\PGlb}{\overline{\text{PGl}}}
\newcommand{\PGlbf}{\overline{\text{\bf PGl}}}
\newcommand{\PGln}{\text{PGl}_n}
\newcommand{\PGlnb}{\overline{\text{PGl}}_n}
\newcommand{\PGlnC}{\text{PGl}_{n,\C}}
\newcommand{\Proj}{\text{Proj}\,}
\newcommand{\Q}{{\mathbb{Q}}}                    
\newcommand{\rk}{\text{rk}\, }   
\newcommand{\Spec}{\text{Spec}\, }               
\newcommand{\tensor}{\otimes}
\newcommand{\T}{\widetilde{T}}
\newcommand{\Tensor}{\bigotimes}
\newcommand{\To}{\longrightarrow}
\newcommand{\ub}{{\text{{\boldmath $u$}}}}

\newcommand{\Wedge}{\bigwedge} 
\newcommand{\x}{{\text{{\boldmath $x$}}}}    

\newcommand{\y}{{\text{{\boldmath $y$}}}}

\newcommand{\z}{{\text{{\boldmath $z$}}}}
\newcommand{\Z}{{\mathbb{Z}}}                    
\newcommand{\one}{{\text{{\boldmath $1$}}}}

\newtheorem{theorem}{Theorem}[section]
\newtheorem{proposition}[theorem]{Proposition}
\newtheorem{lemma}[theorem]{Lemma}
\newtheorem{corollary}[theorem]{Corollary}

\theoremstyle{definition}
\newtheorem{definition}[theorem]{Definition}

\begin{document}

\title[A Compactification of $\Gl_n$]
{A Modular Compactification of the General Linear Group}
\author[Ivan Kausz]{Ivan Kausz}
\date{June 2, 2000}
\address{NWF I - Mathematik, Universit\"{a}t Regensburg, 93040 Regensburg, 
Germany}
\email{ivan.kausz@mathematik.uni-regensburg.de}

\maketitle

\tableofcontents

\section{Introduction}

In this paper we give a modular description of a certain
compactification $\KGln$ of the general linear group $\Gln$. 
The variety $\KGln$ is
constructed as follows: First one embeds $\Gln$ in the obvious way in 
the projective space which contains the affine space of $n\times n$ matrices 
as a standard  open set. Then one successively blows up the closed
subschemes defined by the vanishing of the $r\times r$ subminors
($1\leq r\leq n$),
along with the intersection of these subschemes with the hyperplane 
at infinity. 

We were led to the problem of finding a modular description of $\KGln$ in
the course of our research on the degeneration of moduli spaces of 
vector bundles. Let me explain in some detail 
the relevance of compactifications of 
$\Gln$ in this context. 

Let $B$ be a regular integral one-dimensional base scheme and  $b_0\in B$
a closed point.
Let $C\to B$ be a proper flat familly of curves over $B$ which is smooth
outside $b_0$ and whose fibre $C_0$ at $b_0$ is irreducible with one 
ordinary double point $p_0\in C_0$.
Let $\tilde{C_0}\to C_0$ be the normalization of $C_0$ and let 
$p_1,p_2\in\tilde{C_0}$ the two points lying above the singular point $p_0$.
Thus the situation may be depicted as follows:
$$
\xymatrix@R-2ex{
&& (\text{$\tilde{C_0}$},p_1,p_2) \ar[dll] \ar[drr] && && \\ 
\text{$\tilde{C_0}$}
&& && C_0 \ar[d] \ar[rr] && C \ar[d] \\
&& && b_0 \ar[rr]        && B
}
$$
where the left arrow means ``forgetting the points $p_1,p_2$''.
There is a corresponding diagram of moduli-functors of vector bundles
(v.b.) of rank $n$:

$$
\xymatrix@R-3ex{
&
\hbox{\footnotesize
$ 
\left\{
\begin{array}{lll}
\text{\qquad v.b. $E$ on $\tilde{C_0}$} \\
\text{\quad \ together with an} \\
\text{isomorphism $E[p_1]\isomto E[p_2]$}
\end{array}
\right\}
$
}
\ar[dl]^{f_1} \ar[dr]_{f_2}^(.6){\isomorph} & & \\ 
\hbox{\footnotesize
$ 
\left\{
\begin{array}{ll}
\text{v.b. on}\\
\quad \tilde{C_0}
\end{array}
\right\}
$
}
& & 
\hbox{\footnotesize
$ 
\left\{
\begin{array}{ll}
\text{v.b. on}\\
\quad C_0
\end{array}
\right\}
$
}
\ar[d] \ar[r] & 
\hbox{\footnotesize
$ 
\left\{
\begin{array}{ll}
\text{v.b. on}\\
\ C/B
\end{array}
\right\}
$
}
\ar[d] \\
& & b_0 \ar[r]        & B
}
$$
where $E[p_i]$ denotes the fibre of $E$ at the point $p_i$ 
(cf. section \ref{notation} below). The morphism $f_1$ is ``forgetting
the isomorphism between the fibres'' and $f_2$ is 
``glueing together the fibres at $p_1$ and $p_2$ along the
given isomorphism''. The square on the right is the inclusion
of the special fibre. 
It is clear that $f$ is a locally trivial fibration
with fibre $\Gln$. Consequently, $f_1$ is not proper and thus
$\{\text{v.b. on $C/B$}\}$ is not proper over $B$. It is
desirable to have a diagram $(*)$:

$$
\xymatrix@R-2ex{
&
\hbox{\footnotesize
$ 
\left\{
\begin{array}{lll}
\text{\ generalized} \\
\text{v.b.-data on} \\
(\tilde{C_0},p_1,p_2)
\end{array}
\right\}
$
}
\ar[dl]^{\overline{f_1}} \ar[dr]_{\overline{f_2}} & & \\ 
\hbox{\footnotesize
$ 
\left\{
\begin{array}{ll}
\text{v.b. on}\\
\quad \tilde{C_0}
\end{array}
\right\}
$
}
& & 
\hbox{\footnotesize
$ 
\left\{
\begin{array}{lll}
\text{generalized}\\
\text{\quad v.b. on}\\
\quad \ \ C_0
\end{array}
\right\}
$
}
\ar[d] \ar[r] & 
\hbox{\footnotesize
$ 
\left\{
\begin{array}{lll}
\text{generalized}\\
\text{\quad v.b. on}\\
\quad \ C/B
\end{array}
\right\}
$
}
\ar[d] \\
& & b_0 \ar[r]        & B
}
$$
where the functors of ``generalized'' objects contain
the original ones as open subfunctors and where
$\{\text{generalized v.b. on $C/B$}\}$ 
is proper over $B$ or at least satisfies the existence
part of the valuative criterion for properness.
The motivation is that such a diagram may help to
calculate cohomological invariants of 
$\{\text{v.b. on $Y$}\}$
($Y$ a smooth projective curve) by induction on the genus of
$Y$ (notice that the genus of $\tilde{C_0}$ is one less
than the genus of the generic fibre of $C/B$).

In the current literature there exist two different approaches
for the construction of diagram $(*)$. In the first approach
the ``generalized v.b.'' on $C_0$  are torsion-free sheaves
(cf. \cite{Seshadri1}, \cite{Faltings}, \cite{Narasimhan&Ramadas},
\cite{Sun}). The second approach is by Gieseker who considered
 only the rank-two case (cf. \cite{Gieseker}). 
Here the ``generalized v.b.'' on $C_0$ are certain vector bundles on
$C_0$, $C_1$ or $C_2$, where $C_i$ is built from $C_0$ by inserting a 
chain of $i$ copies of the projective line at $p_0$.
(Cf. also \cite{Teixidor i Bigas} for a discussion of the two 
approaches).
Of course, this is only a very rough picture of what is going on
in these papers since I do not mention concepts of stability
for the various objects nor
the representability of the functors by varieties or by
algebraic stacks.

In both approaches the morphism $\overline{f_2}$ is the normalization
morphism (at least on the complement of a set of small dimension)
and $\overline{f_1}$ is a locally trivial fibration with fibre a 
compactification of $\Gln$.
In the torsion-free sheaves approach this compactification is
Gr$(2n,n)$, the grassmanian of $n$-dimensional subspaces of a
$2n$-dimensional vector space. In Gieseker's construction the relevant 
compactification of $\text{Gl}_2$ is $\text{KGl}_2$.
An advantage of Gieseker's construction is that in contrast to the 
torsion-free sheaves approach, the space 
$\{\text{generalized v.b. on $C/B$}\}$
is regular and its special fibre over $b_0$ is a divisor with normal 
crossings. 

Very recently, Nagaraj and Seshadri have generalized  Gieseker's
construction of the right part of diagram $(*)$, i.e. the diagram
$$
\xymatrix@R-2ex{
\hbox{
$ 
\left\{
\begin{array}{lll}
\text{generalized}\\
\text{\quad v.b. on}\\
\quad \ \ C_0
\end{array}
\right\}
$
}
\ar[d] \ar[r] & 
\hbox{
$ 
\left\{
\begin{array}{lll}
\text{generalized}\\
\text{\quad v.b. on}\\
\quad \ C/B
\end{array}
\right\}
$
}
\ar[d] \\
b_0 \ar[r]        & B
}
$$
to arbitrary rank $n$ (cf. \cite{Nagaraj&Seshadri}, \cite{Seshadri2}).
Nagaraj's and Seshadri's ``generalized vector bundles'' on $C_0$
are certain equivalence classes of vector bundles on one of the 
curves $C_0,\dots,C_n$, whose push-forward to $C_0$ are stable
torsion free sheaves.

Without worrying about stability I have recently (and independently
from Nagaraj and Seshadri) constructed the full diagram $(*)$ at
least at the level of functors (details will appear in a forthcoming 
paper) and I have reasons to believe 
that the fibres of 
the corresponding morphism $\overline{f_1}$ should be represented by
$\KGln$. The present paper is the first step in the proof of this fact.

The compactification $\KGln$ of $\Gln$ has properties similar
to those of the ``wonderful compactification'' of algebraic groups
of adjoint type as studied by De Concini and Procesi 
(cf. \cite{CP}). Namely:
\begin{enumerate}
\item
The group $\Gln\times\Gln$ acts on $\KGln$, extending the operation
 of $\Gln\times\Gln$ on $\Gln$
induced by right and left multiplication (cf. \ref{operation}).
\item
The complement of $\Gln$ in $\KGln$ is a divisor with normal crossings 
with irreducible components $Y_i$, $Z_j$ ($i,j\in\{0,\dots,n-1\}$)
(cf. \ref{nc}).
\item
The orbit closures of the operation of  $\Gln\times\Gln$ 
on $\KGln$ are precisely the intersections $Y_I\cap Z_J$, where $I,J$
are subsets of $\{0,\dots,n-1\}$ with $\min(I)+\min(J)\geq n$ and where 
$Y_I:=\cap_{i\in I} Y_i$, 
$Z_J:=\cap_{j\in J} Z_j$
(cf. \ref{orbits}).
\item
For each $I,J$ as above there exists a natural mapping from
$Y_I\cap Z_J$ to the product of two flag varieties. This mapping
is a locally trivial fibration with standard fibre a product
of copies of $\overline{\text{PGl}}_{n_k}$ (the wonderful compactification
of $\text{PGl}_{n_k}$) for some $n_k\geq 1$ and of one copy of
$\text{KGl}_m$ for some $m\geq 0$
(cf. \ref{strata}).
\end{enumerate}

Our main theorem \ref{KGln modular} says that $\KGln$ parametrizes what 
we call ``generalized isomorphisms'' from the trivial bundle of 
rank $n$ to itself. A generalized isomorphism between vector
bundles $E$ and $F$ is by definition a diagram
$$
\hbox{\footnotesize
\xymatrix@C=1ex{
E & = & E_0
\ar@/^1.0pc/|{\tensor}[rr]
& &
E_1 
\ar[ll]
\ar@/^1.0pc/|{\tensor}[rr]
& &
E_2
\ar[ll]
& 
\dots
& 
E_{n-1}
\ar@/^1.0pc/|{\tensor}[rrr]
& & &
E_n
\ar[lll]
\ar[rr]^\sim
& &
F_n
\ar[rrr]
& & &
F_{n-1}
\ar@/_1.0pc/|{\tensor}[lll]
& 
\dots
&
F_2
\ar[rr]
& &
F_1
\ar[rr]
\ar@/_1.0pc/|{\tensor}[ll]
& &
F_0
\ar@/_1.0pc/|{\tensor}[ll]
& = & F
}}
$$
with certain properties, where the $E_i$ and $F_j$ are vector
bundles of the same rank as $E$ and $F$ and where the arrow
$
\xymatrix{
\ar[r]|{\tensor} &
}
$
indicates a morphims of the source into the target tensored with a
line bundle to be specified. Cf. \ref{generalized isomorphism}
for a precise definition.

The wonderful compactification $\PGlnb$ of $\PGln$ 
is contained as an orbit closure in $\KGln$, in fact
$Y_0\isomorph\PGlnb$. Therefore theorem
\ref{KGln modular} implies a modular description of $\PGlnb$.
One of the reasons why I decided to publish the present paper
separately from my investigations on the degeneration of moduli spaces
of vector bundles on curves is the fact that $\PGlnb$ has been
quite extensively studied in the past (cf. \cite{Laksov1} for
a historical overview and also the recent paper \cite{Thaddeus2}).
Although some efford has been made 
to find a modular description for it, up to now only partial
results in this direction have been obtained  (cf. \cite{Vainsencher},
\cite{Laksov2}, \cite{Thorup&Kleiman}). In section \ref{PGln} we explain
the connection of these results with ours. Recently Lafforgue has used  
$\PGlnb$ to compactify the stack of Drinfeld's 
shtukas (cf. \cite{Lafforgue1}, \cite{Lafforgue2}).

Sections \ref{construction} and \ref{bf&gi} contain the main
definitions:
In section \ref{construction} we give the construction of $\KGln$ and
in section \ref{bf&gi} we define the notion of generalized isomorphisms.
At the end of section \ref{bf&gi} we state our main theorem 
\ref{KGln modular}. 
Its proof is given in sections \ref{exterior} and \ref{proof}. In section
\ref{PGln} we define complete collineations and compare our
notion with the one given by previous authors, in section 
\ref{geometry of the strata}
we study the orbit closures of the operation of $\Gln\times\Gln$
on $\KGln$ and in section \ref{grassmann} we define an equivariant 
morphism of $\KGln$ onto the Grassmannian compactification of $\Gln$
and compute its fibres.

My interest in degeneration of moduli spaces of bundles on curves 
has been greatly stimulated by a workshop on conformal blocks and the Verlinde
formula, organized  in March 1997 by  the physicists J\"{u}rgen Fuchs 
and Christoph Schweigert 
at the Mathematisches Forschungsinstitut in Oberwolfach.
Part of this work has been prepared during a stay at the Mathematical 
Institute of the University of Oxford. Its hospitality is gratefully
acknowledged. 
Thanks are due to Daniel Huybrechts for mentioning to me the work of
Thaddeus, to M. Thaddeus himself for sending me a copy of part of his
thesis and to M. Rapoport for drawing my attention to the work of Laksov
and  Lafforgue. I would also like to thank Uwe Jannsen for his
constant encouragement.

\section{An elementary example}

This section is not strictly necessary for the comprehension
of what follows. But since the rest of the paper is a bit
technical, I felt that a simple example might facilitate
its understanding.

Let $A$ be a discrete valuation ring, $K$ its field of fractions, 
$\m$ its maximal ideal, $t\in\m$ a local parameter and $k:=A/\m$ the
residue class field of A.
Let $E$ and $F$ be two free $A$-modules of rank $n$ and let 
$\varphi_K: E_K\isomto F_K$
be an isomorphism between the generic fibers $E_K:=E\tensor_A K$ and
$F_K:=F\tensor_A K$ of $E$ and $F$.
We can choose $A$-bases of $E$ and $F$ such that $\varphi_K$
has the matrix presentation
$\text{diag}(t^{m_1},\dots,t^{m_n})$
with respect to these bases,
where $m_i\in\Z$ and $m_1\leq\dots\leq m_n$.
Now let $a_0:=0=:b_0$ and for $1\leq i\leq n$ set
$a_i:=-\min(0,m_{n+1-i})$ and $b_i:=\max(0,m_i)$.
Note that we have
\begin{eqnarray*}
& & 
0=a_0=\dots=a_{n-l}\leq a_{n-l+1}\leq\dots\leq a_n
\\ 
\text{and} & &
0=b_0=\dots=b_l\leq b_{l+1}\leq\dots\leq b_n
\end{eqnarray*}
for some $l\in\{0,\dots,n\}$.
Let 
$$
E_n\subseteq\dots\subseteq E_1\subseteq E_0:=E 
\quad \text{and}\quad 
F_n\subseteq\dots\subseteq F_1\subseteq F_0:=F
$$
be the $A$-submodules
defined by
$$
E_{i+1}:=
\left[
\vcenter{
\xymatrix @R-4pc@ C-4pc @! {
t^{a_{i+1}-a_i}\eins_{n-i} & 0 \\
0 & \eins_i
}} \right] E_i
\quad ,
\quad
F_{i+1}:=
\left[
\vcenter{
\xymatrix @R-3.7pc@ C-3.7pc @! {
\eins_i & 0 \\
0 & t^{b_{i+1}-b_i}\eins_{n-i} 
}} \right] F_i
\quad ,
$$
where $\eins_i$ denotes the $i\times i$ unit matrix.
Then $\varphi_K$ induces an isomorphism $\varphi:E_n\isomto F_n$ and
we have the natural injections
\begin{eqnarray*}
&E_i\injto \m^{a_i-a_{i+1}}E_{i+1} \quad,\quad E_i\hookleftarrow E_{i+1} &\\
&F_{i+1}\injto F_i \quad,\quad \m^{b_i-b_{i+1}}F_{i+1}\hookleftarrow F_i &
\quad .
\end{eqnarray*}
Observe that the compositions
$E_{i+1}\injto E_i\injto \m^{a_i-a_{i+1}}E_{i+1}$
and 
$E_i\injto\m^{a_i-a_{i+1}}E_{i+1}\injto\m^{a_i-a_{i+1}}E_i$
are both the injections induced by the inclusion $A\injto\m^{a_i-a_{i+1}}$.
Furthermore, if $a_i-a_{i+1}<0$ then the morphism of $k$-vectorspaces
$E_{i+1}\tensor k\to E_i\tensor k$ is of rank $i$ and the sequence
$$
E_{i+1}\tensor k\to E_i\tensor k\to 
(\m^{a_i-a_{i+1}}E_{i+1})\tensor k\to
(\m^{a_i-a_{i+1}}E_{i})\tensor k
$$
is exact. This shows that the tupel
$$
(\m^{a_i-a_{i+1}},\ 1\in\m^{a_i-a_{i+1}}\ ,\ \ E_{i+1}\injto E_i\ ,\ \ 
 \m^{a_i-a_{i+1}}E_{i+1}\hookleftarrow E_i\ ,\ \ i)
$$
is what we call a ``bf-morphism'' of rank $i$ (cf. definition \ref{bf}).
Observe now that if $a_i-a_{i+1}<0$ and $(f,g)$ is one of the
following two pairs of  morphisms:
\begin{eqnarray*}
& E\tensor k\overset{f}{\longrightarrow}(\m^{-a_i}E_i)\tensor k
      \overset{g}{\longrightarrow}(\m^{-a_{i+1}}E_{i+1})\tensor k  
      \quad , & \\
& E_i\tensor k \overset{g}{\longleftarrow}E_{i+1}\tensor k 
             \overset{f}{\longleftarrow} E_n\tensor k \quad ,     &
\end{eqnarray*}
then $\im(g\comp f)=\im(g)$. The above statements hold true also if we replace
the $E_i$-s by the $F_i$-s and the $a_i$-s by the $b_i$-s.
Observe finally that in the diagram
$$
\xymatrix @C=1ex@R-1ex{
             &             &          0 \ar[d]     &       &        \\
             &             &
\ker(E_n\tensor k\to E_0\tensor k)\ar[d]\ar[dr]&
             &        \\
0\ar[r]&
\ker(F_n\tensor k\to F_0\tensor k)\ar[r]\ar[dr]&
E_n\tensor k\cong F_n\tensor k\ar[d]\ar[r]&
\im(F_n\tensor k\to F_0\tensor k)\ar[r]&0 \\
             &             &\im(E_n\tensor k\to E_0\tensor k)\ar[d] &   &   \\
             &             &          0            &       &
\\ }       
$$
the oblique arrows are injections.

All these properties are summed up in the statement that the tupel
\begin{eqnarray*}
\Phi &:=& ((\m^{b_i-b_{i+1}},\ 1),\ (\m^{a_i-a_{i+1}},\ 1),\  
          E_i\injto \m^{a_i-a_{i+1}}E_{i+1},\ E_i\hookleftarrow
          E_{i+1}, \\
& & F_{i+1}\injto F_i,\ \m^{b_i-b_{i+1}}F_{i+1}\hookleftarrow F_i\ 
    (0\leq i\leq n-1),\ \varphi: E_n\isomto F_n)
\end{eqnarray*}
is a generalized isomorphism from $E$ to $F$ in the sense of definition
\ref{generalized isomorphism}, where for $a\leq 0$ we consider
$\m^a$ as an invertible $A$-module with global section $1\in\m^a$.
Observe that $\Phi$ does not depend on our choice of the bases for $E$ and 
$F$. Indeed, it is well-known that the sequence $(m_1,\dots,m_n)$ is
independent of such a choice and it is easy to see that
$E_n=\varphi_K^{-1}(F)\cap E$,\quad $F_n=\varphi_K(E_n)$ and
$$
E_i=E_n + \m^{a_i}E \quad , \quad F_i=F_n + \m^{b_i}F
$$
for $1\leq i\leq n-1$, where the $+$-sign means generation in $E_K$ 
and $F_K$ respectively.
Observe furthermore that by pull-back the 
generalized isomorphism $\Phi$ induces a generalized
isomorphism $f^*\Phi$ on a scheme $S$ for every morphism $f:S\to\Spec(A)$. 
Of course the morphisms 
$f^*E_{i+1}\to f^*E_i$ etc. will be in general no longer injective, but
this is not required in the definition.

\section{Notations}
\label{notation}
We collect some less common notations, which we will use 
freely in this paper:

\begin{itemize}
\item
For two integers $a\leq b$ we sometimes denote by $[a,b]$ the set
$\{c\in\Z\ |\ a\leq c\leq b\}$.
\item
For a $n\times n$-matrix with entries $a_{ij}$ in some ring, and for
two subsets $A$ and $B$ of cardinality $r$ of $\{1,\dots,n\}$, we will 
denote by $\det_{AB}(a_{ij})$ the determinant of the 
$r\times r$-matrix $(a_{ij})_{i\in A, j\in B}$.
\item
For a scheme $X$ we will denote by $\K_X$ the sheaf of total quotient
rings of $\Oo_X$.
\item
For a scheme $X$, a coherent sheaf $\E$ on $X$ and a point $x\in X$,
we denote by $\E[x]$ the fibre $\E\tensor_{\Oo_X}\kappa(x)$ of $\E$ at 
$x$.
\item
For $n\in\N$, the symbol $S_n$ denotes the symmetric group of 
permutations of the set $\{1,\dots,n\}$.
\end{itemize}

\section{Construction of the compactification}
\label{construction}

Let $X^{(0)}:=\Proj\Z[\text{$x_{00}$, $x_{ij}$ $(1\leq i,j\leq n)$}]$. 
We define closed subschemes 
$$
\xymatrix@R-1.2pc{
Y_0^{(0)} \ar@<-.5ex>@{^{(}->}[r] & 
Y_1^{(0)} \ar@<-.5ex>@{^{(}->}[r] & 
\quad\dots\quad \ar@<-.5ex>@{^{(}->}[r]     & 
Y_{n-1}^{(0)}             & 
\\
                              & 
Z_{n-1}^{(0)} \ar@<-.5ex>@{^{(}->}[r] \ar@{^{(}->}[u]& 
\quad\dots\quad \ar@<-.5ex>@{^{(}->}[r]         & 
Z_1^{(0)} \ar@<-.5ex>@{^{(}->}[r] \ar@{^{(}->}[u]    & 
Z_0^{(0)} }
$$
of $X^{(0)}$, by setting $Y_r^{(0)}:= V^+(\I_r^{(0)})$,\quad 
$Z_r^{(0)}:=V^+(\J_r^{(0)})$,
where
$\I_r^{(0)}$ is the  homogenous ideal in
$\Z[\text{$x_{00}$, $x_{ij}$ $(1\leq i,j\leq n)$}]$,
generated by all $(r+1)\times(r+1)$-subdeterminants of the matrix
$(x_{ij})_{1\leq i,j\leq n}$,
and where $\J_{r}^{(0)}=(x_{00})+\I_{n-r}^{(0)}$ 
for $0\leq r\leq n-1$. For $1\leq k\leq n$ let the scheme $X^{(k)}$
together with closed subschemes 
$Y_r^{(k)}$, $Z_r^{(k)}\subset X^{(k)}$ $(0\leq r\leq n-1)$ 
be inductively defined as follows:

$X^{(k)}\to X^{(k-1)}$ is the blowing up of $X^{(k-1)}$ along the
closed subscheme $Y_{k-1}^{(k-1)}\cup Z_{n-k}^{(k-1)}$. 
The subscheme $Y_{k-1}^{(k)}\subset X^{(k)}$ 
(respectively $Z_{n-k}^{(k)}\subset X^{(k)}$) 
is the inverse image of $Y_{k-1}^{(k-1)}$
(respectively of $Z_{n-k}^{(k-1)}$) under the morphism 
$X^{(k)}\to X^{(k-1)}$, and
for $r\neq k-1$ (respectively $r\neq n-k$) the subscheme 
$Y_r^{(k)}\subset X^{(k)}$
(respectively of $Z_r^{(k)}\subset X^{(k)}$) 
is the complete transform of 
$Y_r^{(k-1)}\subset X^{(k-1)}$ 
(respectively $Z_r^{(k-1)}\subset X^{(k-1)}$).
We set
$$
\KGln := X^{(n)} \qquad
\text{and} \qquad 
Y_r:=Y_r^{(n)}\ ,\
Z_r:=Z_r^{(n)}\quad
(0\leq r\leq n-1) \ .
$$
We are interested in finding a modular description for the compactification 
$\KGln$ of 
$\Gln=\Spec\Z[\text{$x_{ij}/x_{00}$ $(1\leq i,j\leq n)$,
                     $\det(x_{ij}/x_{00})^{-1}$}]$.

Let 
$(\alpha,\beta)\in S_n\times S_n$
and set
$$
x_{ij}^{(0)}(\alpha,\beta):= \frac{x_{\alpha(i),\beta(j)}}{x_{00}}
\qquad
(1\leq i,j\leq n)
\quad.
$$
For $1\leq k\leq n$ we define elements
\begin{eqnarray*}
y_{ji}(\alpha,\beta)\ ,
& 
z_{ij}(\alpha,\beta)
&
(1\leq i\leq k,\quad i<j\leq n)
\\
x_{ij}^{(k)}(\alpha,\beta)
& &
(k+1\leq i,j\leq n)
\end{eqnarray*}
of the function field 
$\Q(X^{(0)})=\Q(x_{ij}/x_{00}\ (1\leq i,j\leq n))$
of $X^{(0)}$ inductively as follows:
\begin{eqnarray*}
y_{ik}(\alpha,\beta):= 
\frac{x_{ik}^{(k-1)}(\alpha,\beta)}
{x_{kk}^{(k-1)}(\alpha,\beta)} 
\qquad\qquad\qquad\qquad\quad
& &
(k+1\leq i\leq n)
\quad,
\\
z_{kj}(\alpha,\beta):= 
\frac{x_{kj}^{(k-1)}(\alpha,\beta)}
{x_{kk}^{(k-1)}(\alpha,\beta)} 
\qquad\qquad\qquad\qquad\quad
& &
(k+1\leq j\leq n)
\quad,
\\
x_{ij}^{(k)}(\alpha,\beta):= 
\frac{x_{ij}^{(k-1)}(\alpha,\beta)}
{x_{kk}^{(k-1)}(\alpha,\beta)}
-
y_{ik}(\alpha,\beta)\ 
z_{kj}(\alpha,\beta)
& &
(k+1\leq i,j\leq n)
\quad.
\end{eqnarray*}
Finally, we set $t_0(\alpha,\beta) := t_0 := x_{00}$ and
$$
t_i(\alpha,\beta) := 
t_0 \cdot
\prod_{j=1}^i
x_{jj}^{(j-1)}(\alpha,\beta)
\qquad\qquad
(1\leq i\leq n) \quad.
$$
Observe, that for each $k\in\{0,\dots,n\}$, we have the following
decomposition of the matrix $[x_{ij}/x_{00}]$:

$$
\label{matrix decomposition}
\left[\frac{x_{ij}}{x_{00}}\right] = n_\alpha
\left[
\vcenter{
\hbox{\tiny
\xymatrix @R-3.7pc@C-4pc@!{
1 \ar@{.}[rrdd]
      &     &     &     &      & 
\ar@{--} '[lllll]+<1.3ex,0ex> '[llldd]+<1.3ex,0ex> '[dd]+0 []+0 \\
      &     &     &     &  0   &     \\
      &     &  1  &     &      &     \\
      & 
      \text{\tiny $y_{ij}(\alpha,\beta)$}
            &     &     &      &     \\
      &     &     &     & \eins_{n-k}    
                               &     \\
\ar@{--} '[uuuuu]+<0ex,-1.3ex> '[uuurr]+<0ex,-1.3ex> '[rr]+0 []+0
      &     &     &     &      &  }}}
\right]
\left[
\vcenter{
\hbox{\tiny
\xymatrix @R-6.8pc@C-6.2pc@! {
\frac{t_1(\alpha,\beta)}{t_0} \ar@{.}[rrdd]
      &     &  0  & 
                        &     &      \\
      &     &     &     &  0  &      \\
  0   &     &
             \frac{t_k(\alpha,\beta)}{t_0} 
                  &     &     &      \\
      &     &     &     &     &      \\
      &  0  &     &     & 
                         \frac{t_k(\alpha,\beta)}{t_0}
                         [x_{ij}^{(k)}(\alpha,\beta)]
                              &      \\
      &     &     &     &     &   }}}
\right]
\left[
\vcenter{
\hbox{\tiny
\xymatrix @R-3.8pc@C-4pc@!{
1 \ar@{.}[rrdd]
      &     &     &     &      & 
\ar@{--} '[lllll]+<1.3ex,0ex> '[llldd]+<1.3ex,0ex> '[dd]+0 []+0  \\
      &     &     &
                  \, z_{ij}(\alpha,\beta)
                        &      &     \\
      &     &  1  &     &      &     \\
      &     &     &     &      &     \\
      &  0  &     &     & \eins_{n-k} &   \\
\ar@{--} '[uuuuu]+<0ex,-1.3ex> '[uuurr]+<0ex,-1.3ex> '[rr]+0 '[]+0 
             &     &     &     &      &  }}}
\right]
n_\beta^{-1}
$$
Here, $n_\alpha$ is the permutation matrix associated to $\alpha$,
i.e. the matrix, whose entry in the $i$-th row and $j$-th column is 
$\delta_{i,\alpha(j)}$.
For convenience, we define for each $l\in\{0,\dots,n\}$ a bijection
$\iota_l:\{1,\dots,n+1\} \isomto \{0,\dots,n\}$, by setting
$$
\iota_l(i)=
\left\{
\begin{array}{lll}
i & \text{if\quad $1\leq i\leq l$} \\
0 & \text{if\quad $i=l+1$} \\
i-1 & \text{if\quad $l+2\leq i\leq n+1$}
\end{array}
\right.
$$
for $1\leq i\leq n+1$. 
With this notaton, we define for each triple 
$(\alpha,\beta,l)\in S_n\times S_n\times [0,n]$
polynomial subalgebras $R(\alpha,\beta,l)$ of $\Q(\KGln)=\Q(X^{(0)})$
together with ideals 
$\I_r(\alpha,\beta,l)$ and
$\J_r(\alpha,\beta,l)$\ 
($0\leq r\leq n-1$)
of $R(\alpha,\beta,l)$ as follows:
{\small
\begin{eqnarray*}
R(\alpha,\beta,l) &:=&
\Z\left[
\frac{t_{\iota_l(i+1)}(\alpha,\beta)}{t_{\iota_l(i)}(\alpha,\beta)} \
(1\leq i\leq n),\
y_{ji}(\alpha,\beta),\
z_{ij}(\alpha,\beta) \
(1\leq i<j\leq n) 
\right] ,
\\
\I_{r}(\alpha,\beta,l) &:=&
\left(
\frac{t_{\iota_l(r+2)}(\alpha,\beta)}{t_{\iota_l(r+1)}(\alpha,\beta)}
\right) \quad 
\text{if} \quad 
l\leq r\leq n-1 \quad
\text{and} \quad
\I_{r}(\alpha,\beta,l) := (1) \
\text{else,}
\\
\J_{r}(\alpha,\beta,l) &:=&
\left(
\frac{t_{\iota_l(n-r+1)}(\alpha,\beta)}{t_{\iota_l(n-r)}(\alpha,\beta)}
\right) \
\text{if} \
n-l\leq r\leq n-1 \
\text{and} \
\J_{r}(\alpha,\beta,l) := (1) \
\text{else.}
\end{eqnarray*}
}
\begin{proposition}
\label{local description}
There is a covering of $\KGln$ by open affine pieces
$X(\alpha,\beta,l)$\ 
$((\alpha,\beta,l)\in S_n\times S_n\times [0,n])$,
such that 
$
\Gamma(X(\alpha,\beta,l),\Oo)=R(\alpha,\beta,l)
$
(equality as subrings of the function field $\Q(\KGln)$).
Furthermore, for $0\leq r\leq n-1$ the ideals
$\I_r(\alpha,\beta,l)$ and
$\J_r(\alpha,\beta,l)$ of
$R(\alpha,\beta,l)$
are the defining ideals for the closed subschemes
$Y_r(\alpha,\beta,l):=Y_r \cap X(\alpha,\beta,l)$ and
$Z_r(\alpha,\beta,l):=Z_r \cap X(\alpha,\beta,l)$
respectively.
\end{proposition}
\begin{proof}
We  make the blowing-up procedure explicit, in terms of open
affine coverings.
For each $k\in\{0,\dots,n\}$ we define a finite 
index set $\Pc_k$, consisting of all pairs
$$
(p,q)=
\left(
\left[
\begin{array}{c}
p_0 \\ : \\ p_k
\end{array}
\right]
\ ,\
\left[
\begin{array}{c}
q_0 \\ : \\ q_k
\end{array}
\right]
\right)
\ \in\ 
\{0,\dots,n\}^{k+1}\times\{0,\dots,n\}^{k+1}
$$
with the property that $p_i\neq p_j$ and $q_i\neq q_j$ for $i\neq j$
and that $p_i=0$ for some $i$, if and only if $q_i=0$.
Observe that for each $k\in\{0,\dots,n\}$ there is a surjection 
$S_n\times S_n\times\{0,\dots,n\} \to \Pc_k$,
which maps the triple $(\alpha,\beta,l)$ to the element
$$
(p,q)=
\left(
\left[
\begin{array}{c}
\alpha(\iota_l(1)) \\ : \\ \alpha(\iota_l(k+1))
\end{array}
\right]
\ ,\
\left[
\begin{array}{c}
\beta(\iota_l(1)) \\ : \\ \beta(\iota_l(k+1))
\end{array}
\right]
\right)
$$
of $\Pc_k$. (Here we have used the convention that $\alpha(0):=0$ for any
permutation $\alpha\in S_n$). Furthermore, this surjection is in fact
a bijection in the case of $k=n$.
Let $(p,q)\in\Pc_k$ and chose an element $(\alpha,\beta,l)$ in its
preimage under the surjection $S_n\times
S_n\times\{0,\dots,n\}\to\Pc_k$.
We define subrings $R^{(k)}(p,q)$ of 
$\Q(x_{ij}/x_{00}\ (1\leq i,j\leq n))$ together with ideals
$\I^{(k)}_r(p,q)$,\ $\J^{(k)}_r(p,q)$\ $(0\leq r\leq n)$, distinguishing
three cases.
\vspace{5mm}

First case: $0\leq l\leq k-1$
{\scriptsize
\begin{eqnarray*}
R^{(k)}(p,q) &:=&
\Z\left[
\frac{t_{\iota_l(i+1)}(\alpha,\beta)}{t_{\iota_l(i)}(\alpha,\beta)}\ 
(1\leq i\leq k),\ 
y_{ji}(\alpha,\beta),\
z_{ij}(\alpha,\beta)\
\text{{\scriptsize
$ 
\left(
\begin{array}{ll}
1\leq i\leq k,\\ 
i<j\leq n
\end{array}
\right),
$
}}
\right.
\\ & &\qquad\qquad\qquad\qquad\qquad\qquad\qquad\qquad
\left.
x^{(k)}_{ij}(\alpha,\beta)\ 
(k+1\leq i,j\leq n)
\right]
\\
\\
\I^{(k)}_r(p,q) &:=&
\left\{
\begin{array}{lll}
(1) 
&\text{if $r\in[0,l-1]$}
\\
\left(
\frac{t_{\iota_l(r+2)}(\alpha,\beta)}{t_{\iota_l(r+1)}(\alpha,\beta)}
\right)
&\text{if $r\in[l,k-1]$}
\\
\left(
\text{det}_{AB}(x^{(k)}_{ij}(\alpha,\beta))\
\text{{\scriptsize
$
\left(
\begin{array}{ll}
A,B\subseteq\{k+1,\dots,n\} \\
\card A=\card B=r+1-k
\end{array}
\right)
$
}}
\right)
&\text{if $r\in[k,n-1]$}
\end{array}
\right.
\\
\\
\J^{(k)}_r(p,q) &:=&
\left\{
\begin{array}{ll}
(1)
&\text{if $r\in[0,n-l-1]$}
\\
\left(
\frac{t_{\iota_l(n-r+1)}(\alpha,\beta)}{t_{\iota_l(n-r)}(\alpha,\beta)}
\right)
&\text{if $r\in[n-l,n-1]$}
\end{array}
\right.
\end{eqnarray*}
}
\vspace{5mm}

Second case: $l=k$
{\scriptsize
\begin{eqnarray*}
R^{(k)}(p,q) &:=&
\Z\left[
\frac{t_{\iota_l(i+1)}(\alpha,\beta)}{t_{\iota_l(i)}(\alpha,\beta)}\ 
(1\leq i\leq k),\ 
y_{ji}(\alpha,\beta),\
z_{ij}(\alpha,\beta)\ 
\text{{\scriptsize
$
\left(
\begin{array}{ll}
1\leq i\leq k,\\ 
i<j\leq n
\end{array}
\right),
$
}}
\right.
\\ & &\qquad\qquad\qquad\qquad\qquad\qquad\qquad
\left.
\frac{t_k(\alpha,\beta)}{t_0}x^{(k)}_{ij}(\alpha,\beta)\ 
(k+1\leq i,j\leq n)
\right]
\\
\\
\I^{(k)}_r(p,q) &:=&
\left\{
\begin{array}{lll}
(1) 
&\text{\!\!\!\!\!\!if $r\in[0,l-1]$}
\\
\left(
\text{det}_{AB}
\left(
\frac{t_k(\alpha,\beta)}{t_0}x^{(k)}_{ij}(\alpha,\beta)
\right)\
\text{{\scriptsize
$
\left(
\begin{array}{ll}
A,B\subseteq\{k+1,\dots,n\} \\
\card A=\card B=r+1-k
\end{array}
\right)
$
}}
\right)
&\text{\!\!\!\!\!\!if $r\in[l,n-1]$}
\end{array}
\right.
\\
\\
\J^{(k)}_r(p,q) &:=&
\left\{
\begin{array}{ll}
(1)
&\text{if $r\in[0,n-l-1]$}
\\
\left(
\frac{t_{\iota_l(n-r+1)}(\alpha,\beta)}{t_{\iota_l(n-r)}(\alpha,\beta)}
\right)
&\text{if $r\in[n-l,n-1]$}
\end{array}
\right.
\end{eqnarray*}
}

\vspace{5mm}
Third case: $k+1\leq l\leq n$
{\scriptsize
\begin{eqnarray*}
R^{(k)}(p,q) &:=&
\Z\left[
\frac{t_{\iota_l(i+1)}(\alpha,\beta)}{t_{\iota_l(i)}(\alpha,\beta)}\ 
(1\leq i\leq k),\
\frac{t_0}{t_{k+1}(\alpha,\beta)},
\right.
\\ & &
\left.
\quad\ 
y_{ji}(\alpha,\beta),\
z_{ij}(\alpha,\beta)\ 
\text{{\scriptsize
$
\left(
\begin{array}{ll}
1\leq i\leq k+1,\\ 
i<j\leq n
\end{array}
\right),
$
}}
\ x^{(k+1)}_{ij}(\alpha,\beta)\ 
(k+2\leq i,j\leq n)
\right]
\\
\\
\I^{(k)}_r(p,q) &:=&
\left\{
\begin{array}{ll}
(1) 
&\text{\!\!\!\!\!\!if $r\in[0,k]$}
\\
\left(
\text{det}_{AB}
\left(
x^{(k+1)}_{ij}(\alpha,\beta)
\right)\
\text{{\scriptsize
$
\left(
\begin{array}{ll}
A,B\subseteq\{k+2,\dots,n\} \\
\card A=\card B=r-k
\end{array}
\right)
$
}}
\right)
&\text{\!\!\!\!\!\!if $r\in[k+1,n-1]$}
\end{array}
\right.
\\
\\
\J^{(k)}_r(p,q) &:=&
\left\{
\begin{array}{lll}
\left(
\frac{t_0}{t_{k+1}(\alpha,\beta)},\ 
\text{det}_{AB}(x^{(k+1)}_{ij}(\alpha,\beta))\
\right. & \left.
\text{{\scriptsize
$
\left(
\begin{array}{ll}
A,B\subseteq\{k+2,\dots,n\} \\
\card A=\card B=n-r-k
\end{array}
\right)
$
}}
\right)
\\
&\qquad\qquad\qquad\text{if $r\in[0,n-k-1]$}
\\
\\
\left(
\frac{t_{\iota_l(n-r+1)}(\alpha,\beta)}{t_{\iota_l(n-r)}(\alpha,\beta)}
\right)
&\qquad\qquad\qquad\text{if $r\in[n-k,n-1]$}
\end{array}
\right.
\end{eqnarray*}
}

Observe that the objects 
$R^{(k)}(p,q)$, $\I^{(k)}_r(p,q)$, $\J^{(k)}_r(p,q)$ thus defined,
depend indeed only on
$(p,q)$ and not on the chosen element $(\alpha,\beta,l)$.
By induction on $k$ one shows that $X^{(k)}$ is covered by open affine
pieces $X^{(k)}(p,q)$\ ($(p,q)\in\Pc_k$), such that
$\Gamma(X^{(k)}(p,q),\Oo)=R^{(k)}(p,q)$\ (equality as subrings of the
function field $\Q(X^{(k)})$), and such that  the ideals
$\I_r^{(k)}(\alpha,\beta)$ and $\J_r^{(k)}(\alpha,\beta)$  are the
defining ideals of the closed subschemes
$Y^{(k)}_r\cap X^{(k)}(p,q)$ and $Z^{(k)}_r\cap X^{(k)}(p,q)$
respectively.
\end{proof}
\begin{corollary}
\label{nc}
The scheme $\KGln$ is smooth and projective over $\Spec \Z$ and contains
$\Gln$ as a dense open subset. The complement of $\Gln$ in $\KGln$ is the
union of the closed subschemes $Y_i$, $Z_i$ $(0\leq i\leq n-1)$, which
is a divisor with normal crossings.
Furthermore, we have $Y_i\cap Z_j=\emptyset$ for $i+j<n$.
\end{corollary}
\begin{proof}
This is immediate from the local description given in 
\ref{local description}.
\end{proof}

We will now define a certain toric scheme, which will
play an important role in the sequel.
Let $M:=\Z^n$, with canonical basis $e_1,\dots,e_n$.
For $m\in M$ we denote by $t^m$ the corresponding monomial in the ring
$\Z[M]$.
Furthermore, we write $t_i/t_0$ for the canonical generator $t^{e_i}$
of $\Z[M]$.
Let $N:=M^{\dual}$ be
the dual of $M$ with the dual basis 
$e_1^{\dual},\dots,e_n^{\dual}$.
For $0\leq l\leq n$ let $\sigma_l\subset N_{\Q}:=N\tensor\Q$ be the cone
generated by the elements
\ $-\sum_{j=1}^i e_j^{\dual}$\ \ $(1\leq i\leq l)$
and the elements
$\sum_{j=i}^n e_j^{\dual}$\ \ $(l+1\leq i\leq n)$. In other words:
$$
\sigma_l=
\sum_{i=1}^l
\Q_+\cdot
\left(-\sum_{j=1}^i e_j^{\dual}\right) +
\sum_{i=l+1}^n
\Q_+\cdot
\left(\sum_{j=i}^n e_j^{\dual}\right) \quad.
$$
Let $\Sigma$ be the fan generated by all $\sigma_l$\ $(0\leq l\leq n)$ and let 
$\T:=X_\Sigma$ the associated toric scheme (over $\Z$).
See e.g. \cite{Danilov} for definitions. $\T$ is
covered by the open sets 
$\T_l:=X_{\sigma_l^{\dual}}=
\Spec\Z[t^m\ (m\in\sigma_l^{\dual}\cap M)]=
\Spec\Z[t_{\iota_l(i+1)}/t_{\iota_l(i)}\ (1\leq i\leq n)]$.
Observe that there are Cartier divisors 
$Y_{r,\T}$, $Z_{r,\T}$\ $(0\leq r\leq n-1)$ on $\T$, such that 
for each $l\in\{0,\dots,n\}$ over the open part $\T_l\subset\T$, 
\begin{eqnarray*}
Y_{r,\T}\  & \text{is given by the equation} &
\left\{
\begin{array}{ll}
1 & \text{if $0\leq r\leq l-1$} \\
t_{\iota_l(r+2)}/t_{\iota_l(r+1)} &  \text{if $l\leq r\leq n-1$}
\end{array}
\right.
\\
Z_{r,\T}\ & \text{is given by the equation} &
\left\{
\begin{array}{ll}
1 &  \text{if $0\leq r\leq n-l-1$} \\
t_{\iota_l(n-r+1)}/t_{\iota_l(n-r)} & \text{if $n-l\leq r\leq n-1$} 
\end{array}
\right.
\end{eqnarray*}
Observe furthermore that 
$Y_{i,\T}\cap Z_{j,\T}=\emptyset$ for $i+j<n$ and that
for each $r\in\{1,\dots,n\}$, multiplication 
by $t_r/t_0$ establishes an isomorphism
$$
\Oo_{\T}\left(-\sum_{i=0}^{n-r}Z_{i,\T}\right) \Isomto
\Oo_{\T}\left(-\sum_{i=0}^{r-1}Y_{i,\T}\right)
\quad.
$$
\begin{lemma}
\label{T modular}
The toric scheme $\T$ together with the ``universal'' tupel
$$
( \Oo_{\T}(Y_{i,\T}),\ \one_{\Oo_{\T}(Y_{i,\T})},\
\Oo_{\T}(Z_{i,\T}),\ \one_{\Oo_{\T}(Z_{i,\T})}\
(0\leq i\leq n-1),\ t_r/t_0\ (1\leq r\leq n)) 
$$
represents the functor, which to each scheme $S$ associates the set of
equivalence classes of tupels
$$
(\Ll_i,\ \lambda_i,\ \M_i,\ \mu_i\ (0\leq i\leq n-1), 
\varphi_r\ (1\leq r\leq n))  \quad,
$$
where the $\Ll_i$ and $\M_i$ are invertible $\Oo_S$-modules
with global sectons $\lambda_i$ and $\mu_i$ respectively, such that
for $i+j<n$ the zero-sets of $\lambda_i$ and $\mu_j$ do not intersect,
and where the $\varphi_r$ are isomorphisms
$$
\Tensor_{i=0}^{n-r}\M_i^{\dual}
\Isomto
\Tensor_{i=0}^{r-1}\Ll_i^{\dual}
\quad.
$$
Here two tupels
$(\Ll_i,\ \lambda_i,\ \M_i,\ \mu_i\ (0\leq i\leq n-1), 
 \varphi_r\ (1\leq r\leq n))$
and
$(\Ll'_i,\ \lambda'_i,\ \M'_i,\ \mu'_i\ (0\leq i\leq n-1), 
 \varphi'_r\ (1\leq r\leq n))$
are called equivalent, if there exist isomorphisms
$\Ll_i\isomto\Ll'_i$ and $\M_i\isomto\M'_i$ for $0\leq i\leq n-1$,
such that all the obvious diagrams commute.
\end{lemma}
\begin{proof}
Let $S$ be a scheme and
$(\Ll_i,\ \lambda_i,\ \M_i,\ \mu_i\ (0\leq i\leq n-1), 
 \varphi_r\ (1\leq r\leq n))$
a tupel defined over $S$, which has the properties stated in the lemma.
Let us first consider the case, where all the sheaves $\Ll_i$, $\M_i$
are trivial and where there exists an $l\in\{0,\dots,n\}$, such that 
$\lambda_i$ and $\mu_j$ is nowhere vanishing  for $0\leq i<l$
and $0\leq j<n-l$ respectively.
Observe that under theses
conditions there exists a unique set of trivializations
$\Ll_i\isomto\Oo_S$,\ $\M_i\isomto\Oo_S$,\ $(0\leq i\leq n)$ such that 
$\lambda_i\mapsto 1$ for $0\leq i<l$,\quad 
$\mu_j\mapsto 1$ for $0\leq j<n-l$, and such that the diagrams
$$
\xymatrixnocompile{
\Tensor_{i=0}^{n-r}\M_i^{\dual}  \ar[rd]^{\sim} \ar[rr]^{\varphi_r}_{\sim} &
                                         &
\Tensor_{i=0}^{r-1}\Ll_i^{\dual} \ar[ld]_{\sim}
\\
      &
\Oo_S & }
$$
commute for $1\leq r\leq n$. Let 
$a_\nu\in\Gamma(S,\Oo_S)$\ $(1\leq \nu\leq n)$ be defined by 
$\lambda_i\mapsto a_{i+1}$ for $l\leq i\leq n-1$ and 
$\mu_j\mapsto a_{n-j}$ for $n-l\leq j\leq n-1$, and let $f_l:S\to\T_l$
be the morphism defined by 
$f_l^*(t_{\iota_l(\nu+1)}/ t_{\iota_l(\nu)})=
a_\nu$\ $(1\leq \nu\leq n)$.
It is straightforward to check that the induced morphism $f:S\to\T$
does not depend on the chosen number $l$ and that it is unique with
the property that the pull-back under $f$ of the universal tupel is
equivalent to the given one on $S$.

Returning to the general case, observe that there is an open covering 
$S=\cup_k U_k$,
such that for each $k$ there exists an $l$ with the property that over
$U_k$ all the  $\Ll_i$, $\M_i$ are trivial and that $\lambda_i$ and $\mu_j$ is
nowhere vanishing over $U_k$ for $0\leq i<l$ and $0\leq j< n-l$. 
The above construction shows that there exists a unique morphism
$f:S\to\T$ such that for each $k$ the restriction to $U_k$ of the
pull-back under $f$ of the universal tupel is equivalent to the
restriction to $U_k$ of the given one. Thus it remains only to show
that the isomorphisms defining the equivalences over the $U_k$ glue
together to give a global equivalence of the pull-back of the
universal tupel with the given one.
However, this is clear, since it is easy to see that there exists at
most one set of isomorphisms $\Ll_i\isomto\Ll'_i$,\ $\M_i\isomto\M'_i$
establishing an equvalence between two tuples 
$(\Ll_i,\ \lambda_i,\ \M_i,\ \mu_i,\ \varphi_r)$ and
$(\Ll'_i,\ \lambda'_i,\ \M'_i,\ \mu'_i,\ \varphi'_r)$.
\end{proof}
For each pair $(\alpha,\beta)\in S_n\times S_n$ we define the open subset
$X(\alpha,\beta)\subseteq\KGln$ as the union of the open affines 
$X(\alpha,\beta,l)$ $(0\leq l\leq n)$. Let
\begin{eqnarray*}
U^- &:=& \Spec\Z[y_{ji}\ (1\leq i<j\leq n)] \quad, \\
U^+ &:=& \Spec\Z[z_{ij}\ (1\leq i<j\leq n)] \quad.
\end{eqnarray*}
Let $y:X(\alpha,\beta)\to U^-$ 
(respectively $z:X(\alpha,\beta)\to U^+$)
be the morphism defined by the property that 
$y^*(y_{ji})=y_{ji}(\alpha,\beta)$
(respectively $z^*(z_{ij})=z_{ij}(\alpha,\beta)$) for
$1\leq i,j\leq n$.
Observe that just as in the case of $\T$, multiplication by the
rational function $t_r(\alpha,\beta)/t_0$ provides an isomorphism
$$
\Oo_{X(\alpha,\beta)}
\left(
-\sum_{i=0}^{n-r}Z_i(\alpha,\beta)
\right)
\Isomto
\Oo_{X(\alpha,\beta)}
\left( 
-\sum_{i=0}^{r-1}Y_i(\alpha,\beta)
\right)
$$
for $1\leq r\leq n$, where $Y_i(\alpha,\beta)$ 
(respectively $Z_i(\alpha,\beta)$) denotes the restriction of $Y_i$
(respectively $Z_i$) to the open set $X(\alpha,\beta)$.
Thus, by lemma \ref{T modular}, the tupel
$$
(\Oo(Y_i(\alpha,\beta)),\ \one,\
 \Oo(Z_i(\alpha,\beta)),\ \one \
 (i\in[0,n-1]),\
 t_r(\alpha,\beta)/t_0\
 (r\in[1,n]))
$$
defines a morphism $t:X(\alpha,\beta)\to\T$.
\begin{lemma}
\label{ytz}
The morphism 
$(y,t,z): X(\alpha,\beta) \to U^-\times\T\times U^+$
is an isomorphism.
\end{lemma}
\begin{proof}
Let $\Omega(\alpha,\beta)\subset X(\alpha,\beta)$ be the preimage of
$\Gln$ under the morphism $X(\alpha,\beta)\injto\KGln\to X^{(0)}$.
By definition of $\KGln$, we have for all $l\in\{0,\dots,n\}$:
\begin{eqnarray*}
\Omega(\alpha,\beta) &=&
X(\alpha,\beta,l) \setminus
\bigcup_{i=0}^{n-1}
\left(
Y_i(\alpha,\beta,l)\cup Z_i(\alpha,\beta,l)
\right)
\\ &=&
\Spec\Z[y_{ji}(\alpha,\beta),\ z_{ij}(\alpha,\beta)\ 
        (1\leq i<j\leq n),
\\ & & \qquad\qquad\qquad\qquad\qquad\qquad\qquad
        (t_i(\alpha,\beta)/t_0)^{\pm 1}\ 
        (1\leq i\leq n)]
\quad.
\end{eqnarray*} 
Let $T:=\Spec\Z[(t_i/t_0)^{\pm 1}]\subset\T$ be the Torus in
$\T$. We have an isomorphism 
$\Omega(\alpha,\beta)\isomto U^-\times T\times U^+$ defined by
$y_{ji}\mapsto y_{ji}(\alpha,\beta)$,\ 
$z_{ij}\mapsto z_{ij}(\alpha,\beta)$,\
$t_i/t_0\mapsto t_i(\alpha,\beta)/t_0$, and a commutative quadrangle
$$
\begin{CD}
X(\alpha,\beta) @>(y,t,z)>> U^-\times\T\times U^+ \quad \\
@AAA                                     @AAA     \\
\Omega(\alpha,\beta) @>\sim>> U^-\times T\times U^+ \quad,
\end{CD}
$$
where the vertical arrows are the natural inclusions.
Furthermore, the map $(y,t,z)$ induces an isomorphism
$X(\alpha,\beta,l)\isomto U^-\times\T_l\times U^+$ for 
$0\leq l\leq n$.
Using the fact that $X(\alpha,\beta)$ is separated and that 
$\Omega(\alpha,\beta)$ dense
in  $X(\alpha,\beta)$, the lemma now follows easily.
\end{proof}

\section{bf-morphisms and generalized isomorpisms}
\label{bf&gi}

\begin{definition}
\label{bf}
Let $S$ be a scheme, $\E$ and $\F$ two localy free $\Oo_S$-modules and
$r$ a nonnegative integer.
A {\em bf-morphism of rank $r$ from $\E$ to $\F$} is a tupel
$$
g=(\text{$\M$, $\mu$,\quad $\E\to\F$,\quad $\M\tensor\E\leftarrow\F$,
   \quad $r$})\ ,
$$
where $\M$ is an invertible $\Oo_S$-module and $\mu$ a global section
of $\M$ such that the following holds:
\begin{enumerate}
\item
The composed morphisms 
$\E\to\F\to\M\tensor\E$ and 
$\F\to\M\tensor\E\to\M\tensor\F$
are both induced by the morphism $\mu:\Oo_S\to\M$.
\item
For every point $x\in S$ with $\mu(x)=0$, the complex
$$
\E[x]\to\F[x]\to(\M\tensor\E)[x]\to(\M\tensor\F)[x]
$$
is exact and the rank of the morphism $\E[x]\to\F[x]$
equals r.
\end{enumerate}
\end{definition}

The letters ``bf'' stand for ``back and forth''.
As a matter of notation, we will sometimes write $g^{\sharp}$ for the
morphism $\E\to\F$ and $g^{\flat}$ for the morphism $\F\to\M\tensor\E$
occuring in the bf-morphism $g$. 
Note that in case $\mu$ is nowhere
vanishing, the number $\rk g:=r$ cannot be deduced from the other 
ingredients of $g$. 
Sometimes we will use the following more suggestive notation for
the bf-morphism $g$:
$$
g=\left(
\xymatrix@C-1pc{
\E
\ar[rr]^{r}_{(\M,\mu)}
& & 
\F
\ar@/_1.2pc/|{\tensor}[ll]
}
\right)\ .
$$
In situations where it is clear, what $(\M,\mu)$ and $r$ are,
we will sometimes omit these data from our notation:
$$
g=\left(
\xymatrix@C-1pc{
\E
\ar[rr]
& & 
\F
\ar@/_1.2pc/|{\tensor}[ll]
}
\right)\ .
$$

\begin{definition}
\label{generalized isomorphism}
Let $S$ be a scheme, $\E$ and $\F$ two locally free $\Oo_S$-modules of
rank $n$. A {\em generalized isomorphism from $\E$ to $\F$} is a tupel
\begin{eqnarray*}
\Phi &=& (\text{$\Ll_i$, $\lambda_i$, $\M_i$, $\mu_i$,\quad 
          $\E_i\to\M_i\tensor\E_{i+1}$,\quad $\E_i\leftarrow\E_{i+1}$},\\
     & &\ \text{$\F_{i+1}\to\F_i$,\quad
          $\Ll_i\tensor\F_{i+1}\leftarrow\F_i$\quad
          $(0\leq i\leq n-1)$,\quad $\E_n\isomto\F_n$})\ ,
\end{eqnarray*}
where 
$\E=\E_0$, $\E_1,\dots,\E_n$, $\F_n,\dots,\F_1$, $\F_0=\F$, 
are localy free
$\Oo_S$-modules of rank $n$ and the tupels
\begin{eqnarray*}
& &(\text{$\M_i$, $\mu_i$,\quad $\E_{i+1}\to\E_i$,\quad
          $\M_i\tensor\E_{i+1}\leftarrow\E_i$,\quad i})\\
\text{and}& &(\text{$\Ll_i$, $\lambda_i$,\quad $\F_{i+1}\to\F_i$,\quad
              $\Ll_i\tensor\F_{i+1}\leftarrow\F_i$,\quad i})\ 
\end{eqnarray*}
are bf-morphisms of rank $i$ for $0\leq i\leq n-1$, such that for each
$x\in S$ the
following holds:
\begin{enumerate}
\item
If $\mu_i(x)=0$ and $(f,g)$ is one of the following two pairs of
morphisms:
\begin{eqnarray*}
&    \E[x]
     \overset{f}{\longrightarrow}
     ((\tensor_{j=0}^{i-1}\M_j)\tensor\E_{i})[x]
     \overset{g}{\longrightarrow}
     ((\tensor_{j=0}^{i}\M_j)\tensor\E_{i+1})[x] \ ,& \\
&   \E_{i}[x]
    \overset{g}{\longleftarrow}
    \E_{i+1}[x]
    \overset{f}{\longleftarrow}
    \E_n[x] \ ,&
\end{eqnarray*}
then $\im(g\comp f)=\im(g)$. 
Likewise, if $\lambda_i(x)=0$ and $(f,g)$ is one of the following two pairs
of morphisms:
\begin{eqnarray*}
&  \F_n[x]
   \overset{f}{\longrightarrow}
   \F_{i+1}[x]
   \overset{g}{\longrightarrow}
   \F_{i}[x] \ , &\\
&  ((\tensor_{j=0}^{i}\Ll_j)\tensor\F_{i+1})[x]
   \overset{g}{\longleftarrow}
   ((\tensor_{j=0}^{i-1}\Ll_j)\tensor\F_{i})[x]
   \overset{f}{\longleftarrow}
   \F[x] \ ,&
\end{eqnarray*}
then $\im(g\comp f)=\im(g)$.
\item
In the diagram:
$$
\xymatrix @C-1ex@R-1ex{
             &             &          0 \ar[d]     &       &        \\
             &             &
\ker(\E_n[x]\to\E_0[x])\ar[d]\ar[dr]&
             &        \\
0\ar[r]&
\ker(\F_n[x]\to\F_0[x])\ar[r]\ar[dr]&
\E_n[x]\cong\F_n[x]\ar[d]\ar[r]&
\im(\F_n[x]\to\F_0[x])\ar[r]&0 \\
             &             &\im(\E_n[x]\to\E_0[x])\ar[d] &       &        \\
             &             &          0            &       &
\\ }       
$$
the oblique arrows are injections.
\end{enumerate}
\end{definition}

\begin{definition}
A {\em quasi-equivalence} between two generalized isomorphisms
\begin{eqnarray*}
\Phi &=& (\text{$\Ll_i$, $\lambda_i$, $\M_i$, $\mu_i$,\quad 
          $\E_i\to\M_i\tensor\E_{i+1}$,\quad $\E_i\leftarrow\E_{i+1}$},\\
     & &\ \text{$\F_{i+1}\to\F_i$,\quad
          $\Ll_i\tensor\F_{i+1}\leftarrow\F_i$\quad
          $(0\leq i\leq n-1)$,\quad $\E_n\isomto\F_n$})\ ,\\
\Phi' &=& (\text{$\Ll'_i$, $\lambda'_i$, $\M'_i$, $\mu'_i$,\quad 
          $\E'_i\to\M'_i\tensor\E'_{i+1}$,\quad $\E'_i\leftarrow\E'_{i+1}$},\\
     & &\ \text{$\F'_{i+1}\to\F'_i$,\quad
          $\Ll'_i\tensor\F'_{i+1}\leftarrow\F'_i$\quad
          $(0\leq i\leq n-1)$,\quad $\E'_n\isomto\F'_n$})
\end{eqnarray*}
from $\E$ to $\F$ consists in 
isomorphisms $\Ll_i\isomto\Ll'_i$ and 
$\M_i\isomto\M'_i$
for $0\leq i\leq n-1$, and 
isomorphisms $\E_i \isomto\E'_i$ and 
$\F_i \isomto\F'_i$
for $0\leq i\leq n$, 
such that all the obvious diagrams are commutative.
A quasi-equivalence between $\Phi$ and $\Phi'$ is called an
{\em equivalence}, if the isomorphisms $\E_0\isomto\E'_0$ and
$\F_0\isomto\F'_0$ are in fact the identity on $\E$ and $\F$
respectively.
\end{definition}
After these general definitions, we now return to our scheme
$\KGln$. The notations are as in the previous section.

From the matrix-decomposition on page \pageref{matrix decomposition}
(for $k=n$) we see
that the matrix $[x_{ij}/x_{00}]_{1\leq i,j \leq n}$ has entries in the
subspace $\Gamma(\KGln,\Oo(\sum_{i=0}^{n-1}Z_i))$ of the function
field  $\Q(\KGln)$ of $\KGln$. Therefore it defines a morphism
$$
\x: E_0 
     \longrightarrow 
     \Oo\left(\sum_{i=0}^{n-1}Z_i\right) \cdot F_0 \quad ,
$$
where $E_0=F_0=\oplus^n\Oo_{\KGln}$.

Let $E_n\subset E_0$ be the preimage under $\x$ of 
$F_0 \subset \Oo(\sum_{i=0}^{n-1}Z_i) \cdot F_0$ and let
$F_n\subset F_0$ be the image under $\x$ of $E_n$.
Thus $\x$ induces a morphism
$$
E_n \longrightarrow F_n \quad ,
$$
which we again denote by $\x$.
For $1\leq i \leq n-1$ we define $\Oo_{\KGln}$-submodules $E_i$ and
$F_i$ of $\oplus^n\K_{\KGln}$ as follows:
\begin{eqnarray*}
E_i &:=& E_n + \Oo\left(-\sum_{j=0}^{i-1}Z_j\right)\cdot E_0 \\
F_i &:=& F_n + \Oo\left(-\sum_{j=0}^{i-1}Y_j\right)\cdot F_0
\end{eqnarray*}
(the plus-sign means generation in $\oplus^n\K_{\KGln}$).
Observe that for $0 \leq i\leq n-1$ we have the following natural injections:
\begin{eqnarray*}
& E_i \injto \Oo(Z_i)\cdot E_{i+1}\quad,\quad  E_i \hookleftarrow E_{i+1} &\\
& F_{i+1} \injto F_{i}\quad,\quad \Oo(Y_i)\cdot F_{i+1} \hookleftarrow F_i &
\end{eqnarray*}

\begin{proposition}
\label{Phi_univ}
The tupel 
\begin{eqnarray*}
\Phi_{\text{univ}} &:=& 
(\text{$\Oo(Y_i)$, $\one_{\Oo(Y_i)}$, $\Oo(Z_i)$, $\one_{\Oo(Z_i)}$,\
       $E_i \injto \Oo(Z_i)\cdot E_{i+1}$, \
       $E_i \hookleftarrow E_{i+1}$,} \\
& &\ \text{$F_{i+1} \injto F_{i}$,\
           $\Oo(Y_i)\cdot F_{i+1} \hookleftarrow F_i$\
           $(0\leq i\leq n-1)$, \ $\x:E_n\to  F_n$})
\end{eqnarray*}
is a generalized isomorphism from $\oplus^n\Oo_{\KGln}$ to itself.
\end{proposition}
\begin{proof}
It suffices to show that for each $(\alpha,\beta)\in S_n\times S_n$
the restriction of $\Phi_{\text{univ}}$ to the open set
$X(\alpha,\beta)$ is a generalized isomorphism from
$\oplus^n\Oo_{X(\alpha,\beta)}$ to itself.
Let $\z(\alpha,\beta)$ ($\y(\alpha,\beta)$) 
be the upper (lower) triangular $n\times n$ matrix with
$1$ on the diagonal and entries $z_{ij}(\alpha,\beta)$\
($y_{ji}(\alpha,\beta)$) over (under) the diagonal $(1\leq i<j\leq n)$. 
For $0\leq i\leq n$ we define 
\begin{eqnarray*}
E_i(\alpha,\beta) &:=&
\z(\alpha,\beta)\cdot n_\beta^{-1}\cdot E_i|_{X(\alpha,\beta)}
\quad,
\\
F_i(\alpha,\beta) &:=&
\y(\alpha,\beta)^{-1}\cdot n_\alpha^{-1}\cdot F_i|_{X(\alpha,\beta)}
\quad.
\end{eqnarray*}
Here we interprete the matrices $\z(\alpha,\beta)\cdot n_\beta^{-1}$ and 
$\y(\alpha,\beta)^{-1}\cdot n_\alpha^{-1}$ as automorphisms of 
$\oplus^n\K_{X(\alpha,\beta)}$. Accordingly we view the sheaves
$E_i(\alpha,\beta)$ and $F_i(\alpha,\beta)$ as subsheaves of
$\oplus^n\K_{X(\alpha,\beta)}$.
We have to show that the tupel
\begin{eqnarray*}
\Phi(\alpha,\beta) &:=&
(\Oo(Y_i(\alpha,\beta))\ ,\ 
 \one_{\Oo(Y_i(\alpha,\beta))}\ ,\ 
 \Oo(Z_i(\alpha,\beta))\ ,\ 
 \one_{\Oo(Z_i(\alpha,\beta))}\ ,\ 
\\
& &\ \ 
 E_i(\alpha,\beta) \injto 
 \Oo(Z_i(\alpha,\beta))\cdot E_{i+1}(\alpha,\beta)\ ,\  
 E_i(\alpha,\beta) \hookleftarrow E_{i+1}(\alpha,\beta),
\\
& &\ \
 F_{i+1}(\alpha,\beta) \injto F_{i}(\alpha,\beta)\ ,\quad 
 \Oo(Y_i(\alpha,\beta))\cdot F_{i+1}(\alpha,\beta) \hookleftarrow 
 F_i(\alpha,\beta) 
\\
& &\ \
 (0\leq i\leq n-1)\ ,
\\
& &\ \ 
  \y(\alpha,\beta)^{-1}n_\alpha^{-1}\x n_\beta\z(\alpha,\beta)^{-1}:
  E_n(\alpha,\beta)\isomto  F_n(\alpha,\beta))
\end{eqnarray*}
is a generalized isomorphism from 
$\oplus^n\Oo_{X(\alpha,\beta)}$ to itself.

We have for $0\leq i\leq n$ the following equality of subsheaves of 
$\oplus^n\K_{X(\alpha,\beta)}$:
\begin{eqnarray*}
E_i(\alpha,\beta) &=&
\bigoplus_{j=1}^{n-i}\ 
\Oo\left(
-\sum_{\nu=0}^{i-1} Z_\nu(\alpha,\beta)
\right)
\oplus
\bigoplus_{j=n-i+1}^{n}
\Oo\left(
-\sum_{\nu=0}^{n-j} Z_\nu(\alpha,\beta)
\right)
\quad,
\\
F_i(\alpha,\beta) &=&
\bigoplus_{j=1}^{i}\ 
\Oo\left(
-\sum_{\nu=0}^{j-1} Y_\nu(\alpha,\beta)
\right)
\oplus
\bigoplus_{j=i+1}^{n}
\Oo\left(
-\sum_{\nu=0}^{i-1} Y_\nu(\alpha,\beta)
\right)
\quad.
\end{eqnarray*}
This is easily checked by restricting both sides of the equations
to the open subsets $X(\alpha,\beta,l)$,\ $(0\leq l\leq n)$ of
$X(\alpha,\beta)$ and using \ref{local description}.
Observe that the morphisms
$$
E_i(\alpha,\beta) \injto 
\Oo(Z_i(\alpha,\beta))\cdot E_{i+1}(\alpha,\beta)\ ,\quad  
E_i(\alpha,\beta) \hookleftarrow E_{i+1}(\alpha,\beta)
$$
are described by the matrices
$$
\left[
\begin{array}{cc}
\eins_{n-i} & 0 \\
0 & \mu_i \eins_i
\end{array}
\right]
\quad\text{and}\quad
\left[
\begin{array}{cc}
\mu_i\eins_{n-i} & 0 \\
0 & \eins_i
\end{array}
\right]
\quad,
$$
and the morphisms
$$
F_{i+1}(\alpha,\beta) \injto F_{i}(\alpha,\beta)\ ,\quad 
\Oo(Y_i(\alpha,\beta))\cdot F_{i+1}(\alpha,\beta) \hookleftarrow 
F_i(\alpha,\beta) 
$$  
by the matrices
$$
\left[
\begin{array}{cc}
\eins_i & 0 \\
0 & \lambda_i \eins_{n-i}
\end{array}
\right]
\quad\text{and}\quad
\left[
\begin{array}{cc}
\lambda_i\eins_i & 0 \\
0 & \eins_{n-i}
\end{array}
\right]
$$
respectively, where we have abbreviated 
$\one_{\Oo(Y_i(\alpha,\beta))}$ by $\lambda_i$, and
$\one_{\Oo(Z_i(\alpha,\beta))}$ by $\mu_i$.
Furthermore the matrix-decomposition on page \pageref{matrix decomposition}
(for $k=n$) shows that 
$\y(\alpha,\beta)^{-1}n_\alpha^{-1}\x n_\beta\z(\alpha,\beta)^{-1}$ 
is the diagonal
matrix with entries
$(t_1(\alpha,\beta)/t_0,\dots,t_n(\alpha,\beta)/t_0)$.
With this information at hand, it is easy to see that
$\Phi(\alpha,\beta)$ is indeed a generalized isomorphism from 
$\oplus^n\Oo_{X(\alpha,\beta)}$ to itself.
\end{proof}
\begin{theorem}
\label{KGln modular}
Let $S$ be a scheme and $\Phi$ a generalized isomorphism from 
$\oplus^n \Oo_S$ to itself. Then there is a unique morphism
$f:S\to\KGln$ such that $f^*\Phi_{\text{univ}}$ is equivalent to
$\Phi$. In other words, the scheme $\KGln$ together with
$\Phi_{\text{univ}}$ represents the functor,
which to each scheme $S$ associates the set of equivalence classes of 
generalized isomorphisms from $\oplus^n \Oo_S$ to itself.
\end{theorem}
The proof of the theorem will be given in section \ref{proof}.

\begin{corollary}
\label{operation}
There is a (left) action of $\Gln\times\Gln$ on $\KGln$, which extends the
action 
$((\varphi,\psi),\Phi)\mapsto\psi\Phi\varphi^{-1}$
of $\Gln\times\Gln$ on $\Gln$.
The divisors $Z_i$ and $Y_i$ are invariant under this action.
\end{corollary}

\begin{proof}
The the morphism $(\Gln\times\Gln)\times\KGln\to\KGln$ 
defining the action is given on $S$-valued points by
$$
((\varphi,\psi),\Phi)\mapsto\Phi'
\quad,
$$
where $\Phi$ is a generalized isomorphism as in definition
\ref{generalized isomorphism} from $\E_0=\oplus^n\Oo_S$
to $\F_0=\oplus^n\Oo_S$ and $\Phi'$ is the generalized isomorphism
where for $2\leq i \leq n$ the bf-morphisms from $\E_i$ to $\E_{i-1}$, 
the ones from $\F_i$ to $\F_{i-1}$ and the isomorphism $\E_n\isomto \F_n$ 
are the same as in the tupel $\Phi$, and where the 
bf-morphisms
\begin{eqnarray*}
&(\M_0,\mu_0,\ \E_1\to\E_0,\ \M_0\tensor\E_0\ot\E_1,\ 0)& \\
\text{and}\quad 
&(\Ll_0,\lambda_0,\ \F_1\to\F_0,\ \Ll_0\tensor\F_0\ot\F_1,\ 0)&
\end{eqnarray*}
in the tupel $\Phi$ are replaced by the bf-morphisms
\begin{eqnarray*}
&(\M_0,\mu_0,\ \E_1\to\E_0\overset{\varphi}{\to}\E_0,\ 
             \M_0\tensor\E_1\ot\E_0\overset{\varphi^{-1}}{\ot}\E_0,\ 0)& \\
\text{and}\quad 
&(\Ll_0,\lambda_0,\ \F_1\to\F_0\overset{\psi}{\to}\F_0,\ 
             \Ll_0\tensor\F_1\ot\F_0\overset{\psi^{-1}}{\ot}\F_0,\ 0)&
\quad
\end{eqnarray*}
respectively.
The invariance of the divisors $Z_i$ and $Y_i$ is clear, since they
are defined by the vanishing of $\mu_i$ and $\lambda_i$ respectively.
\end{proof}

\section{Exterior powers}
\label{exterior}

\begin{lemma}
\label{local frames}
Let S be a scheme and $\E$, $\F$ two locally free $\Oo_S$-modules of
rank $n$. Let
$$
g=(\text{$\M$, $\mu$,\quad $\E\to\F$,\quad $\M\tensor\E\leftarrow\F$,
   \quad $r$}) 
$$
be a bf-morphism of rank $r$ from $\E$ to $\F$. Then each point 
$x\in S$ has an open neighbourhood $U$ such that over $U$ there exist
local frames $(e_1,\dots,e_n)$ and $(f_1,\dots,f_n)$ for $\E$ and $\F$
respectively with the property that the matrices for the morphisms
$$
\text{$\E\longrightarrow\F$\quad and\quad $\M\tensor\E\longleftarrow\F$}
$$
with respect to these frames are
$$
\text{
$
\left[
\begin{array}{cc}
\eins_r & 0 \\
0 & \mu/m\eins_{n-r}
\end{array}
\right]
$
\quad and \quad
$
\left[
\begin{array}{cc}
\mu\eins_r & 0 \\
0 & m\eins_{n-r}
\end{array}
\right]
$
}
$$
respectively, where $m$ is a nowhere vanishing section of $\M$ over $U$.
\end{lemma}
\begin{proof}
Restricting to a neighbourhood of $x$, we may assume 
that the sheaves $\M$, $\E$, $\F$ are free. Let 
$(\tilde{e}_1,\dots,\tilde{e}_n)$ and 
$(\tilde{f}_1,\dots,\tilde{f}_n)$ be global frames for $\E$ and $\F$
respectively. After permutation of their elements, and restricting to
a possibly smaller neighbourhood of $x$, we may further
assume that the morphisms
\begin{eqnarray*}
\langle \tilde{e}_1,\dots,\tilde{e}_r \rangle
& \longrightarrow &
\F/\langle \tilde{f}_{r+1},\dots,\tilde{f}_n \rangle
\\
\langle \tilde{f}_{r+1},\dots,\tilde{f}_n \rangle
& \longrightarrow &
\M\tensor\E/\langle \tilde{e}_1,\dots,\tilde{e}_r \rangle
\end{eqnarray*}
induced by $g^{\sharp}$ and $g^{\flat}$ respectively, are isomorphisms.
Let
\begin{eqnarray*}
&
\tilde{\E} := \langle \tilde{e}_1,\dots,\tilde{e}_r \rangle
\quad,\quad
\tilde{\E}' := 
\ker(\E\to\M\tensor\F/
\langle\tilde{f}_{r+1},\dots,\tilde{f}_n\rangle) 
& \\ &
\tilde{\F} := \ker(\F\to\E/
\langle\tilde{e}_1,\dots,\tilde{e}_r\rangle)
\qquad,\qquad
\tilde{\F}' := \langle\tilde{f}_{r+1},\dots,\tilde{f}_n\rangle
\ .
&
\end{eqnarray*}
Then we have direct-sum decompositions
$
\E=\tilde{\E}\oplus\tilde{\E}'
\ ,\ 
\F=\tilde{\F}\oplus\tilde{\F}'
\ ,
$
which are respected by $g^{\sharp}$ and $g^{\flat}$.
Let $m$ be a nowhere vanishing section of $\M$.
The frames $(e_1,\dots,e_n)$, $(f_1,\dots,f_n)$ 
of $\E$, $\F$, defined by setting
$e_i:=\tilde{e}_i$,\ $f_i:=g^{\sharp}(\tilde{e}_i)$ 
for $1\leq i\leq r$ and
$e_i:=(1/m)g^{\flat}(\tilde{f}_i)$,\ $f_i:=\tilde{f}_i$ 
for $r+1\leq i\leq n$, have the desired property.
\end{proof}
\begin{proposition}
\label{wedge g}
Let S be a scheme and $\E$, $\F$ two locally free $\Oo_S$-modules of
rank $n$. Let
$$
g=(\text{$\M$, $\mu$,\quad $\E\to\F$,\quad $\M\tensor\E\leftarrow\F$,
   \quad $i$}) 
$$
be a bf-morphism of rank $i$ from $\E$ to $\F$ and let $1\leq r \leq n$.
\begin{enumerate}
\item
There exists a unique morphism 
$$
\wedge^r g: 
\wedge^r\E \longrightarrow 
(\M^{\dual})^{\tensor\max(0,r-i)}\tensor\wedge^r\F
$$
with the following property:
If $(e_1,\dots,e_n)$ and $(f_1,\dots,f_n)$ are local frames for $\E$
and $\F$ respectively over an open set $U\subseteq S$, 
such that the matrices for the morphisms
$$
\text{
$\E\longrightarrow\F$ 
\quad and\quad
$\M\tensor\E\longleftarrow\F$
}
$$
with respect to these frames are
$$
\text{
$
\left[
\begin{array}{cc}
\eins_i & 0 \\
0 & \mu/m\eins_{n-i}
\end{array}
\right]
$
\quad and \quad
$
\left[
\begin{array}{cc}
\mu\eins_i & 0 \\
0 & m\eins_{n-i}
\end{array}
\right]
$
}
$$
respectively ($m$ being a nowhere vanishing section of $\M$ over $U$),
then $(\wedge^r g)|_U$ takes the form
$$
e_I\wedge e_J \mapsto 
m^{p-r}\tensor\mu^{\min(i,r)-p}\tensor f_I\wedge f_J \quad ,
$$
where $I\subseteq\{1,\dots,i\}$, $J\subseteq\{i+1,\dots,n\}$
with $\card I=p$, $\card J=r-p$ and where
$e_I=e_{i_1}\wedge\dots\wedge e_{i_p}$, if $I=\{i_1,\dots,i_p\}$ and
$i_1<\dots<i_p$. The $e_J,f_I,f_J$ are defined analoguosly.
\item
Similarly, there exists a unique morphism 
$$
\wedge^{-r} g:
\wedge^r\F \longrightarrow
\M^{\tensor\min(r,n-i)}\tensor\wedge^r\E
$$
with the following property: If $(e_1,\dots,e_n)$ and $(f_1,\dots,f_n)$
are local frames for $\E$ and $\F$ respectively over an open set
$U\subseteq S$, such that the matrices for the morphisms
$$
\text{
$\F\longrightarrow\M\tensor\E$
\quad and\quad
$\E\longleftarrow\F$
}
$$
with respect to these frames are

$$
\text{
$
\left[
\begin{array}{cc}
m\eins_{n-i} & 0 \\
0 & \mu\eins_{i}
\end{array}
\right]
$
\quad and \quad
$
\left[
\begin{array}{cc}
\mu/m\eins_{n-i} & 0 \\
0 & \eins_{i}
\end{array}
\right]
$
}
$$
respectively ($m$ being a nowhere vanishing section of $\M$ over $U$),
then $(\wedge^{-r} g)|_U$ takes the form
$$
f_I\wedge f_J \mapsto
m^p\tensor\mu^{\min(r,n-i)-p}\tensor e_I\wedge e_J \quad,
$$
where $I\subseteq\{1,\dots,n-i\}$, $J\subseteq\{n-i+1,\dots,n\}$
with $\card I=p$, $\card J=r-p$.
\end{enumerate}
\end{proposition}
\begin{proof}
1. An easy calculation shows that the morphism given by the prescription
$$
e_I\wedge e_J \mapsto 
m^{p-r}\tensor\mu^{\min(i,r)-p}\tensor f_I\wedge f_J 
\quad,
$$
does not depend on the chosen local frames $(e_1,\dots,e_n)$,
$(f_1,\dots,f_n)$. Therefore using \ref{local frames}, we may define
$\wedge^r g$ by this local prescription.
 
2. This can be proven along the same lines as 1. Alternatively, it
follows by applying 1. to the bf-morphism
$
(\M,\ \mu,\ \E'\to\F',\ \M\tensor\E'\leftarrow\F',\ n-r)
$
obtained from $g$ by setting $\E':=\F$ and $\F':=\M\tensor\E$.
\end{proof}

In the situation of the above proposition \ref{wedge g}, assume that
$\E=\E_1\oplus\E_2$ and $\F=\F_1\oplus\F_2$, where 
$\rk\E_i=\rk\F_i=:n_i$ for $i=1,2$. Assume furthermore that the
morphisms $\E\to\F$ and $\F\to\M\tensor\E$ both respect these
direct-sum decompositions and that there are $i_1,i_2\geq 0$ with
$i_1+i_2=i$,
such that the tupels
\begin{eqnarray*}
g_1 &:=& (\M,\ \mu,\ \E_1\to\F_1,\ \M\tensor\E_1\leftarrow\F_1,\ i_1) \\
\text{and}\quad
g_2 &:=& (\M,\ \mu,\ \E_2\to\F_2,\ \M\tensor\E_2\leftarrow\F_2,\ i_2)\quad,
\end{eqnarray*}
induced by $g$, are also bf-morphisms. We write $g=g_1\oplus
g_2$. The following lemma says that exterior
powers of bf-morphisms are compatible with direct sums whenever
this makes sense.
\begin{lemma}
\label{g=g1+g2}
Let $1\leq r\leq n$ and $r=r_1+r_2$ for some $r_1,r_2\geq 0$.
\begin{enumerate}
\item
If $\max(0,r-i)=\max(0,r_1-i_1)+\max(0,r_2-i_2)$, then we have for every
$\epsilon_1\in\Gamma(S,\wedge^{r_1}\E_1)$,\ 
$\epsilon_2\in\Gamma(S,\wedge^{r_2}\E_2)$ the following equality:
$$
(\wedge^r g)(\epsilon_1\wedge\epsilon_2)=
(\wedge^{r_1}g_1)(\epsilon_1)\wedge
(\wedge^{r_2}g_2)(\epsilon_2)\quad.
$$
\item
If $\min(i,r)=\min(i_1,r_1)+\min(i_2,r_2)$, then we have for every
$\omega_1\in\Gamma(S,\wedge^{r_1}\F_1)$,\ 
$\omega_2\in\Gamma(S,\wedge^{r_2}\F_2)$ the following equality:
$$
(\wedge^{-r} g)(\omega_1\wedge\omega_2)=
(\wedge^{-r_1}g_1)(\omega_1)\wedge
(\wedge^{-r_2}g_2)(\omega_2)\quad.
$$
\end{enumerate}
\end{lemma}
\begin{proof}
This follows immediately from the local description of $\wedge^r g$ and
$\wedge^{-r} g$ respectively.
\end{proof}
\begin{definition}
Let $S$ be a scheme, $\E$ and $\F$ two localy free $\Oo_S$-modules of
rank $n$ and
\begin{eqnarray*}
\Phi &=& (\text{$\Ll_i$, $\lambda_i$, $\M_i$, $\mu_i$,\quad 
          $\E_i\to\M_i\tensor\E_{i+1}$,\quad $\E_i\leftarrow\E_{i+1}$},\\
     & &\ \text{$\F_{i+1}\to\F_i$,\quad
          $\Ll_i\tensor\F_{i+1}\leftarrow\F_i$\quad
          $(0\leq i\leq n-1)$,\quad $\E_n\isomto\F_n$})\ ,
\end{eqnarray*}
a generalized isomorphism from $\E=\E_0$ to $\F=\F_0$. For 
$1\leq r\leq n$ we define the $r$-th exterior power
$$
\wedge^r\Phi :\quad
\Wedge^r\E\quad \longrightarrow\quad
\Tensor_{\nu=1}^r
\left(
\Tensor_{i=0}^{\nu-1}\Ll_i^{\dual}\tensor
\Tensor_{i=0}^{n-\nu}\M_i
\right)\tensor
\Wedge^r\F
$$
of $\Phi$ as the composition
$$
\wedge^r\Phi :=
(\wedge^{-r} g_0)\comp(\wedge^{-r} g_1)
\comp\dots\comp
(\wedge^{-r} g_{n-1})\comp(\wedge^r h_n)\comp(\wedge^r h_{n-1})
\comp\dots\comp
(\wedge^r h_0)\quad ,
$$
where $g_i$ and $h_i$ are the bf-morphisms
\begin{eqnarray*}
& &(\text{$\M_i$, $\mu_i$,\quad $\E_{i+1}\to\E_i$,\quad
          $\M_i\tensor\E_{i+1}\leftarrow\E_i$,\quad i})\\
\text{and}& &(\text{$\Ll_i$, $\lambda_i$,\quad $\F_{i+1}\to\F_i$,\quad
              $\Ll_i\tensor\F_{i+1}\leftarrow\F_i$,\quad i})\ 
\end{eqnarray*}
respectively for $0\leq i\leq n-1$, and where $h_n$ is the isomorphism 
$\E_n\isomto\F_n$.
\end{definition}

In the situation of the above definition consider especially the case,
where $\E=\oplus^n\Oo_S$ and $\F=\oplus^n\Ll$ for some invertible 
$\Oo_S$-module $\Ll$. Then we have natural direct sum decompositions
$$
\Wedge^r\E=
\bigoplus_{B} \Oo_S
\quad\text{and}\quad
\Wedge^r\F=
\bigoplus_{A}\Ll^r 
\quad,
$$
where $A$ and $B$ run through all subsets of cardinality $r$ of
$\{1,\dots,n\}$. For two such subsets $A$ and $B$, we denote by $\pi_A$ 
(respectively by $\iota_B$) the projection 
$\wedge^r\F\to\Ll^r$ 
onto the $A$-th component 
(respectively the inclusion 
$\Oo_S\injto\wedge^r\E$ 
of the $B$-th component).
Now we define
$$
\text{det}_{A,B}\Phi:=
\pi_A\comp(\wedge^r\Phi)\comp\iota_B:\quad
\Oo_S\quad\longrightarrow\quad
\Tensor_{\nu=1}^r
\left(
\Tensor_{i=0}^{\nu-1}\Ll_i^{\dual}\tensor
\Tensor_{i=0}^{n-\nu}\M_i
\right)\tensor\Ll^r
\quad.
$$
\begin{lemma}
\label{not vanishing}
Let $(\alpha,\beta)\in S_n\times S_n$ and
let $X(\alpha,\beta)$ be the open set of $\KGln$ defined in section
\ref{construction}.  Let $\Phi_{\text{univ}}$ be the generalized
isomorphism defined in \ref{Phi_univ}. Then the sections
$\text{det}_{\alpha[1,r],\beta[1,r]}\Phi_{\text{univ}}$ are nowhere
vanishing on $X(\alpha,\beta)$ for $1\leq r\leq n$.
\end{lemma}
\begin{proof}
From the proof of \ref{Phi_univ} it follows readily that the
restriction of $\text{det}_{\alpha[1,r],\beta[1,r]}\Phi_{\text{univ}}$
to $X(\alpha,\beta)$ is $\prod_{\nu=1}^r(t_\nu(\alpha,\beta)/t_0)$ as
an element of 
$$
\Gamma\left(
X(\alpha,\beta),\ 
\Oo\left(
\sum_{\nu=1}^r\left(
\sum_{i=0}^{n-\nu}Z_i - \sum_{i=0}^{\nu-1}Y_i
\right)\right)\right)
\subset\Gamma(X(\alpha,\beta),\ \K_{\KGln}) \quad.
$$
On the other hand, \ref{local description} tells us that 
$\prod_{\nu=1}^r(t_\nu(\alpha,\beta)/t_0)$ is a generator of
\linebreak[4] 
$\Oo(\sum_{\nu=1}^r(\sum_{i=0}^{n-\nu}Z_i - \sum_{i=0}^{\nu-1}Y_i))$
over $X(\alpha,\beta)$.
\end{proof}

\section{Proof of theorem \ref{KGln modular}}
\label{proof}

Let $S$ be a scheme, $\Ll$ an invertible $\Oo_S$-module. For 
$0\leq i\leq n-1$ let $(\Ll_i,\lambda_i)$, $(\M_i,\mu_i)$ be
invertible $\Oo_S$-modules together with global sections, such that
the zero sets of $\lambda_i$ and $\mu_j$ do not intersect for 
$i+j<n$. Given these data, we associate to every tupel 
$(\varphi_1,\dots,\varphi_n)$ of isomorphisms 
$$
\varphi_r:\quad
\Tensor_{i=0}^{n-r}\M_i^{\dual} 
\overset{\sim}{\longrightarrow}
\Tensor_{i=0}^{r-1}\Ll_i^{\dual}\tensor\Ll
\quad\quad
(1\leq r\leq n)
$$
the following generalized isomorphism
from $\oplus^n\Oo_S$ to $\oplus^n\Ll$\ :
\begin{eqnarray*}
\Phi(\varphi_1,\dots,\varphi_n)
     &:=& (\text{$\Ll_i$, $\lambda_i$, $\M_i$, $\mu_i$,\ 
          $\E_i\to\M_i\tensor\E_{i+1}$,\ $\E_i\leftarrow\E_{i+1}$},\\
     & &\ \text{$\F_{i+1}\to\F_i$,\
          $\Ll_i\tensor\F_{i+1}\leftarrow\F_i$\
          $(0\leq i\leq n-1)$,\ $\E_n\isomto\F_n$})\ ,
\end{eqnarray*}
where the locally free modules $\E_i$ and $\F_i$ are defined as
\begin{eqnarray*}
\E_i &:=&
\bigoplus_{j=1}^{n-i}
\left(
\Tensor_{\nu=0}^{i-1}\M_\nu^{\dual}
\right)
\oplus
\bigoplus_{j=n-i+1}^n
\left(
\Tensor_{\nu=0}^{n-j}\M_\nu^{\dual}
\right)\quad, \\
\F_i &:=&
\left(
\bigoplus_{j=1}^{i}
\left(
\Tensor_{\nu=0}^{j-1}\Ll_\nu^{\dual}
\right) \oplus 
\bigoplus_{j=i+1}^{n}
\left(
\Tensor_{\nu=0}^{i-1}\Ll_\nu^{\dual}
\right)
\right) \tensor \Ll
\quad,
\end{eqnarray*}
the morphisms
$$
\E_i\longrightarrow\M_i\tensor\E_{i+1}
\quad\text{and}\quad
\E_i\longleftarrow\E_{i+1}
$$
are described by the matrices
$$
\left[
\begin{array}{cc}
\eins_{n-i} & 0 \\
0 & \mu_i \eins_i
\end{array}
\right]
\quad\text{and}\quad
\left[
\begin{array}{cc}
\mu_i\eins_{n-i} & 0 \\
0 & \eins_i
\end{array}
\right]
\quad,
$$
the morphisms
$$
\F_{i+1}\longrightarrow\F_i
\quad\text{and}\quad
\Ll_i\tensor\F_{i+1}\longleftarrow\F_i
$$  
by the matrices
$$
\left[
\begin{array}{cc}
\eins_i & 0 \\
0 & \lambda_i \eins_{n-i}
\end{array}
\right]
\quad\text{and}\quad
\left[
\begin{array}{cc}
\lambda_i\eins_i & 0 \\
0 & \eins_{n-i}
\end{array}
\right]
$$
respectively, and the isomorphism $\E_n\isomto\F_n$ is given by the
diagonal matrix with entries $(\varphi_1,\dots,\varphi_n)$.
\begin{definition}
Let $S$ be a scheme, $\Ll$ an invertible $\Oo_S$-module and
\begin{eqnarray*}
\Phi &:=& (\text{$\Ll_i$, $\lambda_i$, $\M_i$, $\mu_i$,\quad 
          $\E_i\to\M_i\tensor\E_{i+1}$,\quad $\E_i\leftarrow\E_{i+1}$},\\
     & &\ \text{$\F_{i+1}\to\F_i$,\quad
          $\Ll_i\tensor\F_{i+1}\leftarrow\F_i$\quad
          $(0\leq i\leq n-1)$,\quad $\E_n\isomto\F_n$})
\end{eqnarray*}
an arbitrary generalized isomorphism from $\E_0=\oplus^n\Oo_S$ to
$\F_0=\oplus^n\Ll$. A {\em diagonalization of $\Phi$} with respect to
a pair $(\alpha,\beta)\in S_n\times S_n$ of permutations is a tupel
($u_i$,$v_i$\ $(0\leq i\leq n)$,\
$(\varphi_1,\dots,\varphi_n)$)
of isomorphisms
\begin{eqnarray*}
u_i&:&\quad \E_i  \Isomto
\bigoplus_{j=1}^{n-i}
\left(
\Tensor_{\nu=0}^{i-1}\M_\nu^{\dual}
\right)
\oplus
\bigoplus_{j=n-i+1}^n
\left(
\Tensor_{\nu=0}^{n-j}\M_\nu^{\dual}
\right)\quad\quad (0\leq i\leq n) 
\\
v_i&:&\quad \F_i \Isomto
\left(
\bigoplus_{j=1}^{i}
\left(
\Tensor_{\nu=0}^{j-1}\Ll_\nu^{\dual}
\right) \oplus 
\bigoplus_{j=i+1}^{n}
\left(
\Tensor_{\nu=0}^{i-1}\Ll_\nu^{\dual}
\right)
\right) \tensor \Ll
\quad\quad (0\leq i\leq n) 
\\
\varphi_r &:& \quad
\Tensor_{i=0}^{n-r}\M_i^{\dual} 
\overset{\sim}{\longrightarrow}
\Tensor_{i=0}^{r-1}\Ll_i^{\dual}\tensor\Ll
\quad\quad
(1\leq r\leq n)
\end{eqnarray*}
such that ($u_i$, $v_i$\ $(0\leq i\leq n)$)
establishes a quasi-equivalence between  
$\Phi$ and $\Phi(\varphi_1,\dots,\varphi_n)$
and such that
$$
u_n\cdot n_\beta:\
\E_0=\oplus^n\Oo_S \Isomto \oplus^n\Oo_S
\quad\text{and}\quad
v_n\cdot n_\alpha:\
\F_0=\oplus^n\Ll \Isomto \oplus^n\Ll
$$
are described by upper and lower triangular matrices respectively, with unit
diagonal entries.
\end{definition}
\begin{definition}
As in the above definition let $\Phi$ be a generalized isomorphism from 
$\oplus^n\Oo_S$ to $\oplus^n\Ll$. A pair
$(\alpha,\beta)\in S_n\times S_n$ 
of permutations is called 
{\em admissible}, if for all $1\leq r\leq n$ the global sections
$
\ \text{det}_{\alpha[1,r],\beta[1,r]}\Phi\ 
$
of
$$
\Tensor_{\nu=1}^r
\left(
\Tensor_{i=0}^{n-\nu}\M_i
\tensor
\Tensor_{i=0}^{\nu-1}\Ll_i^{\dual}
\right)
\tensor \Ll^r
$$
are nowhere vanishing on $S$.
\end{definition}
\begin{proposition}
\label{admissibility and diagonalizability}
Let $S$ be a scheme, $\Ll$ an invertible $\Oo_S$-module and 
$\Phi$ a generalized isomorphism from $\oplus^n\Oo_S$ to
$\oplus^n\Ll$. Then:
\begin{enumerate}
\item
For 
$(\alpha,\beta)\in S_n\times S_n$
the following are equivalent:
\begin{enumerate}
\item
there exists a diagonalization of $\Phi$ with respect to 
$(\alpha,\beta)$
\item
$(\alpha,\beta)$ is admissible for $\Phi$
\end{enumerate}
\item
Every point of $S$ has an open neighbourhood $U$, such that there is a
diagonalization of $\Phi|_U$ with respect to some pair
$(\alpha,\beta)\in S_n\times S_n$.
\item
For a given pair 
$(\alpha,\beta)\in S_n\times S_n$
there is at most one diagonalization of $\Phi$ with respect to 
$(\alpha,\beta)$.
\end{enumerate}
\end{proposition}
\begin{proof}
Let
\begin{eqnarray*}
\Phi &:=& (\text{$\Ll_i$, $\lambda_i$, $\M_i$, $\mu_i$,\quad 
          $\E_i\to\M_i\tensor\E_{i+1}$,\quad $\E_i\leftarrow\E_{i+1}$},\\
     & &\ \text{$\F_{i+1}\to\F_i$,\quad
          $\Ll_i\tensor\F_{i+1}\leftarrow\F_i$\quad
          $(0\leq i\leq n-1)$,\quad $\E_n\isomto\F_n$})\ ,
\end{eqnarray*}
be a generalized isomorphism from $\E_0=\oplus^n\Oo_S$ to
$\F_0=\oplus^n\Ll$. For $0\leq i\leq n-1$ denote bf-morphisms as follows:
\begin{eqnarray*}
g_i &:=& (\M_i,\ \mu_i,\ \E_{i+1}\to\E_i,\ 
          \M_i\tensor\E_{i+1}\leftarrow\E_i,\ i)\quad, 
\\
h_i &:=& (\Ll_i,\ \lambda_i,\ \F_{i+1}\to\F_i,\ 
          \Ll_i\tensor\F_{i+1}\leftarrow\F_i,\ i)\quad, 
\end{eqnarray*}
and let $h_n$ the isomorphism $\E_n\isomto\F_n$.

1. We may assume that $\alpha=\beta=\id$.
We show by induction on $n$ that admissibility of $(\id,\id)$ for
$\Phi$ implies the diagonalizability of $\Phi$ with respect to $(\id,\id)$.
The case $n=1$ is trivial, so assume $n\geq 2$. By assumption, the
morphism
$$
\text{det}_{\{1\},\{1\}}\Phi:\ \Oo_S\To
\Ll_0^{\dual}\tensor\Tensor_{i=0}^{n-1}\M_i\tensor\Ll
$$
is an isomorphism. Let
\begin{eqnarray*}
\label{tilde(E)'}
\tilde{\E}'_i &:=& 
\ker(
\pi_1\comp(\wedge^1h_0)\comp\dots\comp(\wedge^{-1}g_i):\
\E_i\To\Ll_0^{\dual}\tensor\Tensor_{j=i}^{n-1}\M_j\tensor\Ll
)
\\
& &\ \qquad\qquad\qquad\quad\qquad\qquad\qquad\qquad\qquad\qquad\qquad\qquad
(i\in[0,n])
\\
\tilde{\F}'_i &:=& 
\ker(
\pi_1\comp(\wedge^1h_0)\comp\dots\comp(\wedge^1h_{i-1}):\
\F_i\To\Ll_0^{\dual}\tensor\Ll
)\qquad
(i\in[1,n])
\\
\tilde{\F}'_0 &:=&
\ker(\pi_1:\oplus^n\Ll\To\Ll)=\oplus^{n-1}\Ll
\end{eqnarray*}
and
\begin{eqnarray*}
\tilde{\E}_i &:=&
\im(
(\wedge^{-1}g_{i-1})\comp\dots\comp(\wedge^{-1}g_0)\comp\iota_1:\
\Tensor_{j=0}^{i-1}\M_j^{\dual}\To\E_i
)\qquad
(i\in[0,n])
\\
\tilde{\F}_i &:=&
\im(
(\wedge^1h_i)\comp\dots\comp(\wedge^{-1}g_0)\comp\iota_1:\
\Tensor_{j=0}^{n-1}\M_j^{\dual}\To\F_i
)\qquad
(i\in[1,n])
\\
\tilde{\F}_0 &:=&
\im(
(\wedge^1h_0)\comp\dots\comp(\wedge^{-1}g_0)\comp\iota_1:\
\Ll_0^{\dual}\tensor\Tensor_{j=0}^{n-1}\M_j^{\dual}\To\F_0
)\qquad
\end{eqnarray*}
Then we have natural direct sum decompositions
\begin{eqnarray*}
\E_i &=& \tilde{\E}_i\oplus\tilde{\E}'_i \qquad (0\leq i\leq n)\quad, \\
\F_i &=& \tilde{\F}_i\oplus\tilde{\F}'_i \qquad (0\leq i\leq n)\quad.
\end{eqnarray*}
Since the bf-morphisms $g_i$ and $h_i$ respect these decompositions,
we can write 
$g_i=\tilde{g}_i\oplus\tilde{g}'_i$ and 
$h_i=\tilde{h}_i\oplus\tilde{h}'_i$
where $\tilde{g}_i$ 
(respectively $\tilde{g}'_i$, $\tilde{h}_i$, $\tilde{h}'_i$) is a
bf-morphism from $\tilde{\E}_{i+1}$ to $\tilde{\E}_i$
(respectively from $\tilde{\E}'_{i+1}$ to $\tilde{\E}'_i$,\quad
              from $\tilde{\F}_{i+1}$ to $\tilde{\F}_i$,\quad
              from $\tilde{\F}'_{i+1}$ to $\tilde{\F}'_i$)
for $0\leq i \leq n-1$. 
By the same reason, we can write $h_n=\tilde{h}_n\oplus\tilde{h}'_n$,
where $\tilde{h}_n:\tilde{\E}_n\isomto\tilde{\F}$ and 
$\tilde{h}'_n:\tilde{\E}'_n\isomto\tilde{\F}'$.
Observe that $\rk\tilde{g}_i=0$ and $\rk\tilde{g}'_i=i$ for
$0\leq i\leq n-1$ and that
$$
\rk\tilde{h}_i=
\left\{
\begin{array}{ll}
0\ ,\ &\text{if}\quad i=0 \\
1\ ,\ &\text{if}\quad i>0
\end{array}
\right. 
\quad,\quad
\rk\tilde{h}'_i=
\left\{
\begin{array}{ll}
0\ ,\ &\text{if}\quad i=0 \\
i-1\ ,\ &\text{if}\quad i>0
\end{array}
\right.
\quad.
$$
Now we define
\begin{eqnarray*}
\Ll' &:=& \Ll\tensor\Ll_0^{\dual}
\\
\Ll'_i &:=& \Ll_{i+1}\quad,\quad
\lambda'_i:=\lambda_{i+1}\qquad (0\leq i\leq n-2)
\\
\M'_i &:=& \M_i\quad,\quad
\mu'_i:=\mu_i\qquad (0\leq i\leq n-2)
\\
\E'_i &:=& \tilde{\E}'_i \qquad (0\leq i\leq n-1)
\\
\F'_i &:=& \tilde{\F}'_{i+1} \qquad (0\leq i\leq n-1)
\end{eqnarray*}
where we identify $\E'_0$ with $\oplus^{n-1}\Oo_S$ via the isomorphism
$$
\begin{CD}
\E'_0=\tilde{\E}'_0 @>\text{inclusion}>> \E_0=\oplus^n\Oo_S
@>\pi_{[2,n]}>> \oplus^{n-1}\Oo_S
\end{CD}
\quad,
$$
and $\F'_0$ with $\oplus^{n-1}\Ll'$ via the isomorphism
$$
\begin{CD}
\F'_0=\tilde{\F}'_1
@>\wedge^1\tilde{h}'_0>>
\Ll_0^{\dual}\tensor\oplus^{n-1}\Ll=\oplus^{n-1}\Ll'
\end{CD}
\quad.
$$
Let
\begin{eqnarray*}
\label{Phi'}
\Phi' &:=& (\text{$\Ll'_i$, $\lambda'_i$, $\M'_i$, $\mu'_i$,\quad 
          $\E'_i\to\M'_i\tensor\E'_{i+1}$,\quad $\E'_i\leftarrow\E'_{i+1}$},\\
      & &\ \text{$\F'_{i+1}\to\F'_i$,\quad
          $\Ll'_i\tensor\F'_{i+1}\leftarrow\F'_i$\quad
          $(0\leq i\leq n-2)$,\quad $\E'_{n-1}\isomto\F'_{n-1}$})\ ,
\end{eqnarray*}
where $\E'_{n-1}\isomto\F'_{n-1}$ is the composition
$$
\begin{CD}
\E'_{n-1}=\tilde{\E}'_{n-1}
@>\wedge^{-1}\tilde{g}'_{n-1}>>
\tilde{\E}'_n
@>\tilde{h}'_n>>
\tilde{\F}'_n=\F'_{n-1}
\end{CD}
\quad,
$$
and where the other morphisms are the ones from the $\tilde{g}'_i$ and
the $\tilde{h}'_i$. It is easy to see that $\Phi'$ is a generalized
isomorphism from $\oplus^{n-1}\Oo_S$ to
$\oplus^{n-1}\Ll'$. Furthermore, it follows from \ref{g=g1+g2} that
$$
(\wedge^r\Phi)(e_1\wedge\dots\wedge e_r)=
(\wedge^1\Phi)(e_1)\wedge
(\wedge^{r-1}\Phi')(e'_1\wedge\dots\wedge e'_{r-1})
\qquad (2\leq r\leq n)\ ,
$$
where $(e_1,\dots,e_n)\subset\Gamma(S,\E_0)$ and
$(e'_1,\dots,e'_{n-1})\subset\Gamma(S,\E'_0)$ are the 
canonical global frames of $\oplus^{n}\Oo_S$ and $\oplus^{n-1}\Oo_S$
respectively. Therefore we have
$$
\label{factorize det}
\text{det}_{[1.r][1,r]}\Phi=\text{det}_{\{1\}\{1\}}\Phi\tensor
\text{det}_{[1,r-1][1,r-1]}\Phi'
\qquad (2\leq r\leq n) \ .
$$
Since, by assumption, the sections $\text{det}_{[1.r][1,r]}\Phi$ are
nowhere vanishing, the above equation implies that the same is true also
for the sections $\text{det}_{[1,r-1][1,r-1]}\Phi'$\
$(2\leq r\leq n)$. In other words, $(\id,\id)$ is admissible for
$\Phi'$. By induction-hypothesis, we conclude that there exists a
diagonalization 
$
(u'_i,\ v'_i,\ (0\leq i\leq n-1),\
(\varphi'_1,\dots,\varphi'_{n-1}))
$
of $\Phi'$ with respect to $(\id,\id)$.

Let 
\begin{eqnarray*}
\tilde{u}'_i &:=& u'_i\qquad (0\leq i\leq n-1)\ ,\\
\tilde{v}'_i &:=& v'_{i-1}\qquad (1\leq i\leq n)\ ,
\end{eqnarray*}
and
\begin{eqnarray*}
\tilde{u}'_{n} &:&              
\begin{CD}
\tilde{\E}'_n 
@> \tilde{g}_{n-1}^{\sharp} > \sim >
\tilde{\E}'_{n-1}=\E'_{n-1}
@> u'_{n-1} > \sim >
\bigoplus_{j=2}^{n}\left(\Tensor_{\nu=0}^{n-j}\M^{\dual}\right)
\end{CD}
\quad,
\\
\tilde{v}'_0 &:&
\begin{CD}
\tilde{\F}'_0=\bigoplus^{n-1}\Ll=\Ll_0\tensor\F'_0
@> v'_0 > \sim >
\bigoplus^{n-1}\Ll
\end{CD}
\quad.
\end{eqnarray*}
Observe that there are natural isomorphisms
\begin{eqnarray*}
\tilde{u}_i &:&
\begin{CD}
\tilde{\E}_i
@> \sim >>
\Tensor_{j=0}^{i-1}\M_j^{\dual}
\end{CD}
\qquad (0\leq i\leq n)
\\
\tilde{v}_i &:&
\begin{CD}
\tilde{\F}_i
@> \sim >>
\Tensor_{j=0}^{n-1}\M_j^{\dual}
@> \text{det}_{\{1\}\{1\}}\Phi >>
\Ll\tensor\Ll_0^{\dual}
\end{CD}
\qquad (1\leq i\leq n)
\\
\tilde{v}_0 &:&
\begin{CD}
\tilde{\F}_0
@> \sim >>
\Ll_0\tensor\Tensor_{j=0}^{n-1}\M_j^{\dual}
@> \text{det}_{\{1\}\{1\}}\Phi >>
\Ll
\end{CD}
\end{eqnarray*}
We set
$u_i:=\tilde{u}_i\oplus\tilde{u}'_i$\ , 
$v_i:=\tilde{v}_i\oplus\tilde{v}'_i$ for $0\leq i\leq n$, and
$\varphi_r:=\varphi'_{r-1}$ for $2\leq r\leq n$.
Finally, we let
$\varphi_1:\Tensor_{i=0}^{n-1}\M_i^{\dual}\isomto\Ll\tensor\Ll_0^{\dual}$
be the isomorphism induced by $\text{det}_{\{1\}\{1\}}\Phi$.
It is now
easy to see that the tupel 
$(u_i,\ v_i\ (0\leq i\leq n),\ (\varphi_1,\dots,\varphi_n))$
is a diagonalization of $\Phi$ with respect to $(\id,\id)$.

Conversely, assume that there exists a diagonalization
$(u_i,\ v_i\ (0\leq i\leq n),\ (\varphi_1,\dots,\varphi_n))$
of $\Phi$ with respect to $(\id,\id)$. Observe that the diagram
$$
\xymatrix@C=5ex@R-1pc{
            &
\wedge^r\E_0 \ar[dd]^{\wedge^ru_0} \ar[rr]^-{\wedge^r\Phi} & &
\Nc_r\tensor(\wedge^r\F_0) \ar[dd]^{\wedge^rv_0} \ar[dr]^{\pi_{[1,r]}}&
\\
\Oo_S \ar[ur]^{\iota_{[1,r]}} \ar[dr]_{\iota_{[1,r]}}      &
            & 
            &
            &
\Nc_r\tensor\Ll^{\tensor r}            
\\
            &
\wedge^r(\oplus^n\Oo_S)
\ar[rr]^-{\wedge^r\Phi(\varphi_1,\dots,\varphi_n)} & &
\Nc_r\tensor(\wedge^r(\oplus^n\Ll)) \ar[ur]_{\pi_{[1,r]}}&  }      
$$
where 
$
\Nc_r:=
\Tensor_{\nu=1}^r
(\Tensor_{i=0}^{n-\nu}\M_i\tensor\Tensor_{i=0}^{\nu-1}\Ll_i^{\dual})
$,
is commutative for $1\leq r\leq n$. Therefore we may assume that
$\Phi=\Phi(\varphi_1,\dots,\varphi_n)$. But then
$\text{det}_{[1,r][1,r]}\Phi$ is the section induced by 
the isomorphism
$$
\varphi_1\tensor\dots\tensor\varphi_r:
\Tensor_{\nu=1}^r\Tensor_{i=0}^{n-\nu}\M_i^{\dual}
\Isomto
\Tensor_{\nu=1}^r\Tensor_{i=0}^{\nu-1}\Ll_i^{\dual}
$$
for $1\leq r\leq n$. In particular the $\text{det}_{[1,r][1,r]}\Phi$
are nowhere vanishing on $S$, which is precisely what is required for
the admissibility of $(\id,\id)$ for $\Phi$.

2. By 1, it suffices to show that in the case 
$S=\Spec k$\ ($k$ a field) there exists a pair 
$(\alpha,\beta)\in S_n\times S_n$
which is admissible for $\Phi$. We apply induction on $n$, the case $n=1$
being trivial.

It is an easy exercise in linear algebra, to show that the morphism 
$$
\wedge^1\Phi=
(h_0^{\flat})^{-1}\comp h_1^{\sharp}\dots
\comp h_{n-1}^{\sharp}\comp h_n\comp g_{n-1}^{\flat}\comp\dots
\comp g_0^{\flat}
$$ 
has at least rank one. Consequently there exist indices
$i_1,j_1\in\{1,\dots,n\}$,
such that the composition 
$
\text{det}_{\{i_1\}\{j_1\}}\Phi=
\pi_{i_1}\comp(\wedge^1\Phi)\comp\iota_{j_1}
$
is an isomorphism.
Let the sheaves $\tilde{\E}'_i$, $\tilde{\F}'_i$\ $(0\leq i\leq n)$ 
be defined as on page \pageref{tilde(E)'}, with $\pi_1$ replaced 
by $\pi_{i_1}$, and using these sheaves, let $\Phi'$ be defined as on 
page \pageref{Phi'}. This is a generalized isomophism from $k^{n-1}$
to $\oplus^{n-1}(\Ll_0^{\dual}\tensor\Ll)$. By induction-hypothesis
there exists a pair $(\alpha',\beta')\in S_{n-1}\times S_{n-1}$,
which is admissible for $\Phi'$.
Let $\alpha\in S_n$ be defined by
$$
\alpha(r) :=
\left\{
\begin{array}{ll}
i_1\ , &\text{if}\quad r=1\\
\alpha'(r-1)+1\ , &\text{if} \quad 2\leq r\leq n
\end{array}
\right. \qquad,
$$
and let $\beta\in S_n$ be defined analogously.
As on page \pageref{factorize det} we have
$$
\text{det}_{\alpha[1.r],\beta[1,r]}\Phi=
\text{det}_{\{i_1\}\{j_1\}}\Phi\tensor
\text{det}_{\alpha'[1,r-1],\beta'[1,r-1]}\Phi'
\qquad (2\leq r\leq n) 
$$
for $2\leq r\leq n$,\ i.e. the pair $(\alpha,\beta)$ is admissible
for $\Phi$.

3. This follows from the proof of 1, since it is clear that the
construction of the diagonalization there is unique.
\end{proof}
\begin{proposition}
\label{X(alpha,beta) modular}
Let $S$ be a scheme and $\Phi$ a generalized isomorphism from
$\oplus^n\Oo_S$ to itself.
\begin{enumerate}
\item
If $(\alpha,\beta)\in S_n\times S_n$ is admissible for $\Phi$, then
there exists a unique morphism $S\to X(\alpha,\beta)$, such that the
pull-back of $\Phi_{\text{univ}}$ to $S$ by this morphism 
is equivalent to $\Phi$.
\item
\label{description of X(alpha,beta)}
We have the following description of $X(\alpha,\beta)$ as an open subset of
$\KGln$:
$$
X(\alpha,\beta)=
\left\{
x\in\KGln\ |\
(\text{det}_{\alpha[1,r],\beta[1,r]}\Phi_{\text{univ}})(x)\neq 0\quad
\text{for}\ 1\leq r\leq n 
\right\} 
\quad.
$$
\item
If $(\alpha',\beta')\in S_n\times S_n$ is a further admissible pair
for $\Phi$, then the above morphism $S\to X(\alpha,\beta)$ factorizes
over the inclusion 
$X(\alpha,\beta)\cap X(\alpha',\beta') \injto X(\alpha,\beta)$.
\end{enumerate}
\end{proposition}
\begin{proof}
1. Let
\begin{eqnarray*}
\Phi &=& (\text{$\Ll_i$, $\lambda_i$, $\M_i$, $\mu_i$,\quad 
          $\E_i\to\M_i\tensor\E_{i+1}$,\quad $\E_i\leftarrow\E_{i+1}$},\\
     & &\ \text{$\F_{i+1}\to\F_i$,\quad
          $\Ll_i\tensor\F_{i+1}\leftarrow\F_i$\quad
          $(0\leq i\leq n-1)$,\quad $\E_n\isomto\F_n$})\ .
\end{eqnarray*}
By proposition \ref{admissibility and diagonalizability} there exists a
diagonalization $(u_i,\ v_i\ (0\leq i\leq n),\ (\varphi_1,\dots,\varphi_n))$
of $\Phi$ with respect to $(\alpha,\beta)$. Let
$a_{ji}\in\Gamma(S,\Oo_S)$ (respectively $b_{ij}\in\Gamma(S,\Oo_S)$)
($1\leq i<j\leq n$) be the nontrivial entries of the lower
(respectively upper) triangular matrix $(v_0\cdot n_\alpha)^{-1}$
(respectively $u_0\cdot n_\beta$). Let $a:S\to U^-$ and $b:S\to U^+$
be the morphisms defined by $a^*(y_{ji})= a_{ji}$ and
$b^*(z_{ij})=b_{ij}$ respectively. Furthermore, let 
$\varphi:S\to\T$ be the morphism induced by the tupel 
$$
(\Ll_i,\ \lambda_i,\ \M_i,\ \mu_i\ (0\leq i\leq n-1), 
\varphi_r\ (1\leq r\leq n))  \quad,
$$
(cf lemma \ref{T modular}). Thus we have a morphism
$$
f: S \overset{(a,\varphi,b)}{\longrightarrow} U^-\times\T\times U^+
     \Isomto X(\alpha,\beta)
\quad,
$$
where the right isomorphism is the inverse of the one in lemma \ref{ytz}.
It is clear that
\begin{eqnarray*}
f^*\Oo(Y_i(\alpha,\beta)) \isomorph \Ll_i\ &,& 
f^* \one_{\Oo(Y_i(\alpha,\beta))} = \lambda_i\, 
\\
f^*\Oo(Z_i(\alpha,\beta)) \isomorph \M_i\ &,& 
f^* \one_{\Oo(Z_i(\alpha,\beta))} = \mu_i\qquad (0\leq i\leq n-1)\ .
\end{eqnarray*} 
Denote by 
$(u'_i,\ v'_i\ (0\leq i\leq n),\ (\varphi'_1,\dots,\varphi'_n))$ the
pull-back under $f$ of the diagonalization of
$\Phi_{\text{univ}}|_{X(\alpha,\beta)}$, 
which exists by \ref{not vanishing} and 
\ref{admissibility and diagonalizability}.
By the uniquness of diagonalizations 
(cf \ref{admissibility and diagonalizability}), we have $u_0=u'_0$,\ 
$v_0=v'_0$ and 
$(\varphi_1,\dots,\varphi_n)=(\varphi'_1,\dots,\varphi'_n)$.
Therefore the isomorphisms
\begin{eqnarray*}
(u'_i)^{-1}\comp u_i &:& \E_i \Isomto f^*E_i \\
(v'_i)^{-1}\comp v_i &:& \F_i \Isomto f^*F_i 
\end{eqnarray*}
induce an equivalence between $\Phi$ and $f^*\Phi_{\text{univ}}$. This
proves the existence part of the proposition.
For uniqueness, assume that $\tilde{f}$ is a further morphism from $S$ to
$X(\alpha,\beta)$, such that $\Phi$ is equivalent to 
$\tilde{f}^*\Phi_{\text{univ}}$. 
Let 
$\overline{u}_i:\E_i \isomto \tilde{f}^*E_i$,\ 
$\overline{v}_i:\F_i \isomto \tilde{f}^*F_i$,\ 
$(0\leq i\leq n)$
be an equivalence. Note that by definition 
$\overline{u}_0=\text{id}_{\oplus^n\Oo_S}=\overline{u}_0$.
Let 
$(\tilde{u}_i,\ \tilde{v}_i\ (0\leq i\leq n), 
 (\tilde{\varphi}_1,\dots,\tilde{\varphi}_n))$ 
be the pull-back under $\tilde{f}$ of the diagonalization with respect
to $(\alpha,\beta)$ of 
$\Phi_{\text{univ}}|_{X(\alpha,\beta)}$.
Then
$(\tilde{u}_i\comp\overline{u}_i,\ 
  \tilde{v}_i\comp\overline{v}_i\ (0\leq i\leq n), 
 (\tilde{\varphi}_1,\dots,\tilde{\varphi}_n))$ 
is a diagonalization of $\Phi$ with respect to $(\alpha,\beta)$.
By \ref{admissibility and diagonalizability}.3 we conclude that
$\tilde{u}_0=\tilde{u}_0\comp\overline{u}_0=u_0$,\quad
$\tilde{v}_0=\tilde{v}_0\comp\overline{v}_0=v_0$ and
$(\tilde{\varphi}_1,\dots,\tilde{\varphi}_n)=(\varphi_1,\dots,\varphi_n)$.
But this implies that the composite morphism
$$
S\overset{\tilde{f}}{\longrightarrow}
X(\alpha,\beta)\Isomto U^-\times\T\times U^+
$$
equals $(a,\varphi,b)$ and thus that $\tilde{f}=f$.

2. Denote for a moment by $U$ the open subset of $\KGln$, defined by
the nonvanishing of
$\text{det}_{\alpha[1,r],\beta[1,r]}\Phi_{\text{univ}}$ for
$1\leq r\leq n$. We have already seen in \ref{not vanishing} that
$X(\alpha,\beta)$ is contained in $U$. Let $x\in U$. Since
$X(\alpha,\beta)$ is dense in $U$, there exists a generalization 
$y\in X(\alpha,\beta)$ of $x$. Then there exists a morphism 
$f:S\to U$, where $S$ is the Spec of a valuation ring, such that the
special point of $S$ is mapped to $x$ and its generic point to $y$.
By definition of $U$, the pair $(\alpha,\beta)$ is admissible for the 
generalized isomorphism $f^*\Phi_{\text{univ}}$. Therefore 1
tells us that there exists a morphism $f':S\to X(\alpha,\beta)$, which
coincides with $f$ at the generic point of $S$. Since $\KGln$ is
separabel, it follows that $f=f'$ and thus that $x\in X(\alpha,\beta)$.

3. This follows immediatelly from \ref{description of X(alpha,beta)}.
\end{proof}
\begin{proof}(Of theorem \ref{KGln modular}).
Let $S$ be a scheme and $\Phi$ a generalized isomorphism from
$\oplus^n\Oo_S$ to itself. By proposition 
\ref{admissibility and diagonalizability}, there is a covering of $S$
by open sets $U_i$ $(i\in I)$, and for every index $i\in I$ a pair
$(\alpha_i,\beta_i)\in S_n\times S_n$, which is admissible for
$\Phi|_{U_i}$. Proposition  \ref{X(alpha,beta) modular}.1 now
tells us that there exists for each 
$i\in I$ a unique morhpism 
$f_i:U_i\to X(\alpha_i,\beta_i)$ with the property
that there is an equivalence, say $u_i$,
from $\Phi|_{U_i}$ to $f_i^*\Phi_{\text{univ}}$. 
By proposition \ref{X(alpha,beta) modular}.2, the $f_i$ glue
together, to give a morphism $f:S\to\KGln$. It remains to show that
also the $u_i$ glue together, to give an overall equivalence from
$\Phi$ to $f^*\Phi_{\text{univ}}$. For this, it suffices to show that
for two generalized isomorphisms $\Phi$ and $\Phi'$ from 
$\oplus^n\Oo_S$ to itself there exists at most one equivalence from 
$\Phi$ to $\Phi'$. The question being local, we may assume by
proposition \ref{admissibility and diagonalizability}.2 that
$\Phi'$ is diagonalizable with respect to some pair
$(\alpha,\beta)\in S_n\times S_n$. Composing the diagonalization of
$\Phi'$ with any equivalence from $\Phi$ to $\Phi'$ gives a
diagonalization of $\Phi$ with respect to $(\alpha,\beta)$. Since
different equivalences from $\Phi$ to $\Phi'$ would yield different
diagonalizations of $\Phi$, proposition 
\ref{admissibility and diagonalizability}.3 tells us that there exists
at most one equivalence.
\end{proof}

\section{Complete collineations}
\label{PGln}

In this section we prove a modular property for the compactification
$\PGlnb$ of $\PGln$ and compare it with the results of other 
authors. 

The scheme $\PGlnb$ together with closed subschemes
$\overline{\Delta}_r$ ($1\leq r\leq n-1$) is defined by successive
blow ups as follows.
Let
$\overline{\Omega}^{(0)}:=\Proj(\Z[x_{i,j}\ (1\leq i,j\leq n)])$
and  let 
$\overline{\Delta}_r^{(0)}:=
V^+((\det_{AB}(x_{ij})\ |\ A,B\subseteq\{1,\dots,n\} 
                                     ,\  \card A=\card B=r+1))$
($1\leq r\leq n-1$).
Inductively, define  
$\overline{\Omega}^{(\nu)}$
as the blowing up of $\overline{\Omega}^{(\nu-1)}$
along  $\overline{\Delta}^{(\nu-1)}_\nu$.
The closed subscheme
$\overline{\Delta}^{(\nu)}_r\subset\overline{\Omega}^{(\nu)}$ is by definition the strict
(resp. total) transform of $\overline{\Delta}^{(\nu-1)}_r$ for $r\neq\nu$
(resp. $r=\nu$). 
By definition, $\PGlnb:=\overline{\Omega}^{(n-1)}$ and 
$\overline{\Delta}_r:=\overline{\Delta}^{(n-1)}_r$ for $1\leq r\leq n-1$.

The variety $\PGlnb\times\Spec(\C)$ is the
so-called ``wonderful compactification'' of the homogenuos space 
$\PGlnC=(\PGlnC\times\PGlnC)/\PGlnC$ 
(cf. \cite{CP}). Vainsencher
\cite{Vainsencher}, Laksov \cite{Laksov2} and Thorup-Kleiman
\cite{Thorup&Kleiman} have given a modular description for (some of)
the $S$-valued points of $\PGlnb$. We will give a brief account
of their results.

Let $R\subseteq[1,n-1]$ and let $S$ be a scheme.
Following the terminology of Vainsencher, an $S$-valued {\em complete
collineation of type $R$} from a rank-$n$ vector bundle $E$ to a rank-$n$
vector bundle $F$ is 
a collection of morphisms
$$
v_i:E_i\to\Nc_i\tensor F_i \quad (0\leq i\leq k),
$$
where $R=\{r_1,\dots,r_k\}$, $0=:r_0<r_1<\dots<r_k<r_{k+1}:=n$,
the $\Nc_i$ are line bundles,
the $E_i$, $F_i$ are vector bundles on S and
$v_i$ has overall rank $r_{i+1}-r_i$; furthermore it is required that
$E_0=E$, $F_0=F$, 
and
$E_i=\ker(v_{i-1})$, $F_i=\Nc_{i-1}^{\dual}\tensor\coker(v_{i-1})$ 
for $1\leq i\leq k$. Vainsencher proves that the locally closed subscheme
$(\cap_{r\in R} \overline{\Delta}_r)\setminus
 \cup_{r\not\in R}\overline{\Delta}_r$ of $\PGlnb$
represents the functor which to each scheme $S$ associates the set of 
isomorphism classes of $S$-valued complete collineations of type $R$
from $\oplus^n\Oo_S$ to itself.

Laksov went further. He succeeded to give a modular description for those
$S$-valued points of 
$\overline{\Delta}(R):=\cap_{r\in R} \overline{\Delta}_r$ for which
the pull-back of the divisor 
$\sum_{r\not\in R}\overline{\Delta}_r|_{\overline{\Delta}(R)}$
on $\overline{\Delta}(R)$ is a well-defined divisor on $S$.
We refer the reader to \cite{Laksov2} for more details.

Finally, Thorup and Kleiman gave the following description for {\em all}
$S$-valued points of $\PGlnb$.
A morphism $u$ from $(\oplus^n\Oo_S)\tensor(\oplus^n\Oo_S)^{\dual}$ to a line
bundle $\Ll$ is called a {\em divisorial form}, if for each $i\in[1,n]$
the image $\M_i(u)$ of the induced map
$\wedge^i(\oplus^n\Oo_S)\tensor\wedge^i(\oplus^n\Oo_S)^{\dual}\to
 \Ll^{\tensor i}$
is an invertible sheaf. In this case denote by $u^i$ the induced
surjection 
$\wedge^i(\oplus^n\Oo_S)\tensor\wedge^i(\oplus^n\Oo_S)^{\dual}\to
 \M_i(u)$. 
Following Thorup and Kleiman, we  define a 
{\em projectively complete bilinear form} as a tupel 
$\ub=(u_1,\dots,u_n)$, where 
$u_i:\wedge^i(\oplus^n\Oo_S)\tensor\wedge^i(\oplus^n\Oo_S)^{\dual}\to\M_i$ 
is a surjection onto an invertible sheaf $\M_i$
for $1\leq i\leq n$, such that 
$\ub$ is ``locally the pull-back of a divisorial form''.
The last requirement means the following: For each point $x\in S$, there 
exists an open neighborhood $U$ of $x$, a morphism from $U$ to some
scheme $S'$ and a divisorial form 
$u:(\oplus^n\Oo_{S'})\tensor(\oplus^n\Oo_{S'})^{\dual}\to\Ll'$ on $S'$ 
such that the restriction of $\ub$ to $U$ is isomorphic to the 
pull-back of $(u^1,\dots,u^n)$. Thorup and Kleiman show that
$\PGlnb$ represents the functor that to each scheme 
$S$ associates the
set of isomorphism classes of projectively complete bilinear forms
on $S$.

None of these descriptions is completely satisfactory: Those of
Vainsencher and Laksov deal only with special $S$-valued points and
the description of Thorup-Kleiman is not explicit and is not 
truely modular, since the
condition ``to be locally pull-back of a divisorial form'' makes
reference to the existence of morphisms between schemes.

The terminology in the following definition will be justified by the
corollary \ref{corollary to KGln modular} below.

\begin{definition}
\label{complete collineations}
Let $S$ be a scheme and $\E$, $\F$ two locally free $\Oo_S$-modules 
of rank $n$. A {\em complete collineation} from $\E$ to $\F$ is a tupel
$$
\Psi=(\Ll_i,\ \lambda_i,\ \F_{i+1}\to\F_i,\ 
      \Ll_i\tensor\F_{i+1}\leftarrow\F_i\ (0\leq i\leq n-1))
\quad ,
$$
where $\E=\F_n,\F_{n-1},\dots,\F_1,\F_0=\F$ 
are locally free $\Oo_S$-modules of rank $n$, the tupels
$$
(\Ll_i,\ \lambda_i,\ \F_{i+1}\to\F_i,\ 
      \Ll_i\tensor\F_{i+1}\leftarrow\F_i,\ i)
$$
are bf-morphisms of rank $i$ for $0\leq i\leq n-1$ 
and $\lambda_0=0$, 
such that for each point $x\in S$ and index
$i\in\{0,\dots,n-1\}$ with the property that $\lambda_i(x)=0$, 
the following holds:

If $(f,g)$ is one of the following two pairs of morphisms:
\begin{eqnarray*}
&  \F_n[x]
   \overset{f}{\longrightarrow}
   \F_{i+1}[x]
   \overset{g}{\longrightarrow}
   \F_{i}[x] \ , &\\
&  ((\tensor_{j=0}^{i}\Ll_j)\tensor\F_{i+1})[x]
   \overset{g}{\longleftarrow}
   ((\tensor_{j=0}^{i-1}\Ll_j)\tensor\F_{i})[x]
   \overset{f}{\longleftarrow}
   \F_0[x] \ ,&
\end{eqnarray*}
then $\im(g\comp f)=\im(g)$.

Two complete collineations $\Psi$ and $\Psi'$  
from $\E$ to $\F$ are called {\em equivalent}, 
if there are isomorphisms $\Ll_i\isomto\Ll'_i$, $\F_i\isomto\F'_i$, such that
all the obvious diagrams commute and such that $\F_n\isomto\F'_n$ and
$\F_0\isomto\F'_0$ 
is the identity on $\E$ and $\F$.
\end{definition}

\begin{corollary} 
\label{corollary to KGln modular}
On $\PGlnb$ there exists a
universal complete collineation 
$\Psi_{\text{univ}}$ 
from $\oplus^n\Oo$ to itself, 
such that the pair
$(\PGlnb,\Psi_{\text{univ}})$
represents the functor, which to every scheme $S$ associates the set of
equivalence classes of complete
collineations from $\oplus^n\Oo_S$ to itself. 
\end{corollary}

\begin{proof}
Observe that $\PGlnb$ is naturally
isomorphic to the closed subscheme $Y_0$ of $\KGln$.
The restriction of $\Phi_{\text{univ}}$ to 
$\PGlnb$ induces in an obvious way a complete 
collineation $\Psi_{\text{univ}}$ 
of $\oplus^n\Oo$ to itself 
on $\PGlnb$.
The corollary now follows from theorem \ref{KGln modular}.
\end{proof}

We conclude this section by indicating how one can recover
Vainsencher's and Thorup-Kleiman's description from corollary
\ref{corollary to KGln modular}. Let $S$ be a scheme and let
$$
{\Psi}=(\Ll_i,\ \lambda_i,\ \F_{i+1}\to\F_i,\ 
      \Ll_i\tensor\F_{i+1}\leftarrow\F_i\ (0\leq i\leq n-1))
$$
be a complete collineation from $\oplus^n\Oo_S$ to itself
in the sense of definition \ref{complete collineations}.

First assume that there exists a subset $R$ of $[1,n-1]$, such that
the map $S\to\PGlnb$ corresponding to ${\Psi}$ factors
through 
$(\cap_{r\in R} \overline{\Delta}_r)\setminus
 \cup_{r\not\in R}\overline{\Delta}_r$.
This means that $\lambda_r$ is zero for $r\in R$ and is nowhere
vanishing for $r\in[1,n-1]\setminus R$.
As above, let $R=\{r_1,\dots,r_k\}$,
$0=:r_0<r_1<\dots<r_k<r_{k+1}:=n$.
For $0\leq i\leq k$ let
\begin{eqnarray*}
E_i &:=& \ker(\F_n\to\F_{n-1}\to\dots\to\F_{r_i}) \\
F_i &:=& \Nc_i^{\dual}\tensor
         \ker(\F_{r_{i+1}}\to\F_{r_{i+1}-1}\to\dots\to\F_{r_i})
\quad ,
\end{eqnarray*}
where $\Nc_i:=\tensor_{j=1}^i\Ll_{r_i}^{\dual}$.
Observe that the data in ${\Psi}$ provide natural maps
$$
v_i: E_i\to \Nc_i\tensor F_i 
$$
of overall rank $r_{i+1}-r_i$ for $0\leq i\leq k$.
Furthermore we have natural isomorphisms
$E_0=\F_n=\oplus^n\Oo$, $F_0\isomorph\F_0=\oplus^n\Oo$, and for 
$1\leq i\leq k$:
\begin{eqnarray*}
\ker(v_{i-1}) &=& E_{i} \\
\coker(v_{i-1}) &\isomorph& \coker(\F_n\to\F_{r_{i}}) =
                \coker(\F_{r_{i+1}}\to\F_{r_{i}}) \isomorph   \\
            &\isomorph& \ker(\Ll_{r_{i}}\tensor\F_{r_{i+1}}\to
                     \Ll_{r_{i}}\tensor\F_{r_{i}}) \isomorph
                \Nc_{i-1}\tensor F_{i}
\quad .
\end{eqnarray*}
Thus, $(v_i)_{0\leq i\leq k}$ is a complete homomorphism
of type $R$ in the sense of Vainsencher.

Now let ${\Psi}$ be arbitrary.
As in section \ref{exterior}, ${\Psi}$ induces
nowhere vanishing morphisms 
$$
\wedge^r{\Psi}:
\wedge^r\F_n \to
\left( 
\Tensor_{\nu=1}^r\Tensor_{i=0}^{\nu-1}\Ll_i^{\dual} 
\right)
\tensor \wedge^r\F_0
$$
and thus surjections
$$
u_r:\wedge^r\F_n\tensor\wedge^r\F_0^{\dual}\to
    \Tensor_{\nu=1}^r\Tensor_{i=0}^{\nu-1}\Ll_i^{\dual}
$$
for $1\leq r\leq n$.
The tupel $(u_r)_{1\leq r\leq n}$ is a projectively complete bilinear
form in the sense of Thorup-Kleiman. This follows from the fact that 
$\wedge^1\Psi_{\text{univ}}$ induces a divisorial form
on $\PGlnb$.

\section{Geometry of the strata}
\label{geometry of the strata}

In this section we need relative versions of the
varieties $\KGln$, $\PGlnb$ and
$\overline{O}_{I,J}:=\cap_{i\in I}Z_i\cap\cap_{j\in J}Y_j$, where
$I$ and $J$ are subsets of $[0,n-1]$.
They are defined in the following theorem.
\begin{theorem}
Let $T$ be a scheme and let $\E$ and $\F$ be two locally free
$\Oo_T$-modules of rank $n$. For a $T$-scheme $S$ we write
$\E_S$ and $\F_S$ for the pull-back of $E$ and $F$ to $S$.
Let $I$, $J$ be two subsets of $[0,n-1]$ 
Consider the following contravariant functors from the 
category of $T$-schemes to the category of sets:
\begin{eqnarray*}
\KGlbf(\E,\F):\ S &\mapsto&
\left\{
\begin{array}{lll}
\text{equivalence classes of}\\
\text{generalized isomorphisms}\\
\text{from $\E_S$ to $\F_S$}
\end{array}
\right\} 
\\
\PGlbf(\E,\F):\ S &\mapsto&
\left\{
\begin{array}{lll}
\text{equivalence classes of}\\
\text{complete collineations}\\
\text{from $\E_S$ to $\F_S$}
\end{array}
\right\}
\\
{\bf \overline{O}}_{I,J}(\E,\F):\ S &\mapsto&
\left\{
\begin{array}{lllll}
\text{equivalence classes of}\\
\text{generalized isomorphisms $\Phi$}\\
\text{as in definition \ref{generalized isomorphism} from}\\
\text{$\E_S$ to $\F_S$, with $\mu_i=0$ for $i\in I$}\\
\text{and $\lambda_j=0$ for $j\in J$}
\end{array}
\right\}
\end{eqnarray*}
These functors are representable by smooth projective $T$-schemes,
which we will call
$\KGl(\E,\F)$, $\PGlb(\E,\F)$ and $\overline{O}_{I,J}(\E,\F)$
respectively.
\end{theorem}

\begin{proof}
In the case of $T=\Spec \Z$ and $\E=\F=\oplus^n\Oo_{\Spec \Z}$,
the theorem is a consequence of \ref{KGln modular} and
\ref{corollary to KGln modular}, where the representing objects
are $\KGln$, $\PGlnb$ and $\overline{O}_{I,J}$ respectively.
Let $T=\cup U_i$ an open covering such that there exist trivializations
\begin{eqnarray*}
\xi_i &:& \E|_{U_i} \isomto \oplus^n \Oo_{U_i}\\
\zeta_i &:& F|_{U_i} \isomto \oplus^n \Oo_{U_i}
\quad .
\end{eqnarray*}
Let $\KGl_{U_i}:=\KGl_n\times_{\Spec\Z}U_i$ and 
$\pi_i:\KGl_{U_i}\to U_i$ the projection.
By corollary \ref{operation},
over the intersections $U_i\cap U_j$ the pairs 
$(\xi_i\xi_j^{-1},\zeta_i\zeta_j^{-1})$ induce isomorphisms
$\pi_i^{-1}(U_i\cap U_j)\isomto\pi_j^{-1}(U_i\cap U_j)$.
These provide the data for the pieces $\KGl_{U_i}$
to glue together to define $\KGl(\E,\F)$.
Using theorem \ref{KGln modular} it is easy to check that
$\KGl(\E,\F)$ has the required universal property.
This proves the existence of $\KGl(\E,\F)$. 
The existence of $\PGlb(\E,\F)$ and of $\overline{O}_{I,J}(\E,\F)$ is proved 
analogously.
\end{proof}

\begin{definition}
Let $T$ be a scheme and $\E$ a locally free $\Oo_T$-module
of rank $n$. Let ${\bf d}:=(d_0,\dots,d_t)$, where
$0\leq d_0\leq\dots\leq d_{t}\leq n$
Let $\Fl^{\bf d}(\E)$ be the flag variety which represents 
the following contravariant functor from the category of
$T$-schemes to the category of sets:
$$
S\mapsto
\left\{
\begin{array}{lll}
\text{All filtrations $F_0\E\subseteq\dots\subseteq F_t\E$, where} \\
\text{$F_p\E$ is a subbundle of rank $d_p$ of $\E_S$}\\
\text{for $0\leq p\leq t$}
\end{array}
\right\}
$$
Here as usual, a subbundle of $\E_S$ means a locally free subsheaf
of $\E_S$, which is locally a direct summand.
\end{definition}

After these preliminaries we can state the main result of this section, 
which descibes the structure of the schemes $\overline{O}_{I,J}$ defined above.

\begin{theorem}
\label{strata}
Let $T$ be a scheme and let $\E$ and $\F$ be two locally free 
$\Oo_T$-modules of rank $n$. Let
$I:=\{i_1,\dots,i_r\}$ and $J:=\{j_1,\dots,j_s\}$,
where $i_1+j_1\geq n$ and
$0\leq i_1<\dots<i_{r+1}:=n$, $0\leq j_1<\dots<j_{s+1}:=n$.
Let
${\bf d}:=(d_0,\dots,d_{r+s+1})$ and
$\deltabf:=(\delta_0,\dots,\delta_{r+s+1})$,
where
$$
d_p:=
\left\{
\begin{array}{ll}
\text{$n-j_{s+1-p}$\quad\  for\quad  $0\leq p\leq s$}\\
\text{$i_{p-s}$\qquad\qquad for\quad  $s+1\leq p\leq r+s+1$}
\end{array}
\right.
$$
and $\delta_q:=n-d_{r+s+1-q}$ for $0\leq q\leq r+s+1$.
Let 
\begin{eqnarray*}
&
0=U_0\subseteq U_1\subseteq\dots\subseteq U_{r+s+1}=
\E_{\Fl^{\bf d}(\E)\times\Fl^{\deltabf}(\F)}
& \\
\text{and}\quad &
0=V_0\subseteq V_1\subseteq\dots\subseteq V_{r+s+1}=
\F_{\Fl^{\bf d}(\E)\times\Fl^{\deltabf}(\F)}
\end{eqnarray*}
be the pull back to $\Fl^{\bf d}(\E)\times\Fl^{\deltabf}(\F)$ 
of the universal flag on 
$\Fl^{\bf d}(\E)$ and $\Fl^{\deltabf}(\F)$
respectively.
Then there is a natural isomorphism
$$
\overline{O}_{I,J}(\E,\F)\isomto
P_1
\underset{\Fl}{\times}
\dots
\underset{\Fl}{\times}
P_r
\underset{\Fl}{\times}
Q_s
\underset{\Fl}{\times}
\dots
\underset{\Fl}{\times}
Q_1
\underset{\Fl}{\times}
K'
\quad,
$$
where $\Fl:=\Fl^{\bf d}(\E)\times\Fl^{\deltabf}(\F)$ and where
\begin{eqnarray*}
P_p &:=& \PGlb(V_{r-p+1}/V_{r-p},U_{s+p+1}/U_{s+p})
\quad (1\leq p\leq r)\\
Q_q &:=& \PGlb(U_{s-q+1}/U_{s-q},V_{r+q+1}/V_{r+q})
\quad (1\leq q\leq s)\\
K' &:=& \KGl(U_{s+1}/U_{s},V_{r+1}/V_{r})
\quad.
\end{eqnarray*}
\end{theorem}

\begin{proof}
The isomorphism 
$$
\overline{O}_{I,J}\isomorph
P_1
\underset{\Fl}{\times}
\dots
\underset{\Fl}{\times}
P_r
\underset{\Fl}{\times}
Q_s
\underset{\Fl}{\times}
\dots
\underset{\Fl}{\times}
Q_1
\underset{\Fl}{\times}
K'
$$
is given on $S$-valued points by the bijectiv correspondence
$$
\Phi\longleftrightarrow
((F_\bullet\E,F_\bullet\F),
\varphi_1,\dots,\varphi_r,\psi_s,\dots,\psi_1,\Phi')
\quad,
$$
where
{\footnotesize
$$
\Phi = 
\left(
\xymatrix@C-4ex{
\E_{S}
\ar@/^1.2pc/|{\tensor}[rr]
\ar @/_0.65pc/ @{} [rr]|{(\M_0,\mu_0)}
& &
\E_1 
\ar[ll]_0
\ar @/_0.65pc/ @{} [rr]|{(\M_1,\mu_1)}
\ar@/^1.2pc/|{\tensor}[rr]
& &
\E_2
\ar[ll]_1
& 
\dots
& 
\E_{n-1}
\ar @/_0.65pc/ @{} [rr]|{(\M_{n-1},\mu_{n-1})}
\ar@/^1.2pc/|{\tensor}[rr]
& & 
\E_n
\ar[ll]_{n-1}
\ar[rr]^\sim
& & 
\F_n
\ar[rr]^{n-1}
& & 
\F_{n-1}
\ar@/_1.2pc/|{\tensor}[ll]
\ar @/^0.65pc/ @{} [ll]|{(\Ll_{n-1},\lambda_{n-1})}
& 
\dots
&
\F_2
\ar[rr]^1
& &
\F_1
\ar[rr]^0
\ar @/^0.65pc/ @{} [ll]|{(\Ll_1,\lambda_1)}
\ar@/_1.2pc/|{\tensor}[ll]
& &
\F_{S}
\ar @/^0.65pc/ @{} [ll]|{(\Ll_0,\lambda_0)}
\ar@/_1.2pc/|{\tensor}[ll]
}
\right)
$$
}
is a generalized isomorphism from $\E_S$ to $\F_S$ with
$\mu_i=\lambda_j=0$ for $i\in I$ and $j\in J$,
\begin{eqnarray*}
& &
F_\bullet\E=(0=F_0\E\subseteq\dots\subseteq F_{r+s+1}\E=\E_S)
\\ \text{and} & &
F_\bullet\F=(0=F_0\F\subseteq\dots\subseteq F_{r+s+1}\F=\F_S)
\end{eqnarray*}
are flags of type ${\bf d}$ and $\deltabf$ in $\E_S$ and
$\F_S$ respectively,
$$
\varphi_p =
\left(
\xymatrix@C-1pc{
\E_0^{(p)}
\ar@/^1.2pc/|{\tensor}[rr]
& & 
\E_1^{(p)}
\ar[ll]_{0}
\ar @/^1.0pc/ @{} [ll]|{(\M_0^{(p)},\mu_0^{(p)})}
&
\dots
&
\E_{m_p-1}^{(p)}
\ar@/^1.2pc/|{\tensor}[rr]
& &
\E_{m_p}^{(p)}
\ar[ll]_{m_p}
\ar @/^1.0pc/ @{} [ll]|{(\M_{m_p}^{(p)},\mu_{m_p}^{(p)})}
}
\right)
$$
is a complete collineation from 
$\E_{m_p}^{(p)}=F_{r-p+1}\F/F_{r-p}\F$ to
$\E_0^{(p)}=F_{s+p+1}\E/F_{s+p}\E$
for $1\leq p\leq r$,
$$
\psi_q =
\left(
\xymatrix@C-1pc{
\F_{n_q}^{(q)}
\ar[rr]^{n_q}
\ar @/_1.0pc/ @{} [rr]|{(\Ll_{n_q}^{(q)},\lambda_{n_q}^{(q)})}
& &
\F_{n_q-1}^{(q)}
\ar@/_1.2pc/|{\tensor}[ll]
&
\dots
&
\F_1^{(q)}
\ar[rr]^{0}
\ar @/_1.0pc/ @{} [rr]|{(\Ll_0^{(q)},\lambda_0^{(q)})}
& & 
\F_0^{(p)}
\ar@/_1.2pc/|{\tensor}[ll]
}
\right)
$$
is a complete collineation from 
$\F_{n_q}^{(q)}=F_{s-q+1}\F/F_{s-q}\E$ to
$\F_0^{(q)}=F_{r+q+1}\F/F_{r+q}\F$
for $s\geq q\geq 1$ and
$\Phi' =$ 
{\footnotesize
$$
\left(
\xymatrix@C-4ex{
\E'_0
\ar@/^1.2pc/|{\tensor}[rr]
\ar @/_0.8pc/ @{} [rr]|{(\M'_0,\mu'_0)}
& &
\E'_1 
\ar[ll]_0
\ar @/_0.8pc/ @{} [rr]|{(\M'_1,\mu'_1)}
\ar@/^1.2pc/|{\tensor}[rr]
& &
\E'_2
\ar[ll]_1
& 
\dots
& 
\E'_{n'-1}
\ar @/_0.8pc/ @{} [rr]|{(\M'_{n'-1},\mu'_{n'-1})}
\ar@/^1.2pc/|{\tensor}[rr]
& & 
\E'_{n'}
\ar[ll]_{n'-1}
\ar[rr]^\sim
& & 
\F'_{n'}
\ar[rr]^{n'-1}
& & 
\F'_{n'-1}
\ar@/_1.2pc/|{\tensor}[ll]
\ar @/^0.8pc/ @{} [ll]|{(\Ll'_{n'-1},\lambda'_{n'-1})}
& 
\dots
&
\F'_2
\ar[rr]^1
& &
\F'_1
\ar[rr]^0
\ar @/^0.8pc/ @{} [ll]|{(\Ll'_1,\lambda'_1)}
\ar@/_1.2pc/|{\tensor}[ll]
& &
\F'_0
\ar @/^0.8pc/ @{} [ll]|{(\Ll'_0,\lambda'_0)}
\ar@/_1.2pc/|{\tensor}[ll]
}
\right)
$$
}
is a generalized isomorphism from
$\E'_0=F_{s+1}\E/F_s\E$
to
$\F'_0=F_{r+1}\F/F_r\F$.

The mapping
$$
\Phi\mapsto
((F_\bullet\E,F_\bullet\F),
\varphi_1,\dots,\varphi_r,\psi_s,\dots,\psi_1,\Phi')
$$
is defined as follows:
Let $\Phi$ as above be given. 
For convenience we set $\E_{n+1}:=\F_n$,
$\F_{n+1}:=\E_n$ and we let
$\F_{n+1}\to\F_n$ 
and
$\E_n\ot\E_{n+1}$ 
be the isomorphism $\E_n\isomto\F_n$ and its inverse respectively,
whereas we let
$
\xymatrix@C-1ex{
\E_{n}
\ar|(.45){\tensor}[r]
&
\E_{n+1}
}
$ 
and 
$
\xymatrix@C-1ex{
\F_{n+1}
&
\F_{n}
\ar|(.45){\tensor}[l]
}
$
both be the zero morphism.
For what follows, the picture below may help to keep 
track of the indices:
$$
\overunderbraces%
{&\br{13}{\Phi}&}
{
\overset{0}{|}
&
\dots\dots
\overset{n-j_1}{|}
&
\dots
&
\overset{i_1}{|}
&
\dots
&
\overset{i_2}{|}
\dots\dots
\overset{i_r}{|}
&
\dots
&
\overset{n}{|}
\overset{n}{|}
&
\dots
&
\overset{j_s}{|}
\dots\dots
\overset{j_2}{|}
&
\dots
&
\overset{j_1}{|}
&
\dots
&
\overset{n-i_1}{|}
\dots\dots
&
\overset{0}{|}
}{ 
& &
\br{1}{\Phi'}
& &
\br{1}{\varphi_1}
& &
\br{1}{\varphi_r}
& &
\br{1}{\psi_s}
& &
\br{1}{\psi_1}
& &
\br{1}{\Phi'}
& &
}
$$

Let
$F_0\E=F_0\F:=0$, $F_{r+s+1}\E:=\E$, $F_{r+s+1}\F:=\F$ and
\begin{eqnarray*}
F_p\E &:=&
\left\{
\begin{array}{ll}
\left(
\begin{array}{ll}
\text{image of $\ker(\E_n\to\F_{j_{s-p+1}})$ by}\\ 
\text{the morphism $\E_S\ot\E_n$}
\end{array}
\right)\ ,\ 
\text{if $1\leq p\leq s$}
\\
\quad
\text{$\ker(
\xymatrix@C-1ex{
\E_S
\ar|(.35){\tensor}[r]
&
\E_{i_{p-s}+1}
}
)$\ ,\quad 
if $s+1\leq p\leq r+s$}
\end{array}
\right.
\\
F_q\F &:=&
\left\{
\begin{array}{ll}
\left(
\begin{array}{ll}
\text{image of $\ker(\E_{i_{r-q+1}}\ot\F_n)$ by}\\ 
\text{the morphism $\F_n\to\F_S$}
\end{array}
\right)\ ,\ 
\text{if $1\leq q\leq r$}
\\
\quad
\text{$\ker(
\xymatrix@C-1ex{
\F_{j_{q-r}+1}
&
\F_S
\ar|(.35){\tensor}[l]
}
)$\ ,\quad  
if $r+1\leq q\leq r+s$}
\end{array}
\right.
\end{eqnarray*}
It is then clear from the definition of generalized isomorphisms
that
\begin{eqnarray*}
& &F_\bullet\E:=(0=F_0\E\subseteq\dots\subseteq F_{r+s+1}\E=\E_S)
\\
\text{and} & &
F_\bullet\F:=(0=F_0\F\subseteq\dots\subseteq F_{r+s+1}\F=\F_S)
\end{eqnarray*}
are flags of type ${\bf d}$ and $\deltabf$ in $\E_S$ and
$\F_S$ respectively.
Let $1\leq p\leq r$. We set
\begin{eqnarray*}
\E_0^{(p)} &:=&
\ker(
\xymatrix@C-1ex{
\E
\ar|(.3){\tensor}[r]
&
\E_{i_{p+1}+1}
})/
\ker(
\xymatrix@C-1ex{
\E
\ar|(.3){\tensor}[r]
&
\E_{i_{p}+1}
}) =
F_{s+p+1}\E/F_{s+p}\E
\\
\M_0^{(p)} &:=& \Tensor_{i=0}^{i_p}\M_i\ ,
\qquad
\mu_0^{(p)}:=0
\end{eqnarray*}
and
\begin{eqnarray*}
\E_k^{(p)} &:=&
\ker(\E_{i_p}\ot\E_{i_p+k})
\cap
\ker(
\xymatrix@C-1ex{
\E_{i_p+k}
\ar|(.45){\tensor}[r]
&
\E_{i_{p+1}+1}
})
\quad
(1\leq k\leq m_p)
\\
\M_k^{(p)} &:=& \M_{i_p+k}\ ,
\qquad
\mu_k^{(p)}:=\mu_{i_p+k}
\qquad\qquad\qquad\quad
(1\leq k\leq m_p-1)
\end{eqnarray*}
where $m_p=i_{p+1}-i_p$. 
Observe that the sheaves $\E_k^{(p)}$ thus defined are locally
free of rank $m_p$. Indeed, this is clear for $k=0$.
For $k\geq 1$ it suffices to show that
$\E_{i_p+k}$ is generated by the two subsheaves
$\ker(\E_{i_p}\ot\E_{i_p+k})$ and 
$
\ker(
\xymatrix@C-1ex{
\E_{i_p+k}
\ar|(.45){\tensor}[r]
&
\E_{i_{p+1}+1}
})
$.
For this in turn, it suffices to show that the image of 
$
\ker(
\xymatrix@C-1ex{
\E_{i_p+k}
\ar|(.45){\tensor}[r]
&
\E_{i_{p+1}+1}
})
$
by the morphism $\E_{i_p}\ot\E_{i_p+k}$ is 
$\im(\E_{i_p}\ot\E_{i_p+k})$.
But this is clear, since by the definition of generalized
isomorphisms we have
\begin{eqnarray*}
& &
\ker(
\xymatrix@C-1ex{
\E_{i_p+k}
\ar|(.45){\tensor}[r]
&
\E_{i_{p+1}+1}
})
\supseteq
\im(\E_{i_p+k}\ot\E_{i_{p+1}+1})
\\
\text{and}\quad & &
\im(\E_{i_p}\ot\E_{i_p+k})
=
\im(\E_{i_p}\ot\E_{i_{p+1}+1})
\quad.
\end{eqnarray*}
Since 
$$
\Tensor_{i=0}^{i_p}\M_i^{\dual}\tensor(\E_S/F_{s+p}\E)
=
\im(\Tensor_{i=0}^{i_p}\M_i^{\dual}\tensor\E_S\to\E_{i_p+1})
=
\ker(\E_{i_p}\ot\E_{i_p+1})
\quad,
$$
we have a natural isomorphism
$
\E_0^{(p)}=F_{s+p+1}\E/F_{s+p}\E
\isomto
\M_0^{(p)}\tensor\E_1^{(p)}
$.
We define
$
\E_0^{(p)} \ot \E_1^{(p)}
$
to be the zero morphism.
Thus we have a bf-morphism
$
\xymatrix@C-1pc{
\E_0^{(p)}
\ar@/^1.2pc/|{\tensor}[rr]
& & 
\E_{1}^{(p)}
\ar[ll]_{0}
\ar @/^1.0pc/ @{} [ll]|{(\M_0^{(p)},\mu_0^{(p)}=0)}
}
$
of rank zero.
For $1\leq k\leq m_p-1$ let
$
\xymatrix@C-1pc{
\E_k^{(p)}
\ar@/^1.2pc/|{\tensor}[rr]
& & 
\E_{k+1}^{(p)}
\ar[ll]_{k}
\ar @/^1.0pc/ @{} [ll]|{(\M_k^{(p)},\mu_k^{(p)})}
}
$
be the bf-morphism induced by the bf-morphism
$
\xymatrix@C-1pc{
\E_{i_p+k}
\ar@/^1.2pc/|{\tensor}[rr]
& & 
\E_{i_p+k+1}
\ar[ll]_{i_p+k}
\ar @/^1.0pc/ @{} [ll]|{(\M_{i_p+k},\mu_{i_p+k})}
}
$.
Observe that
$
\ker(
\xymatrix@C-1ex{
\E_{i_{p+1}}
\ar|(.45){\tensor}[r]
&
\E_{i_{p+1}+1}
})
=
\im(\E_{i_{p+1}}\ot\E_n)
=
\E_n/\ker(\E_{i_{p+1}}\ot\E_n)
$
and that the morphism $\E_n\to\F_S$ maps
$\ker(\E_{i_p}\ot\E_n)$ injectively into $\F_S$.
Therefore we have a natural isomorphism
$
\E_{m_p}^{(p)}\isomorph F_{r-p+1}\F/F_{r-p}\F
$
by which we identify these two sheaves.
It is not difficult to see that
$$
\varphi_p :=
\left(
\xymatrix@C-1pc{
\E_0^{(p)}
\ar@/^1.2pc/|{\tensor}[rr]
& & 
\E_1^{(p)}
\ar[ll]_{0}
\ar @/^1.0pc/ @{} [ll]|{(\M_0^{(p)},\mu_0^{(p)})}
&
\dots
&
\E_{m_p-1}^{(p)}
\ar@/^1.2pc/|{\tensor}[rr]
& &
\E_{m_p}^{(p)}
\ar[ll]_{m_p}
\ar @/^1.0pc/ @{} [ll]|{(\M_{m_p}^{(p)},\mu_{m_p}^{(p)})}
}
\right)
$$
is a complete collineation in the 
sense of \ref{complete collineations} from 
$\E_{m_p}^{(p)}=F_{r-p+1}\F/F_{r-p}\F$ to \linebreak[4]
$\E_0^{(p)}=F_{s+p+1}\E/F_{s+p}\E$.
In a completely symmetric way the generalized isomorphism
$\Phi$ induces also complete collineations
$$
\psi_q =
\left(
\xymatrix@C-1pc{
\F_{n_q}^{(q)}
\ar[rr]^{n_q}
\ar @/_1.0pc/ @{} [rr]|{(\Ll_{n_q}^{(q)},\lambda_{n_q}^{(q)})}
& &
\F_{n_q-1}^{(q)}
\ar@/_1.2pc/|{\tensor}[ll]
&
\dots
&
\F_1^{(q)}
\ar[rr]^{0}
\ar @/_1.0pc/ @{} [rr]|{(\Ll_0^{(q)},\lambda_0^{(q)})}
& & 
\F_0^{(p)}
\ar@/_1.2pc/|{\tensor}[ll]
}
\right)
$$
from
$\F_{n_q}^{(q)}=F_{s-q+1}\F/F_{s-q}\E$ to
$\F_0^{(q)}=F_{r+q+1}\F/F_{r+q}\F$
for $s\geq q\geq 1$.
It remains to construct the generalized isomorphism 
$\Phi'$. We set
\begin{eqnarray*}
\E'_k &:=&
\ker(
\xymatrix@C-1ex{
\E_{n-j_1+k}
\ar|(.55){\tensor}[r]
&
\E_{i_1+1}
})
/
\im(\E_{n-j_1+k}\ot\ker(\E_n\to\F_{j_1}))
\\
\F'_k &:=&
\ker(
\xymatrix@C-1ex{
\F_{j_1+1}
&
\F_{n-i_1+k}
\ar|(.55){\tensor}[l]
})
/
\im(\ker(\E_{i_1}\ot\F_n)\to\F_{n-i_1+k}))
\end{eqnarray*}
for $0\leq k\leq n':=i_1+j_1-n$ and
\begin{eqnarray*}
\M'_k &:=& \M_{n-j_1+k}\ ,
\quad
\mu'_k:=\mu_{n-j_1+k}
\\
\Ll'_k &:=& \Ll_{n-i_1+k}\ \ \ ,
\quad
\lambda'_k:=\lambda_{n-i_1+k}
\end{eqnarray*}
for $0\leq k\leq n'-1$.
It is then clear that the $\E'_k$ and $\F'_k$ are
locally free of rank $n'=i_1+j_1-n$.
It follows from definition \ref{generalized isomorphism}.2.
that the $\mu_i$ and $\lambda_j$ are nowhere vanishing
for $0\leq i\leq n-j_1-1$ and $n-i_1-1\geq j\geq 0$.
Therefore we may identify $\E_i$ with $\E_S$ and
$\F_j$ with $\F_S$ 
for $0\leq i\leq n-j_1-1$ and $n-i_1-1\geq j\geq 0$
respectively. This implies in particular that we have
$\E'_0=F_{s+1}\E/F_s\E$ and $\F'_0=F_{r+1}\F/F_r\F$.
Let
$$
\xymatrix@C-1.5pc{
\E'_k
\ar@/^1.2pc/|{\tensor}[rr]
& & 
\E'_{k+1}
\ar[ll]_{k}
\ar @/^0.9pc/ @{} [ll]|{(\M'_k,\mu'_k)}
}
\quad\text{and}\quad
\xymatrix@C-1.5pc{
\F'_{k+1}
\ar[rr]^{k}
\ar @/_0.9pc/ @{} [rr]|{(\Ll'_{k},\lambda'_{k})}
& &
\F'_{k}
\ar@/_1.2pc/|{\tensor}[ll]
}
$$
be the bf-morphisms induced by the bf-morphisms
$$
\xymatrix@C-1pc{
\E_{n-j_1+k}
\ar@/^1.2pc/|{\tensor}[rr]
& & 
\E_{n-j_1+k+1}
\ar[ll]_{n-j_1+k}
\ar @/^0.9pc/ @{} [ll]|{(\M_{n-j_1+k},\mu_{n-j_1+k})}
}
\quad\text{and}\quad
\xymatrix@C-1pc{
\F_{n-i_1+k+1}
\ar[rr]^{n-i_1+k}
\ar @/_0.9pc/ @{} [rr]|{(\Ll_{n-i_1+k},\lambda_{n-i_1+k})}
& &
\F_{n-i_1+k}
\ar@/_1.2pc/|{\tensor}[ll]
}
$$
respectively.
We have 
$$
\ker(
\xymatrix@C-1ex{
\E_{i_1}
\ar|(.45){\tensor}[r]
&
\E_{i_1+1}
})
=
\im(\E_{i_1}\ot\E_n)
=
\E_n/\ker(\E_{i_1}\ot\E_n)
$$
and therefore 
$$
\E'_{n'}=\E_n/(\ker(\E_{i_1}\ot\E_n) + \ker(\E_n\to\F_{j_1}))
\quad.
$$
By the same argument:
$$
\F'_{n'}=\F_n/(\ker(\E_{i_1}\ot\F_n) + \ker(\F_n\to\F_{j_1}))
\quad.
$$
Thus the isomorphism $\E_n\isomto\F_n$ induces an isomorphism
$\E'_{n'}\isomto\F'_{n'}$. Again it is not difficult to check
that $\Phi' :=$ 
{\footnotesize
$$
\left(
\xymatrix@C-4ex{
\E'_0
\ar@/^1.2pc/|{\tensor}[rr]
\ar @/_0.8pc/ @{} [rr]|{(\M'_0,\mu'_0)}
& &
\E'_1 
\ar[ll]_0
\ar @/_0.8pc/ @{} [rr]|{(\M'_1,\mu'_1)}
\ar@/^1.2pc/|{\tensor}[rr]
& &
\E'_2
\ar[ll]_1
& 
\dots
& 
\E'_{n'-1}
\ar @/_0.8pc/ @{} [rr]|{(\M'_{n'-1},\mu'_{n'-1})}
\ar@/^1.2pc/|{\tensor}[rr]
& & 
\E'_{n'}
\ar[ll]_{n'-1}
\ar[rr]^\sim
& & 
\F'_{n'}
\ar[rr]^{n'-1}
& & 
\F'_{n'-1}
\ar@/_1.2pc/|{\tensor}[ll]
\ar @/^0.8pc/ @{} [ll]|{(\Ll'_{n'-1},\lambda'_{n'-1})}
& 
\dots
&
\F'_2
\ar[rr]^1
& &
\F'_1
\ar[rr]^0
\ar @/^0.8pc/ @{} [ll]|{(\Ll'_1,\lambda'_1)}
\ar@/_1.2pc/|{\tensor}[ll]
& &
\F'_0
\ar @/^0.8pc/ @{} [ll]|{(\Ll'_0,\lambda'_0)}
\ar@/_1.2pc/|{\tensor}[ll]
}
\right)
$$
}
is a generalized isomorphism from
$\E'_0=F_{s+1}\E/F_s\E$
to
$\F'_0=F_{r+1}\F/F_r\F$.
This completes the construction of the mapping
$$
\Phi\mapsto
((F_\bullet\E,F_\bullet\F),
\varphi_1,\dots,\varphi_r,\psi_s,\dots,\psi_1,\Phi')
\quad.
$$

We proceed by constructing the inverse
of this mapping.
Let data
$
((F_\bullet\E,F_\bullet\F),
\varphi_1,\dots,\varphi_r,\psi_s,\dots,\psi_1,\Phi')
$
be given.
Let $\E_i:=\E_S$ and $\F_j:=\F_S$ for $0\leq i\leq n-j_1$ and
$n-i_1\geq j\geq 0$ respectively.
Now let $n-j_1+1\leq i\leq i_1$ and $j_1\geq j\geq n-i_1+1$.
Let $\tilde{\E}_i$ and $\tilde{\F}_j$ be defined by
the cartesian diagrams
$$
\vcenter{
\xymatrix{
\text{$\tilde{\E}_i$}
\ar[d]
\ar[r]
&
\E'_{i+j_1-n}
\ar[d]
\\
F_{s+1}\E
\ar@{->>}[r]
&
F_{s+1}\E/F_s\E
}}
\quad \text{and} \quad
\vcenter{
\xymatrix{
\text{$\tilde{\F}_j$}
\ar[d]
\ar[r]
&
\F'_{j+i_1-n}
\ar[d]
\\
F_{r+1}\F
\ar@{->>}[r]
&
F_{r+1}\F/F_r\F
}}
$$
respectively.
For a moment let 
$\M:=\Tensor_{k=0}^{i+j_1-n-1}\M'_k$.
We have a commutative diagram 
$$
\vcenter{
\xymatrix{
\M^{\dual}\tensor F_{s+1}\E
\ar[d]
\ar[r]
&
\E'_{i+j_1-n}
\ar[d]        
\\
F_{s+1}\E
\ar@{->>}[r]
&
F_{s+1}\E/F_s\E
}}
\eqno(*)
$$
where the left vertical arrow is induced by
$
\Tensor_{k=0}^{i+j_1-n-1}\mu'_k:\Oo_S\to\M
$
and the upper horizontal arrow is the composition
$$
\M^{\dual}\tensor F_{s+1}\E
\to
\M^{\dual}\tensor F_{s+1}\E/F_s\E
=
\M^{\dual}\tensor\E'_0
\to
\E'_{i+j_1-n}
$$
The diagram $(*)$ induces a morphism
$\M^{\dual}\tensor F_{s+1}\E\to\tilde{\E}_i$.
Analogously, we have a morphism
$\Ll^{\dual}\tensor F_{r+1}\F\to\tilde{\F}_j$,
where we have employed the abbreviation 
$\Ll:=\Tensor_{k=0}^{j+i_1-n-1}\Ll'_k$.
Let $\E_i$ and $\F_j$ be defined by the cocartesian diagrams
$$
\vcenter{
\xymatrix{
\M^{\dual}\tensor F_{s+1}\E
\ar[d]
\ar@<-.5ex>@{^{(}->}[r]
&
\M^{\dual}\tensor \E_S
\ar[d]
\\
\text{$\tilde{\E}_i$}
\ar[r]
&
\E_i
}}
\quad\text{and}\quad
\vcenter{
\xymatrix{
\Ll^{\dual}\tensor F_{r+1}\F
\ar[d]
\ar@<-.5ex>@{^{(}->}[r]
&
\Ll^{\dual}\tensor \F_S
\ar[d]
\\
\text{$\tilde{\F}_j$}
\ar[r]
&
\F_j
}}
$$
respectively.

We define $\E_n=\F_n$ by the cartesian diagram
$$
\xymatrix@C-2ex{
\E_n=\F_n
\ar[d]
\ar[rrr]
& & &
\text{$\tilde{\F}_{j_1}$}
\ar@{->>}[d]
\\
\text{$\tilde{\E}_{i_1}$}
\ar@{->>}[rr]
& &
\E'_{i_1+j_1-n}
\ar^{\isomorph}[r]
&
\F'_{i_1+j_1-n}
}
$$
Observe that the composed morphism
$\E_n\to\tilde{\E}_{i_1}\to F_{s+1}\E$
maps the submodule $\ker(\E_n\to\tilde{\F}_{j_1})$
of $\E_n$ isomorphically onto the submodule
$F_s\E$ of $F_{s+1}\E$.
Therefore we have canonical injections
$$
F_p\E\injto F_s\E\isomto\ker(\E_n\to\tilde{\F}_{j_1})\injto\E_n
$$
for $0\leq p\leq s$.
Analogously, we have canonical injections
$$
F_q\F\injto F_r\F\isomto\ker(\F_n\to\tilde{\E}_{i_1})\injto\F_n
$$
for $0\leq q\leq r$.

Now let $1\leq p\leq r,\ i_p+1\leq i\leq i_{p+1}$
and $s\geq q\geq 1,\ j_{q+1}\geq j\geq j_q+1$.
We want to define $\E_i$ and $\F_j$ in this case.
Let first $\tilde{\E}_i$ and $\tilde{\F}_j$ be defined
by the cocartesian diagrams
$$
\vcenter{
\xymatrix@C-3ex{
\M^{\dual}\tensor F_{s+p+1}\E/F_{s+p}\E
\ar[r]
\ar@<-.5ex>@{^{(}->}[d]
&
\E_{i-i_p}^{(p)}
\ar[d]
\\
\M^{\dual}\tensor \E_S/F_{s+p}\E
\ar[r]
&
\text{$\tilde{\E}_i$}
}}
\quad \text{and}\quad
\vcenter{
\xymatrix@C-3ex{
\Ll^{\dual}\tensor F_{r+q+1}\F/F_{r+q}\E
\ar[r]
\ar@<-.5ex>@{^{(}->}[d]
&
\F_{j-j_q}^{(q)}
\ar[d]
\\
\Ll^{\dual}\tensor \F_S/F_{r+q}\F
\ar[r]
&
\text{$\tilde{\F}_j$}
}}
$$
where we have set
$\M:=\Tensor_{k=0}^{i-i_p-1}\M_k^{(p)}$
and
$\Ll:=\Tensor_{k=0}^{j-j_q-1}\Ll_k^{(q)}$.
Let furthermore $\hat{\E}_i$ and $\hat{\F}_j$ be defined
by the cocartesian diagrams
$$
\vcenter{
\xymatrix{
F_{r-p+1}\F/F_{r-p}\F
\ar[r]
\ar@<-.5ex>@{^{(}->}[d]
&
\E_{i-i_p}^{(p)}
\ar[d]
\\
\F_n/F_{r-p}\F
\ar[r]
&
\text{$\hat{\E}_i$}
}}
\quad\text{and}\quad
\vcenter{
\xymatrix{
F_{s-q+1}\E/F_{s-q}\E
\ar[r]
\ar@<-.5ex>@{^{(}->}[d]
&
\F_{j-j_q}^{(q)}
\ar[d]
\\
\E_n/F_{s-q}\E
\ar[r]
&
\text{$\hat{\F}_j$}
}}
$$
Now we define $\E_i$ and $\F_j$ by the cocartesian diagrams
$$
\vcenter{
\xymatrix@R-3ex{
\E_{i-i_p}^{(p)}
\ar[d]
\ar[r]
&
\text{$\tilde{\E}_i$}
\ar[d]
\\
\text{$\hat{\E}_i$}
\ar[r]
&
\E_i
}}
\quad\text{and}\quad
\vcenter{
\xymatrix@R-3ex{
\F_{j-j_q}^{(q)}
\ar[d]
\ar[r]
&
\text{$\tilde{\F}_j$}
\ar[d]
\\
\text{$\hat{\F}_j$}
\ar[r]
&
\F_j
}}
$$
respectively.
For $p=r$ and $i=i_{r+1}=n$ this gives formally a new definition
of $\E_n$, but it is clear that we have a canonical isomorphism
between the two $\E_n$'s. A similar remark applies to $\F_n$.

We define the invertible sheaves $\M_i$ together
with their respective sections $\mu_i$  as follows:
\begin{eqnarray*}
\M_i &:=& \Oo_S
\ ,\quad
\mu_i:=1
\ \ \quad
(0\leq i\leq n-j_1-1)
\\
\M_i &:=& \M'_{i+j_1-n}
\ ,\quad
\mu_i:=\mu'_{i+j_1-n}
\ \ \quad
(n-j_1\leq i\leq i_1-1)
\\
\M_i &:=& \M^{(p)}_{i-i_p}
\ ,\quad
\mu_i:=\mu^{(p)}_{i-i_p}
\ \ \quad
(1\leq p\leq r\ ,\ \ i_p<i<i_{p+1})
\end{eqnarray*}
\begin{eqnarray*}
\M_{i_1} &:=& \M_0^{(1)}\tensor\Tensor_{k=0}^{i_1+j_1-n-1}(\M'_k)^{\dual}
\ ,\quad
\mu_{i_1}:=0
\ \ \quad
\\
\M_{i_p} &:=& \M_0^{(p)}\tensor\Tensor_{k=0}^{i_p-i_{p-1}-1}(\M_k^{(p-1)})^{\dual}
\ ,\quad
\mu_{i_p}:=0
\ \ \quad
(2\leq p\leq r)
\end{eqnarray*}
Let the $\Ll_j$ and $\lambda_j$ be defined symmetrically
(i.e. by replacing in the above definition the letter $\M$ with $\Ll$,
$\mu$ with $\lambda$, $i$ with $j$, $j$ with $i$ and $r$ with $s$).

It remains to define the bf-morphisms
$$
\xymatrix@C-1.5pc{
\E_i
\ar@/^1.2pc/|{\tensor}[rr]
& & 
\E_{i+1}
\ar[ll]_{i}
\ar @/^.85pc/ @{} [ll]|{(\M_i,\mu_i)}
}
\quad\text{and}\quad
\xymatrix@C-1.5pc{
\F_{j+1}
\ar[rr]^{j}
\ar @/_.85pc/ @{} [rr]|{(\Ll_{j},\lambda_{j})}
& &
\F_{j}
\ar@/_1.2pc/|{\tensor}[ll]
}
$$
for $n-j_1\leq i\leq n-1$ and $n-1\geq j\geq n-i_1$.
Again we restrict ourselves to the left hand side, since
the right hand side is obtained by the symmetric construction.
For $n-j_1\leq i\leq i_1-1$ 
(respectively for $1\leq p\leq r\ ,\ i_p\leq i\leq i_{p+1}-1$) 
the bf-morphism
$
\xymatrix@C-1.5pc{
\E_i
\ar@/^1.2pc/|{\tensor}[rr]
& & 
\E_{i+1}
\ar[ll]
}
$
is induced in an obvious way by the bf-morphism
$
\xymatrix@C-1.5pc{
\E'_{i+j_1-n}
\ar@/^1.2pc/|{\tensor}[rr]
& & 
\E'_{i+j_1-n+1}
\ar[ll]
}
$
(respectively by the bf-morphism
$
\xymatrix@C-1.5pc{
\E^{(p)}_{i-i_p}
\ar@/^1.2pc/|{\tensor}[rr]
& & 
\E^{(p)}_{i-i_p+1}
\ar[ll]
}
$
).
For the definition of the bf-morphism
$
\xymatrix@C-1.5pc{
\E_{i_1}
\ar@/^1.2pc/|{\tensor}[rr]
& & 
\E_{i_1+1}
\ar[ll]
}
$
consider the two canonical exact sequences
$$
\xymatrix@C=3ex@R-4ex{
0 \ar[r] &
\text{$\tilde{\E}$}_{i_1} \ar[r] &
\E_{i_1} \ar[r] &
\M^{\dual}\tensor\E_S/F_{s+1}\E \ar[r] &
0
\\
0 \ar[r] &
\text{$\tilde{\E}$}_{i_1+1} \ar[r] &
\E_{i_1+1}\ar[r] &
\text{$\hat{\E}$}_{i_1+1}/\E_1^{(1)}\ar[r] & 
0
}
$$
where $\M:=\Tensor_{k=0}^{i_1-1}\M_k$.
Observe that we have canonical isomorphisms
$$
\xymatrix@C=3ex@R-4ex{
\text{$\tilde{\E}$}_{i_1}\ar[r]^(.4){a}_(.4){\isomorph} &
\F_n/F_{r}\F\ar[r]^{b}_{\isomorph} &
\text{$\hat{\E}$}_{i_1+1}/\E_1^{(1)}
\\
\text{$\tilde{\E}$}_{i_1+1}\ar[rr]^(.4){c}_(.4){\isomorph} & &
\M_{i_1}^{\dual}\tensor\M^{\dual}\tensor \E_S/F_{s+1}\E
}
$$
The isomorphism $a$ follows from the observation we made
after the definition of $\E_n=\F_n$, namely that the composed
morphism $\F_n\to\tilde{\F}_{j_1}\to F_{r+1}\F$ maps
$\ker(\E_n\to\tilde{\E}_{i_1})$ isomorphically to $F_r\F$.
The isomorphism $b$ comes from the fact that for $i=i_1+1$ 
the left vertical arrow in the defining diagram for $\hat{\E}_i$
vanishes, and the isomorphism $c$ follows since for $i=i_1+1$
the left vertical arrow in the defining diagram for $\tilde{\E}_i$
is an isomorphism. Thus we have morphisms
$$
\xymatrix@R-4ex{
\E_{i_1}\ar@{->>}[r] &
\M^{\dual}\tensor\E_S/F_{s+1} \ar[r]^{c^{-1}} &
\M_{i_1}\tensor\text{$\tilde{\E}$}_{i_1+1}\ar@{^{(}->}[r] &
\M_{i_1}\tensor\E_{i_1+1}
\\
\E_{i_1} &
\text{$\tilde{\E}$}_{i_1}\ar@{_{(}->}[l] &
\text{$\hat{\E}$}_{i_1+1}/\E_1^{(1)} \ar[l]_{a^{-1}b^{-1}} &
\E_{i_1+1}\ar@{->>}[l] 
}
$$
which make up the bf-morphism 
$
\xymatrix@C-1.5pc{
\E_{i_1}
\ar@/^1.2pc/|{\tensor}[rr]
& & 
\E_{i_1+1}
\ar[ll]
}
$.
For $2\leq p\leq r$ the bf-morphism
$
\xymatrix@C-1.5pc{
\E_{i_p}
\ar@/^1.2pc/|{\tensor}[rr]
& & 
\E_{i_p+1}
\ar[ll]
}
$
is constructed similarly from the exact sequences
$$
\xymatrix@R-4ex{
0 \ar[r] &
\text{$\hat{\E}$}_{i_p} \ar[r] &
\E_{i_p} \ar[r] &
\text{$\tilde{\E}$}_{i_p}/\E_{i_p-i_{p-1}}^{(p-1)} \ar[r] &
0
\\
0 \ar[r] &
\text{$\tilde{\E}$}_{i_p+1} \ar[r] &
\E_{i_p+1}\ar[r] &
\text{$\hat{\E}$}_{i_p+1}/\E_1^{(p)}\ar[r] & 
0
}
$$
and the canonical isomorphisms
$$
\xymatrix@C=3ex@R-4ex{
\text{$\hat{\E}$}_{i_p}\ar[r]^(.4){\isomorph} &
\F_n/F_{r-p+1}\F\ar[r]^{\isomorph} &
\text{$\hat{\E}$}_{i_p+1}/\E_1^{(p)}
\\
\text{$\tilde{\E}$}_{i_p+1}\ar[r]^(.25){\isomorph} & 
\M_{i_p}^{\dual}\tensor\M^{\dual}\tensor \E/F_{s+p}\E
\ar[r]^(.55){\isomorph} & 
\M_{i_p}^{\dual}\tensor\text{$\tilde{\E}$}_{i_p}/\E_{i_p-i_{p-1}}^{(p-1)}
}
$$
where $\M:=\Tensor_{k=0}^{i_p-1}\M_k$.

This completes the construction of
{\footnotesize
$$
\Phi = 
\left(
\xymatrix@C-4ex{
\E_{S}
\ar@/^1.2pc/|{\tensor}[rr]
\ar @/_0.65pc/ @{} [rr]|{(\M_0,\mu_0)}
& &
\E_1 
\ar[ll]_0
\ar @/_0.65pc/ @{} [rr]|{(\M_1,\mu_1)}
\ar@/^1.2pc/|{\tensor}[rr]
& &
\E_2
\ar[ll]_1
& 
\dots
& 
\E_{n-1}
\ar @/_0.65pc/ @{} [rr]|{(\M_{n-1},\mu_{n-1})}
\ar@/^1.2pc/|{\tensor}[rr]
& & 
\E_n
\ar[ll]_{n-1}
\ar[rr]^\sim
& & 
\F_n
\ar[rr]^{n-1}
& & 
\F_{n-1}
\ar@/_1.2pc/|{\tensor}[ll]
\ar @/^0.65pc/ @{} [ll]|{(\Ll_{n-1},\lambda_{n-1})}
& 
\dots
&
\F_2
\ar[rr]^1
& &
\F_1
\ar[rr]^0
\ar @/^0.65pc/ @{} [ll]|{(\Ll_1,\lambda_1)}
\ar@/_1.2pc/|{\tensor}[ll]
& &
\F_{S}
\ar @/^0.65pc/ @{} [ll]|{(\Ll_0,\lambda_0)}
\ar@/_1.2pc/|{\tensor}[ll]
}
\right)
$$
}
It is not difficult to see that $\Phi$ 
is a generalized isomorphism from $\E_S$ to $\F_S$ and that the mapping
$
((F_\bullet\E,F_\bullet\F),
\varphi_1,\dots,\varphi_r,\psi_s,\dots,\psi_1,\Phi')
\mapsto
\Phi
$
is inverse to the mapping
$
\Phi
\mapsto
((F_\bullet\E,F_\bullet\F),
\varphi_1,\dots,\varphi_r,\psi_s,\dots,\psi_1,\Phi')
$
constructed before.
We leave the details to the reader.
\end{proof}

In the situation of theorem \ref{strata} we denote by $\Gl(\E)$ 
the group scheme over $T$, whose $S$-valued points are
the automorphisms of $\E_S$. There is a natural left operation of 
$\Gl(\E)\times_{T}\Gl(\F)$ on $\KGl(\E,\F)$, which is given
on $S$-valued points by
$$
(f,g)
\left(
\hbox{\small
\xymatrix@C-1.5pc{
\E_{S}
\ar@/^0.8pc/|{\tensor}[rr]^{u^{\flat}}
& &
\E_1 
\ar[ll]^{u^{\sharp}}
& &
\dots
& &
\F_1
\ar[rr]_{v^{\sharp}}
& &
\F_{S}
\ar@/_0.8pc/|{\tensor}[ll]_{v^{\flat}}
}}
\right)
:=
\left(
\hbox{\small
\xymatrix@C-1.5pc{
\E_{S}
\ar@/^0.8pc/|{\tensor}[rr]^{u^{\flat}\comp f^{-1}}
& &
\E_1 
\ar[ll]^{f\comp u^{\sharp}}
& &
\dots
& &
\F_1
\ar[rr]_{g\comp v^{\sharp}}
& &
\F_{S}
\ar@/_0.8pc/|{\tensor}[ll]_{v^{\flat}\comp g^{-1}}
}}
\right)
$$

\begin{corollary}
\label{orbits}
The orbits of the $\Gl(\E)\times_T\Gl(\F)$-operation on
$\KGl(\E,\F)$ are the locally closed subvarieties
$$
O_{I,J}(\E,\F):=
\overline{O}_{I,J}(\E,\F)
\setminus
\left(
\bigcup_{i\not{\in}\I}Z_i(\E,\F)
\cup
\bigcup_{j\not{\in}\J}Y_j(\E,\F)
\right)
\quad,
$$
where $I,J\subseteq[0,n-1]$ with $\min I + \min J \geq n$,
and where 
$Z_i(\E,\F):=\overline{O}_{\{i\},\emptyset}(\E,\F)$ and 
$Y_j(\E,\F):=\overline{O}_{\emptyset,\{j\}}(\E,\F)$.
\end{corollary}

\begin{proof}
The $S$-valued points of $O_{I,J}(\E,\F)$
are the generalized isomorphisms
{\footnotesize
$$
\Phi = 
\left(
\xymatrix@C-4ex{
\E_{S}
\ar@/^1.2pc/|{\tensor}[rr]
\ar @/_0.65pc/ @{} [rr]|{(\M_0,\mu_0)}
& &
\E_1 
\ar[ll]_0
\ar @/_0.65pc/ @{} [rr]|{(\M_1,\mu_1)}
\ar@/^1.2pc/|{\tensor}[rr]
& &
\E_2
\ar[ll]_1
& 
\dots
& 
\E_{n-1}
\ar @/_0.65pc/ @{} [rr]|{(\M_{n-1},\mu_{n-1})}
\ar@/^1.2pc/|{\tensor}[rr]
& & 
\E_n
\ar[ll]_{n-1}
\ar[rr]^\sim
& & 
\F_n
\ar[rr]^{n-1}
& & 
\F_{n-1}
\ar@/_1.2pc/|{\tensor}[ll]
\ar @/^0.65pc/ @{} [ll]|{(\Ll_{n-1},\lambda_{n-1})}
& 
\dots
&
\F_2
\ar[rr]^1
& &
\F_1
\ar[rr]^0
\ar @/^0.65pc/ @{} [ll]|{(\Ll_1,\lambda_1)}
\ar@/_1.2pc/|{\tensor}[ll]
& &
\F_{S}
\ar @/^0.65pc/ @{} [ll]|{(\Ll_0,\lambda_0)}
\ar@/_1.2pc/|{\tensor}[ll]
}
\right)
$$
}
where $\mu_i=\lambda_j=0$ for $i\in I$ and $j\in J$ and
where $\mu_i$, $\lambda_j$ are nowhere vanishing for $i\not\in I$,
$j\not\in J$. 
It is clear that $O_{I,J}(\E,\F)$ is invariant under
the operation of $\Gl(\E)\times_T\Gl(\F)$.
From the proof of theorem \ref{strata} it follows that we have the
following isomorphism
$$
O_{I,J}(\E,\F)
\isomorph
\overset{o}{P_1}
\underset{\Fl}{\times}
\dots
\underset{\Fl}{\times}
\overset{o}{P_r}
\underset{\Fl}{\times}
\overset{o}{Q_s}
\underset{\Fl}{\times}
\dots
\underset{\Fl}{\times}
\overset{o}{Q_1}
\underset{\Fl}{\times}
\overset{o}{K'}
\quad,
$$
where 
\begin{eqnarray*}
\overset{o}{P_p} &:=& \PGl(V_{r-p+1}/V_{r-p},U_{s+p+1}/U_{s+p})
\quad (1\leq p\leq r)\\
\overset{o}{Q_q} &:=& \PGl(U_{s-q+1}/U_{s-q},V_{r+q+1}/V_{r+q})
\quad (1\leq q\leq s)\\
\overset{o}{K'} &:=& \text{Isom}(U_{s+1}/U_{s},V_{r+1}/V_{r})
\quad.
\end{eqnarray*}
There is a left 
$\Gl(\E)\times_T\Gl(\F)$-operation on
the right-hand side of this isomorphism, given on $S$-valued points
by
\begin{eqnarray*}
& &
(f,g)
(
(F_{\bullet}\E,F_{\bullet}\F),
\varphi_1,\dots,\varphi_r,
\psi_s,\dots,\psi_1,
\Phi'
)
:=  \\
& &\ 
(
(f(F_{\bullet}\E),g(F_{\bullet}\F)),
f^{-1}\varphi_1g,\dots,f^{-1}\varphi_rg,
g\psi_sf^{-1},\dots,g\psi_1f^{-1},
g\Phi'f^{-1}
) ,
\end{eqnarray*}
where $\varphi_p$ is an isomorphism 
(up to multiplication by an invertible section of $\Oo_S$) 
from $F_{r-p+1}\F/F_{r-p}\F$ to $F_{s+p+1}\E/F_{s+p}\E$
for $1\leq p\leq r$,\ 
$\psi_q$ an isomorphism 
(up to multiplication by an invertible section of $\Oo_S$) 
from $F_{s-q+1}\E/F_{s-q}\E$ to $F_{r+q+1}\F/F_{r+q}\F$ 
for $s\geq q\geq 1$
and
$\Phi'$ is an isomorphism from 
$F_{s+1}\E/F_{s}\E$ to $F_{r+1}\F/F_r\F$.
It is easy to see that this operation is transitiv and that 
the isomorphism 
$$
O_{I,J}(\E,\F)
\isomorph
\overset{o}{P_1}
\underset{\Fl}{\times}
\dots
\underset{\Fl}{\times}
\overset{o}{P_r}
\underset{\Fl}{\times}
\overset{o}{Q_s}
\underset{\Fl}{\times}
\dots
\underset{\Fl}{\times}
\overset{o}{Q_1}
\underset{\Fl}{\times}
\overset{o}{K'}
\quad,
$$
is $\Gl(\E)\times_T\Gl(\F)$-equivariant.
\end{proof}

\section{A morphism of $\KGln$ onto the Grassmannian compactification
         of the general linear group}
\label{grassmann}

Let $V$ be an $n$-dimensional vector space over some field.
As mentioned in the introduction,
there is another natural compactification of the general linear
group $\Gl(V)$: The Grassmannian $\Gr_n(V\oplus V)$ of 
$n$-dimensional subspaces of a $V\oplus V$-dimensional vector space. 
The embedding $\Gl(V)\injto\Gr_n(V\oplus V)$
is given by associating to an automorphism $V\isomto V$ its graph
in $V\oplus V$. We will see in this section that there exists a
natural morphism from $\KGl(V)$ to $\Gr_n(V\oplus V)$.
Our motivation here is to obtain a better understanding of the relation
between the Gieseker-type degeneration of moduli spaces of vector
bundles and the torsion-free sheaves approach as developed in 
\cite{Nagaraj&Seshadri} and \cite{Seshadri2}.

As in the previous section, we work over an arbitrary base scheme $T$.
Let $\E$, $\F$ be two locally free $\Oo_T$-modules of rank $n$.
Denote by $\Gr_n(\E\oplus \F)$ the Grassmanian variety over $T$ which
parametrizes  subbundles of rank $n$ of $\E\oplus \F$.
Let $S$ be a $T$-scheme and let
{\footnotesize
$$
\Phi = 
\left(
\xymatrix@C-4ex{
\E_{S}
\ar@/^1.2pc/|{\tensor}[rr]
\ar @/_0.65pc/ @{} [rr]|{(\M_0,\mu_0)}
& &
\E_1 
\ar[ll]_0
\ar @/_0.65pc/ @{} [rr]|{(\M_1,\mu_1)}
\ar@/^1.2pc/|{\tensor}[rr]
& &
\E_2
\ar[ll]_1
& 
\dots
& 
\E_{n-1}
\ar @/_0.65pc/ @{} [rr]|{(\M_{n-1},\mu_{n-1})}
\ar@/^1.2pc/|{\tensor}[rr]
& & 
\E_n
\ar[ll]_{n-1}
\ar[rr]^\sim
& & 
\F_n
\ar[rr]^{n-1}
& & 
\F_{n-1}
\ar@/_1.2pc/|{\tensor}[ll]
\ar @/^0.65pc/ @{} [ll]|{(\Ll_{n-1},\lambda_{n-1})}
& 
\dots
&
\F_2
\ar[rr]^1
& &
\F_1
\ar[rr]^0
\ar @/^0.65pc/ @{} [ll]|{(\Ll_1,\lambda_1)}
\ar@/_1.2pc/|{\tensor}[ll]
& &
\F_{S}
\ar @/^0.65pc/ @{} [ll]|{(\Ll_0,\lambda_0)}
\ar@/_1.2pc/|{\tensor}[ll]
}
\right)
$$
}
be a generalized isomorphism from $\E_S$ to $\F_S$.
By \ref{generalized isomorphism}.2, the
morphism $\E_n\to\E_S\oplus\F_S$ induced by the two
composed morphisms
\begin{eqnarray*}
&
\E_n\to\E_{n-1}\to\dots\to\E_1\to\E_S
& \\
&
\E_n\isomto\F_n\to\F_{n-1}\to\dots\to\F_1\to\F_S
\end{eqnarray*} 
is a subbundle of $\E_S\oplus\F_S$.
Let 
$$
\KGl(\E,\F)\to\Gr_n(\E\oplus\F)
$$ 
be the morphism,
which on $S$-valued points is given by 
$\Phi\mapsto (\E_n\to\E_S\oplus\F_S)$. 
Observe that the following diagram commutes
$$
\xymatrix{
&
\Isom(\E,\F)
\ar@{^{(}->}[dl]
\ar@{_{(}->}[dr]
&
\\
\KGl(\E,\F)
\ar[rr]
& &
\Gr_n(\E\oplus\F)
}
$$
and that furthermore all the arrows in this diagram are
equivariant with respect to the natural action of
$\Gl(\E)\times_T\Gl(\F)$ on the three schemes.
In the next proposition we compute the fibres of
the morphism
$\KGl(\E,\F)\to\Gr_n(\E\oplus\F)$.

\begin{proposition}
Let $S'$ be a $T$-scheme and let $\Hh\injto\E_{S'}\oplus\F_{S'}$
be an $S'$-valued point of $\Gr_n(\E,\F)$ such that  
$\im(\Hh\to\E_{S'})$ and $\im(\Hh\to\F_{S'})$ are subbundles of
$\E_{S'}$ and $\F_{S'}$ respectively. Then the fibre product
$$
\KGl(\E,\F)\underset{\Gr_n(\E\oplus\F)}{\times}S'
$$
is isomorphic to
$$
\PGlb(\text{ker}(\Hh\to\E_{S'}),\coker(\Hh\to\E_{S'}))
\times_{S'}
\PGlb(\text{ker}(\Hh\to\F_{S'}),\coker(\Hh\to\F_{S'}))
,
$$
where by convention $\PGlb(\Nn,\Nn):=S'$ for the zero-sheaf 
$\Nn=0$ on $S'$.
\end{proposition}

\begin{proof}
Let $S$ be an $S'$-scheme.
An $S$-valued point of the fibre product
$
\KGl(\E,\F)\times_{\Gr_n(\E\oplus\F)}S'
$
is given by a generalized isomorphism 
{\footnotesize
$$
\Phi = 
\left(
\xymatrix@C-4ex{
\E_{S}
\ar@/^1.2pc/|{\tensor}[rr]
\ar @/_0.65pc/ @{} [rr]|{(\M_0,\mu_0)}
& &
\E_1 
\ar[ll]_0
\ar @/_0.65pc/ @{} [rr]|{(\M_1,\mu_1)}
\ar@/^1.2pc/|{\tensor}[rr]
& &
\E_2
\ar[ll]_1
& 
\dots
& 
\E_{n-1}
\ar @/_0.65pc/ @{} [rr]|{(\M_{n-1},\mu_{n-1})}
\ar@/^1.2pc/|{\tensor}[rr]
& & 
\E_n
\ar[ll]_{n-1}
\ar[rr]^\sim
& & 
\F_n
\ar[rr]^{n-1}
& & 
\F_{n-1}
\ar@/_1.2pc/|{\tensor}[ll]
\ar @/^0.65pc/ @{} [ll]|{(\Ll_{n-1},\lambda_{n-1})}
& 
\dots
&
\F_2
\ar[rr]^1
& &
\F_1
\ar[rr]^0
\ar @/^0.65pc/ @{} [ll]|{(\Ll_1,\lambda_1)}
\ar@/_1.2pc/|{\tensor}[ll]
& &
\F_{S}
\ar @/^0.65pc/ @{} [ll]|{(\Ll_0,\lambda_0)}
\ar@/_1.2pc/|{\tensor}[ll]
}
\right)
$$}%
from $\E_S$ to $\F_S$ such that the induced
morphism $\E_n\injto\E_S\oplus\F_S$ identifies
$\E_n$ with the subbundle $\Hh_S$.
Let $i_1$ and $j_1$ be the ranks of $\im(\Hh\to\E_{S'})$
and $\im(\Hh\to F_{S'})$ respectively.
Observe that $i_1+j_1\geq n$. We restrict ourselves to the case,
where $i_1$ and $j_1$ are both strictly smaller than $n$.
(The cases where one or both of $i_1,j_1$ are equal to $n$ 
are proved analogously).
Then the sections $\mu_0,\dots,\mu_{i_1-1}$ and 
$\lambda_{j_1-1},\dots,\lambda_0$ are invertible
and $\mu_{i_1}=\lambda_{j_1}=0$.
From the proof of theorem \ref{strata} it follows that
such a $\Phi$ may be given by a tupel
$
((F_{\bullet}\E,F_{\bullet}\F),\varphi,\psi,\Phi')
$
where
\begin{eqnarray*}
F_{\bullet}\E &=&
(0=F_0\E\subseteq F_1\E\subseteq F_2\E\subseteq F_3\E=\E_S)
\\
F_{\bullet}\F &=&
(0=F_0\F\subseteq F_1\F\subseteq F_2\F\subseteq F_3\F=\F_S)
\end{eqnarray*}
are the filtrations given by
\begin{eqnarray*}
F_1\E &:=&
\im(\ker(\Hh_S\to\F_S)\to\E_S)
\\
F_2\E &:=&
\im(\Hh_S\to\E_S)
\\
F_1\F &:=&
\im(\ker(\Hh_S\to\E_S)\to\F_S)
\\
F_2\F &:=&
\im(\Hh_S\to\F_S)
\quad,
\end{eqnarray*}
$\varphi$ is a complete collineation from 
$F_1\F/F_0\F\isomorph\ker(\Hh_S\to\E_S)$ 
to  $F_3\E/F_2\E\isomorph\coker(\Hh_S\to\E_S)$, 
$\psi$ is a complete collineation from
$F_1\E/F_0\E\isomorph\ker(\Hh_S\to\F_S)$ to  
$F_3\F/F_2\F\isomorph\coker(\Hh_S\to\F_S)$,
and $\Phi'$ is the  isomorphism
$$
F_2\E/F_1\E
\isomto
\Hh_S/(\ker(\Hh_S\to\E_S)+\ker(\Hh_S\to\F_S))
\isomto
F_2\F/F_1\F
\quad.
$$
We see in particular that the tupel
$
((F_{\bullet}\E,F_{\bullet}\F),\varphi,\psi,\Phi')
$
is already determined by the subbundle 
$\Hh_S\injto\E_S\oplus\F_S$ (i.e. by the morphism $S\to S'$)
and the pair $(\varphi,\psi)$.
\end{proof}


\end{document}